\documentclass[11pt,a4paper,openbib]{article}
\usepackage[ansinew]{inputenc}
\usepackage{amssymb}
\usepackage{amsmath}
\usepackage{amsxtra} 
\usepackage{graphicx}
{\tiny }
\newcommand {\rel} {{\mathbb R}}
\newcommand {\com} {{\mathbb C}}

\newcommand {\Mod} {{\mathcal M}}
\newcommand {\Dom} {{\mathcal D}}

\newcommand {\nat} {{\mathbb N}}

\newcommand {\ganz} {{\mathbb Z}}
\newcommand {\quat} {{\mathbb H}}
\newcommand {\sphere} {{\mathbb S}}
\newcommand {\Area} {{\mathcal{A} }}
\newcommand {\Will} {{\mathcal{W} }}
\newcommand {\Wil} {{\mathcal{E} }}
\newcommand {\Stereo} {{\mathcal{P} }}

\newcommand {\CC}   {{\mathcal{C} }} 
\newcommand {\PP} {{\mathcal{P} }}

\newcommand {\Mill} {{\mathcal{M} }} 

\newcommand {\LL} {{\mathcal{L} }}

\def    \Hn     {{ {\cal H}^2 }}
\def    \Hnm    {{ {\cal H}^{1} }}

\def    \Lno    {{ {\cal L}^2 }}

\begin{document}

\newtheorem{theorem}{Theorem}[section]
\newtheorem{definition}{Definition}[section]
\newtheorem{proposition}{Proposition}[section]
\newtheorem{lemma}{Lemma}[section]
\newtheorem{corollary}{Corollary}[section]
\newtheorem{remark}{Remark}[section]
\newtheorem{example}{Example}

\author{Ruben Jakob}

\title{Singularities and full convergence of the M\"obius-invariant Willmore flow in the $3$-sphere}

\maketitle

\begin{abstract}
\noindent 	
In this article we continue our investigation 
of the M\"obius-invariant Willmore flow (MIWF), 
starting to move in arbitrary $C^{\infty}$-smooth and umbilic-free initial immersions $F_0$ which map some compact torus $\Sigma$ into $\rel^n$ respectively $\sphere^n$. Here we investigate the behaviour of flow lines $\{F_t\}$ of the MIWF in $\sphere^3$ starting with relatively low Willmore energy, as the time $t$ approaches the maximal time of existence $T_{\textnormal{max}}(F_0)$ of $\{F_t\}$. 
We particularly construct divergent flow lines and we investigate the formation of ``limit surfaces'' of both divergent and convergent flow lines of the MIWF. 
Such a limit surface is the support of a certain integral $2$-varifold $\mu$ in $\rel^4$, arising as a measure-theoretic limit of the sequence of varifolds $\{\Hn\lfloor_{F_{t_{l}}(\Sigma)}\}$, 
for an appropriately chosen sequence  
$t_{l} \nearrow  T_{\textnormal{max}}(F_0)$. 
The support $\textnormal{spt}(\mu)$ of such a limit 
$2$-varifold $\mu$ is either empty or homeomorphic to some compact closed manifold of genus either $0$ or $1$. In the ``non-degenerate case'', in which the limit surface $\textnormal{spt}(\mu)$ is a compact surface of genus $1$, $\textnormal{spt}(\mu)$ can be parametrized by a uniformly conformal bi-Lipschitz homeomorphism $f$ of class 
$(W^{2,2}\cap W^{1,\infty})(\Sigma)$,  
and under certain additional conditions on the sequence $\{F_{t_{l}}\}$ such a uniformly conformal parametrization 
$f$ is a diffeomorphism of class $W^{4,2}(\Sigma,\rel^4)$. Finally, if the initial immersion $F_0$ of a flow line $\{F_t\}$ is assumed to parametrize a smooth Hopf-torus in $\sphere^3$ with Willmore energy not bigger than 
$8 \pi$, then we obtain more precise statements 
about the flow line $\{F_t\}$ as $t \nearrow T_{\textnormal{max}}(F_0)$, especially stronger types of convergence of particular subsequences of $\{F_t\}_{t \geq 0}$ to uniformly conformal $W^{4,2}$-parametrizations of limit Hopf-tori. This insight will finally yield a new criterion for full smooth convergence of such flow lines of the MIWF to smooth diffeomorphisms parametrizing the Clifford torus - up to M\"obius-transformations of $\sphere^3$ - as 
$t \nearrow \infty$.      
\end{abstract}


\section{Introduction and main results}

In this article, the author continues his investigation
of the M\"obius-invariant Willmore flow (MIWF) -
evolution equation \eqref{Moebius.flow} below - in $\sphere^3$ 
respectively in $\rel^3$, which is geometrically motivated  
by the prominent {\bf Willmore-functional} 
\begin{equation} \label{Willmore.functional}
\Will(F):= \int_{\Sigma} K^M_F + \frac{1}{4} \, \mid \vec H_F \mid^2 \, d\mu_{F^*(g_{\textnormal{euc}})}.
\end{equation}
This functional of fourth order can be considered on $C^{4}$-immersions $F:\Sigma \longrightarrow M$, from any smooth compact orientable surface $\Sigma$ 
into an arbitrary smooth Riemannian
manifold $M$, where $K^M_F(x)$ denotes the
sectional curvature of $M$ with respect to the ``immersed tangent plane'' $DF_x(T_x\Sigma)$ in $T_{F(x)}M$. 
In this article, there are only two cases relevant, 
namely the cases $M=\rel^n$ and $M=\sphere^n$, 
in which we simply have $K_F\equiv 0$ respectively 
$K_F\equiv 1$. For ease of exposition, we will only consider the case $M=\sphere^n$, for $n \geq 3$, in the sequel of this introduction. In the author's article \cite{Jakob_Moebius_2016}, the author proved short-time existence and uniqueness of the flow
\begin{equation}  \label{Moebius.flow}
\partial_t F_t = -\frac{1}{2} \frac{1}{|A^0_{F_t}|^4} \,
\Big{(} \triangle_{F_t}^{\perp} \vec H_{F_t} + Q(A^{0}_{F_t})
(\vec H_{F_t}) \Big{)}
\equiv -\frac{1}{|A^0_{F_t}|^4} \,\nabla_{L^2} \Will(F_t),
\end{equation}
which is well-defined on differentiable families of $C^4$-immersions $F_t$ mapping some arbitrarily fixed 
{\bf smooth compact torus} $\Sigma$ into either $\rel^n$ or $\sphere^n$, for $n\geq 3$, {\bf without any umbilic points}. 
As already pointed out in the author's article \cite{Jakob_Moebius_2016}, the ``umbilic free condition'' $|A^0_{F_t}|^2>0$ on $\Sigma$ implies
$\chi(\Sigma)=0$ for the Euler-characteristic of $\Sigma$,
which forces the flow \eqref{Moebius.flow} to be only 
well-defined on families of sufficiently smooth umbilic-free tori, being immersed into $\rel^n$ or $\sphere^n$.\\
Now, given some immersion $F:\Sigma \longrightarrow \sphere^n$, we endow the torus $\Sigma$ with the pullback $g_F:=F^*g_{\textnormal{euc}}$ of the Euclidean metric on $\sphere^n$, i.e. with coefficients 
$g_{ij}:=\langle \partial_i F, \partial_j F \rangle$, and we let $A_F$ denote the second fundamental form of the immersion 
$F$, defined on pairs of tangent vector fields $X,Y$ on $\Sigma$ by:
\begin{equation}  \label{second.fundam.form}
A_F(X,Y):= D_X(D_Y(F))-P^{\textnormal{Tan}(F)}(D_X(D_Y(F)))
\equiv (D_X(D_Y(F)))^{\perp_{F}},
\end{equation}
where $D_X(V)\lfloor_{x}$ denotes the projection of the usual 
derivative of a vector field $V:\Sigma \longrightarrow \rel^{n+1}$ in direction of the vector field $X$ into the respective fiber $T_{F(x)}\sphere^n$ of 
$T\sphere^n$, $P^{\textnormal{Tan}(F)}:\bigcup_{x \in \Sigma} \{x\} \times \textnormal{T}_{F(x)}\sphere^n \longrightarrow
\bigcup_{x \in \Sigma} \{x\} \times \textnormal{T}_{F(x)}(F(\Sigma))=:\textnormal{Tan}(F)$
denotes the bundle morphism which projects
the entire tangent space $\textnormal{T}_{F(x)}\sphere^n$ orthogonally onto its subspace $\textnormal{T}_{F(x)}(F(\Sigma))$, the tangent space of the immersion $F$ in $F(x)$, for every $x \in \Sigma$, and where $^{\perp_{F}}$ abbreviates the bundle morphism $\textnormal{Id}_{\textnormal{T}_{F(\cdot)}\sphere^n}
-P^{\textnormal{Tan}(F)}$. Furthermore, $A^{0}_F$ 
denotes the tracefree part of $A_F$, i.e.
$$
A^{0}_F(X,Y):= A_F(X,Y) - \frac{1}{2} \,g_F(X,Y)\, \vec H_F
$$
and $\vec H_F:=\textnormal{Trace}(A_F) \equiv A_F(e_i,e_i)$
(``Einstein's summation convention'') denotes the mean curvature vector of $F$, where $\{e_i\}$ is a local orthonormal frame of the tangent bundle $T\Sigma$. Finally, $Q(A_F)$ respectively $Q(A^0_F)$ operates on vector fields $\phi$ which are sections into the normal bundle of $F$, i.e. which are normal along $F$, by assigning $Q(A_F)(\phi):=A_F(e_i,e_j) \langle A_F(e_i,e_j),\phi \rangle$, which is by definition again a section into the normal bundle of $F$.
Weiner computed in his seminal paper \cite{Weiner} 
that the first variation of the Willmore functional 
$\nabla_{L^2} \Will$ in some smooth immersion $F$, 
in direction of a smooth section 
$\phi$ into the normal bundle of $F$, is in both cases $M=\sphere^n$ and $M=\rel^n$, $n\geq 3$, given by:
\begin{eqnarray} \label{first_variation}
\langle \nabla_{L^2} \Will(F),\phi \rangle_{L^2(\Sigma,\mu_{g_F})} 
\equiv \int_{\Sigma} \big{\langle} \nabla_{L^2} \Will(F),\phi \big{\rangle} \, d\mu_{g_F}    
= \frac{1}{2} \int_{\Sigma}
\big{\langle} \triangle^{\perp}_{F} \vec H_{F} + Q(A^{0}_{F})
(\vec H_{F}), \phi \big{\rangle}\, d\mu_{g_F}.           
\end{eqnarray}
The decisive difference between the flow \eqref{Moebius.flow} and the classical $L^2$-Willmore-gradient-flow, i.e. the $L^2$-gradient-flow of functional 
\eqref{Willmore.functional}, is the
factor $\frac{1}{|A^0_{F_t}|^4(x)}$, which is finite in 
$x\in \Sigma$, if and only if $x$ is not a umbilic point of the immersion $F_t$. It is this additional factor which 
on the one hand makes the analytic investigation of 
flow \eqref{Moebius.flow} significantly more difficult, 
but on the other hand turns the classical Willmore gradient 
flow into {\bf a conformally invariant flow}, 
thus correcting the scaling behaviour of the 
classical Willmore flow in exactly the right way. 
More precisely, there holds the following important lemma 
for the case $M=\sphere^n$, $n\geq 3$, which is equivalent 
to Lemma 1 in \cite{Jakob_Moebius_2016} dealing with 
the case $M=\rel^n$:
\begin{lemma}  \label{Transformation.of.Willmore}
Let $\Sigma$ be a smooth torus, $\Phi$ an arbitrary M\"obius transformation of $\sphere^n$, $n\geq 3$, and $F:\Sigma \longrightarrow \sphere^n$ a $C^{4}$-immersion satisfying $\mid A^0_{F} \mid^2>0$ on $\Sigma$.  
If we substitute $F$ by the composition $\Phi \circ F$,
then the differential operator 
$F\mapsto \mid A^0_F \mid^{-4} \,\nabla_{L^2} \Will(F)$
transforms like: 
\begin{equation}  \label{trafo_heat_expression}
\mid A^0_{\Phi(F)} \mid^{-4} \, \nabla_{L^2} \Will(\Phi(F))
=D\Phi(F) \cdot \Big{(} \mid A^0_{F} \mid^{-4} \, \nabla_{L^2} \Will(F)  \Big{)}   \qquad \textnormal{on}\,\, \Sigma.
\end{equation}
\end{lemma}
\noindent  
Since for the differential operator $\partial_t$ applied to  $C^1$-families $\{F_t\}$ of $C^{4}$-immersions the chain rule yields the same transformation formula as in \eqref{trafo_heat_expression}, i.e.
$\partial_t(\Phi(F_t))= D\Phi(F_t) \cdot \partial_t(F_t)$,
we can indeed derive the conformal invariance of flow
\eqref{Moebius.flow} from Lemma \ref{Transformation.of.Willmore},
just as in Corollary 1 in the author's article \cite{Jakob_Moebius_2016}:
\begin{corollary} \label{trafo_heat_expression_2}
Any family $\{F_t\}$ of $C^4$-immersions $F_t:\Sigma \longrightarrow \sphere^n$ satisfying 
$\mid A^0_{F_t} \mid^2>0$ on $\Sigma$\, $\forall \, t\in[0,T]$, solves flow equation \eqref{Moebius.flow}
$\forall \, t\in [0,T]$ if and only if its composition $\Phi(F_t)$ with an arbitrary M\"obius transformation $\Phi$ of $\sphere^n$ solves the same flow equation again, 
thus, if and only if
$$
\partial_t (\Phi(F_t))=- \mid A^0_{\Phi(F_t)} \mid^{-4} \, \nabla_{L^2} \Will(\Phi(F_t))
$$
holds $\forall \, t \in [0,T]$ and $\forall \, \Phi 
\in \textnormal{M\"ob}(\sphere^n)$.
\end{corollary}
\noindent 
It is this conformal invariance property of flow \eqref{Moebius.flow} which explains its name:
the {\bf M\"obius-invariant Willmore flow}, or shortly 
{\bf MIWF}. First of all, one might guess that this stark  
difference between the MIWF and the classical Willmore flow was extremely powerful in view of Theorem 4.2 in \cite{Kuwert.Schaetzle.2012} respectively in view of
Theorem 4.2 in \cite{Schaetzle.Conf.factor.2013}, 
which seem to guarantee us here, that the induced metrics 
$g_{\textnormal{euc}}\lfloor_{F_{t_j}(\Sigma)}$ along any 
fixed flow line $\{F_t\}$ of the MIWF would be conformally 
equivalent to smooth metrics $g_{\textnormal{poin},j}$ 
of vanishing scalar curvature, such that the resulting 
conformal factors $u_{t_j}$ can be uniformly estimated in 
$L^{\infty}(\Sigma)$, for some sequence of times 
$t_j \nearrow T_{\textnormal{max}}(F_0)$. 
But since the conformal invariance of the MIWF lets 
us apply {\bf only finitely many} 
conformal transformations to a fixed flow line 
$\{F_t\}$ of the MIWF - in order to either estimate its 
lifespan or to determine its behaviour 
as $t\nearrow T_{\textnormal{max}}(F_0)$ -  
neither Theorem 4.2 in \cite{Kuwert.Schaetzle.2012}, nor Proposition 2.2 in \cite{Schaetzle.Conf.factor.2013}, nor     
Theorems 3.2 and 4.2 in \cite{Schaetzle.Conf.factor.2013}, 
nor Theorem 4.1 in \cite{Kuwert.Li.2012}   
can be applied here, in order to obtain any 
valuable information about the limiting behaviour 
of a fixed flow line $\{F_t\}$ of the MIWF, as 
$t \nearrow T_{\textnormal{max}}(F_0)$. 
The only accessible information in this general 
situation appears to be given by Theorem 5.2 in 
\cite{Kuwert.Li.2012} respectively by Theorem 1.1 in 
\cite{Riviere.Lip.conf.imm.2013} telling us, that the complex structures $S(F_t(\Sigma))$ corresponding to the conformal classes of the induced metrics
$g_{\textnormal{euc}}\lfloor_{F_t(\Sigma)}$ are contained
in some compact subset of the moduli space $\Mill_1$, i.e. 
cannot diverge to the boundary of 
$\Mill_1 \cong \quat/\textnormal{PSL}_2(\ganz)$; 
see here also Theorem 5.1 in \cite{Kuwert.Schaetzle.2012} 
for the earliest reference concerning this type of 
result. But this information does not suffice, neither in order to exclude ``loss of topology'' nor ``degeneration to a constant map'' for an arbitrarily chosen sequence of immersions $\{F_{t_j}\}$ belonging to some fixed flow line $\{F_t\}$ of the MIWF, as $t_j \nearrow T_{\textnormal{max}}(F_0)$; see here Section 5.1 in \cite{Riviere.Park.City.2013} for some illustrative examples and explanations, and compare also to  
the statements of Theorems 1.4, 3.2 and 4.2 in \cite{Schaetzle.Conf.factor.2013} and of Theorem 5.1 
in \cite{Kuwert.Li.2012}.  \\        
Moreover, Proposition 4 and Theorem 4 in the Appendix of 
\cite{Jakob_Moebius_2016} show us that 
the geometrically motivated adjustment of the classical 
Willmore flow can only be achieved by means of a particular 
change of its scaling behaviour. One can easily compute, 
that for any smooth immersion $F:\Sigma \longrightarrow \rel^n$ and for any $\rho>0$ there holds:  
\begin{eqnarray*}      
	A_{\rho \,F} = \rho\, A_{F}, \quad    
	A^0_{\rho \,F} = \rho\, A^0_{F},   \quad            
   	\vec H_{\rho \,F} = \rho^{-1}\,\vec H_{F}  \quad \textnormal{on} \,\, \Sigma,          \\                  
	\triangle_{\rho F}(\vec H_{\rho F}) =
	\rho^{-3} \,\triangle_{F}(\vec H_{F}) 
    \quad  \textnormal{and}  \quad 
	Q(A^{0}_{\rho F}).\vec H_{\rho F}= \rho^{-3}\,
	Q(A^{0}_{F}).\vec H_{F}   \quad \textnormal{on} \,\, \Sigma,   
\end{eqnarray*}
which implies that for any flow line $\{F_t\}_{t\geq 0}$ 
of the classical Willmore flow in $\rel^n$ and for any $\rho>0$ the family $\{\rho\, F_{\rho^{-4}t}\}_{t\geq 0}$ is again a flow line of the classical Willmore flow in $\rel^n$, whereas any flow line $\{F_t\}$ of the MIWF in $\rel^n$ can be scaled with any factor $\rho>0$ in the ambient space $\rel^n$, without losing its property of being a flow line of the MIWF, but no time-gauge is necessary here ! 
In particular, we are not able here, to adopt 
the ``blow-up construction'' of Section 4 in 
\cite{Kuwert_Schaetzle_2001} and its powerful combination 
with the lower bound on the lifespan of any flow line
of the classical Willmore flow from Theorem 1.2 of \cite{Kuwert_Schaetzle_2002}, leading to the 
first satisfactory \underline{full convergence result} 
for the classical Willmore flow in Theorem 5.1 of \cite{Kuwert_Schaetzle_2001} and later to its optimization 
in Theorem 5.2 of \cite{Kuwert.Schaetzle.Annals}. 
Indeed, only a few years after the publication of the paper 
\cite{Kuwert.Schaetzle.Annals} Blatt constructed in \cite{Blatt.2009} some concrete example of a singular flow line of the classical Willmore flow moving \emph{rotationally symmetric} surfaces in $\rel^3$ of genus $0$ whose initial Willmore energy is only slightly bigger than $8\pi$, thus proving optimality of Theorem 5.2 in \cite{Kuwert.Schaetzle.Annals}. At this point one should also mention the more recent contribution   \cite{Dall.Acqua.Schaetzle.Mueller.2024}, in which the 
authors prove that the number $8\pi$ is again 
\emph{the optimal energy threshold} below which any rotationally symmetric surface in $\rel^3$ of genus $1$ would certainly initiate a global flow line converging fully in every $C^k$-norm to the Clifford torus in $\rel^3$, possibly rescaled and translated.  
Below in Theorem \ref{Singularities.do.exist} we 
will supply a modest counterpart  
of Lemma 3.8 in \cite{Dall.Acqua.Schaetzle.Mueller.2024} respectively of Theorems 4.1 and 5.1 in \cite{Blatt.2009}, 
detecting some particular \emph{divergent 
flow lines} of the degenerate elastic energy flow \eqref{elastic.energy.flow} respectively of the MIWF; 
see here also Definition \ref{maximal.interval} (c) 
and Remark \ref{divergent.flow} (1) below. 
Although the \emph{existence of divergent flow lines} 
of the MIWF was strongly expected - in view of its 
apparently singular evolution equation \eqref{Moebius.flow} -  
their concrete detection respectively construction 
turned out to be rather challenging. \\ 
The first mathematical indication for the 
\emph{existence of divergent flow lines} of the MIWF 
was the tedious discovery that any attempt to obtain 
some estimate on the lifespan of a general flow line
of the MIWF, following the lines of the fundamental paper \cite{Kuwert_Schaetzle_2002} or of its   
adaption to the \emph{inverse Willmore flow} in \cite{Mayer},
would end up in some sort of ``computational chaos''. 
More precisely, so far any estimate on the lifespan 
of some $2$nd or $4$th order geometric flow - only depending on geometric data at time $t=0$ - follows from a criterion for the singular time $T_{\textnormal{max}}$ in terms of ``blow up'' of appropriately chosen geometric data along an arbitrarily chosen flow line as $t\nearrow T_{\textnormal{max}}$, and the technical tool behind 
such a characterization is a suitable substitute of 
\emph{Bernstein-Bando-Shi-estimates} 
\footnote{Bernstein-Bando-Shi-estimates were 
developed in the investigation 
of long-time behaviour of Hamilton's Ricci flow on 
compact, closed manifolds.} 
estimating covariant derivatives of any order of 
certain particular tensor fields in $L^{\infty}$ both with respect to space- and time-variables, and this can only be done by induction over the order of covariant differentiation; 
compare for example with Theorem 8.1 in Chapter 8 of \cite{Andrews.Hopper.2011} regarding the Ricci flow, 
with Section 8 of \cite{Huisken.1984} examining the mean-curvature flow, or with Theorem 3.5 in 
\cite{Kuwert_Schaetzle_2001} respectively Section 4 in 
\cite{Kuwert_Schaetzle_2002} dealing with the Willmore flow. But the factor 
$|A^0_{F_t}|^{-4}$ in front of the $L^2$-gradient 
of the Willmore energy in evolution equation \eqref{Moebius.flow} ``seems to produce too many covariant derivatives of $|A^0_{F_t}|^2$'' in each step of such an attempted induction.
In order to investigate and thoroughly explain this 
unpleasant phenomenon of the MIWF, one should  
not struggle with evolution equation \eqref{Moebius.flow} itself, which would be to exchange the classical or inverse Willmore flow equation with equation \eqref{Moebius.flow} and then to try to adapt the procedures in \cite{Kuwert_Schaetzle_2002} respectively \cite{Mayer}. Instead it turned out to be much easier to focus on the 
restriction of the MIWF to Hopf-tori in $\sphere^3$ 
and to reduce those special flow lines of the MIWF 
- by means of the Hopf-fibration 
$\pi:\sphere^3\longrightarrow \sphere^2$ - 
to flow lines of the degenerate version \eqref{elastic.energy.flow} of the classical elastic energy flow \eqref{classical.elastic.energy.flow} moving closed 
smooth curves in $\sphere^2$.   
The latter geometric flow was thoroughly investigated in \cite{Dall.Acqua.Pozzi.2014} adapting classical methods from \cite{Dall.Acqua.Pozzi.2018}, where covariant derivatives 
$(\nabla^{\perp}_s)^k(\vec \kappa_{\gamma_t})$ of any order 
$k$ of the curvature vector $\vec \kappa_{\gamma_t}$ 
along some arbitrarily chosen flow line $\{\gamma_t\}$ 
of the elastic energy flow were uniformly estimated - 
similarly to the inductive procedure as mentioned above - 
i.e. combining the flow equation governing the behaviour of
$\{(\nabla^{\perp}_s)^k(\vec \kappa_{\gamma_t})\}$ with 
certain interpolation- and absorption techniques and 
with Gronwall's Lemma.       
Now, below in Remark \ref{no.classical.approach}  
we will quickly infer from the mentioned Theorem \ref{Singularities.do.exist} the decisive reason 
``why indeed the standard technique from \cite{Dall.Acqua.Pozzi.2014} and \cite{Dall.Acqua.Pozzi.2018} \underline{cannot be successfully applied to flow \eqref{elastic.energy.flow}}'', although at first sight flow \eqref{elastic.energy.flow} appears to be a uniformly parabolic and subcritical flow of fourth order, just as its classical counterpart \eqref{classical.elastic.energy.flow}.  
This striking fact about the \emph{degenerate} elastic energy flow \eqref{elastic.energy.flow} respectively about the MIWF 
already indicates that we \emph{will not be able 
to identify familiar geometric criteria} 
for flow lines of the MIWF to be \emph{divergent in 
finite or infinite time}.\\
Furthermore, the \emph{reduced MIWF} starting in Hopf-tori in $\sphere^3$ does not only considerably simplify the construction of divergent flow lines of the MIWF - see 
Theorem \ref{Singularities.do.exist} -   
but we can also prove below in Theorem \ref{singular.time.MIWF.Hopf.tori} (i) rather elementarily, 
that along any chosen flow line $\{F_{t}\}$ 
of the \emph{reduced MIWF} there is at least some 
subsequence of any chosen sequence 
$t_j \nearrow T_{\textnormal{max}}(F_0)$, 
such that the embedded Hopf-tori $F_{t_{j_l}}(\Sigma)$ 
converge in Hausdorff-distance to some 
embedded Hopf-torus in $\sphere^3$. 
Hence, the genus along any such \emph{reduced  
flow line} $\{F_{t}\}$ is actually preserved in the limit, 
provided one focuses here on appropriate subsequences  
$t_{j_l} \nearrow T_{\textnormal{max}}(F_0)$.  
This important fact will enable us to directly apply 
here the result of Theorem \ref{limit.MIWF} (ii) below, 
estimating the conformal factors of the induced metrics 
$g_{\textnormal{euc}}\lfloor_{F_{t_{j_l}}(\Sigma)}$ of 
the considered Hopf-tori $F_{t_{j_l}}(\Sigma)$ 
without any further condition on the given flow line $\{F_{t}\}$.   \\
In accordance with the qualitative difficulties described 
above and with the additional open question ``how to bound 
$|A^0_{F_t}|^2$ away from zero for general flow lines 
$\{F_t\}$'' as $t\nearrow T_{\textnormal{max}}(F_0)$, 
the first two main theorems of this article, 
Theorems \ref{limit.MIWF} and \ref{singular.time.MIWF.Hopf.tori}, 
\underline{do not support} the expectation that we might be able to rule out \underline{singularities} of flow lines of the MIWF in terms of a ``no curvature concentration-condition'' along their flowing immersions, 
as $t$ approaches the respective maximal time of existence - see here statements \eqref{no.A.concentration} and \eqref{no.concentration.Mill} and Remarks 
\ref{no.better} and \ref{no.better.2} in the appendix - indicating a stark contrast to the behaviour 
of the classical Willmore flow in any $\rel^n$, $n \geq 3$, 
on account of Theorem 1.2 in \cite{Kuwert_Schaetzle_2002}.  \\
Moreover, those technical challenges and 
the existence of divergent flow lines of the MIWF 
- by Theorem \ref{Singularities.do.exist} - force us 
to impose rather strong a-priori conditions 
on flow lines of the MIWF in order to prove their full 
and smooth convergence to certain Willmore tori in $\sphere^3$.   
For example, in Theorem \ref{limit.at.infinity} 
we will assume \emph{global existence} 
and \emph{uniformly bounded mean curvature} 
$\vec H_{F_t,\sphere^3}$ in $L^{\infty}(\Sigma)$ of a flow line $\{F_t\}$ of the MIWF which starts moving in a smooth parametrization of some Hopf-torus with Willmore energy 
below $4 \pi^2$. 
Only in this rather restricted framework we will be 
able to make reasonably precise predictions on the 
``destiny'' of $\{F_t\}$ as $t \nearrow \infty$. 
Here we can combine the full strength of the 
reduction of the MIWF to the subcritical flow 
in \eqref{elastic.energy.flow} via the Hopf-fibration, 
the classification and analysis of elastic curves 
on $\sphere^2$ in Proposition 6 of \cite{Ruben.MIWF.II}, 
the results of Theorem \ref{singular.time.MIWF.Hopf.tori} below and the first \emph{full convergence result} in \cite{Ruben.MIWF.IV}, i.e. Theorem 1 in \cite{Ruben.MIWF.IV}, in order to prove that any such global flow line $\{F_t\}$ 
of the MIWF either ``diverges'' as $t \nearrow \infty$
while its energy descends to some Willmore energy between 
$8 \pi$ and $4 \pi^2$ - or fully converges 
- up to smooth reparametrization -  
in each $C^m(\Sigma,\rel^4)$-norm to a smooth and diffeomorphic parametrization of the Clifford torus 
- up to some conformal transformation of $\sphere^3$.    \\ 
In Theorem \ref{limit.MIWF} we will firstly aim at a rough understanding of ``how either singular or global flow lines of the MIWF behave as $t\nearrow T_{\textnormal{max}}(F_0)$'', i.e. how \underline{singularities} can only look like, 
in either finite or infinite time, and how topological  
and regularity properties of such singularities might depend on certain geometric or analytic quantities along 
the considered flow line $\{F_t\}$, provided 
$\{F_t\}$ starts moving in a umbilic-free immersion $F_0$ 
with Willmore energy not bigger than $8\pi$. 
We should stress here the fact that our techniques of examination will only produce statements about certain subsequences $\{F_{t_{j_l}}\}$ - and their 
\underline{limit surfaces} -  
of arbitrarily chosen sequences $\{F_{t_j}\}$ 
along some flow line $\{F_t\}$, 
both in Theorems \ref{limit.MIWF} and  \ref{singular.time.MIWF.Hopf.tori}; 
see here also Remarks \ref{no.better.2} and  
\ref{no.better.3} in the appendix.    \\
The proof of the first three parts of 
Theorem \ref{limit.MIWF} is based  
on a combination of Kuwert's and Sch\"atzle's \cite{Kuwert.Schaetzle.2012}, \cite{Kuwert.Schaetzle.conf.class.2013}, 
\cite{Schaetzle.Conf.factor.2013} and Rivi\`ere's \cite{Riviere.Lip.conf.imm.2013}, \cite{Riviere.Var.principle.2014}
investigation of sequences of immersions of a compact 
Riemann surface $\Sigma$ into some $\rel^n$ of fixed 
genus $p\geq 1$ which either have sufficiently small 
Willmore energy and whose conformal classes cannot approach the boundary of the moduli space $\Mod_p$ or 
which minimize the Willmore energy under fixed conformal class. In the proofs of the fourth part of Theorem \ref{limit.MIWF} and of Theorem \ref{singular.time.MIWF.Hopf.tori} we will complement 
several techniques and results of Theorem   
\ref{limit.MIWF} (1)--(3) with Rivi\`ere's \cite{Riviere.2008}, \cite{Riviere.2011}, \cite{Riviere.Park.City.2013}
and Bernard's \cite{Bernard.2016} discovery 
of certain conservation laws induced by the conformal 
invariance of the Willmore functional and with many 
of Palmurella's and Rivi\`ere's tricks in \cite{Palmurella.2022}, \cite{Palmurella.2024}, where
a new theory of \underline{weak flow lines} 
of the classical Willmore flow has been developed. 
All mentioned papers by Kuwert, Sch\"atzle, Bernard, 
Rivi\`ere and Palmurella share one key-aspect, 
namely a rather novel combination of 
\underline{Gauge theory} with 
suitably adapted versions of 
\underline{Wente's \cite{Wente.1969} respectively Brezis'
\cite{Brezis.Coron.1984} 
$(L^{\infty} \cap W^{1,2})$-estimates} - following the pioneering work by M\"uller, Sver\'ak \cite{Mueller.Sverak.1995} and H\'elein \cite{Helein.2004} - 
and their explicit motivation was to parametrize 
certain families of embedded surfaces in $\rel^n$ 
- up to M\"obius transformations - 
by means of uniformly conformal immersions, which has 
far reaching applications in the context of 
minimization of the Willmore functional under 
some prescribed geometric constraint 
\footnote{Here we should concretely think of fixing the conformal class \cite{Kuwert.Schaetzle.conf.class.2013} or 
the isoperimetric ratio \cite{Keller.Modino.Riviere}.}
or in the context of a distributional formulation 
of the Willmore flow \cite{Palmurella.2022}, \cite{Palmurella.2024}. \\ 
Let's have a look at our first result, 
Theorem \ref{limit.MIWF}, and let's find out
whether \emph{limit surfaces} of both singular and global flow lines of the MIWF do always exist - in some appropriate concrete sense - and what topological type respectively regularity they must have. 
\noindent 
\begin{theorem}  \label{limit.MIWF}
	Let $\Sigma$ be an arbitrary smooth compact torus, 
	and let $\{F_t\}$ be some flow line of the 
	MIWF starting in a smooth and umbilic-free immersion 
	$F_0:\Sigma \longrightarrow \sphere^3$ with 
	$\Will(F_0)\leq 8\pi$, and let $t_j \nearrow 
	T_{\textnormal{max}}(F_0)$ be arbitrarily chosen. 
	\begin{itemize}
	\item[1)] There is a subsequence $\{F_{t_{j_l}}\}$ and 
	some integral, $2$-rectifiable varifold $\mu$ with 
	unit Hausdorff-$2$-density such that 
	\begin{equation}  \label{weak.convergence.mu} 
	\Hn\lfloor_{F_{t_{j_l}}(\Sigma)} \longrightarrow \mu 
	\quad \textnormal{weakly as Radon measures on}\,\,\,\rel^4.   
	\end{equation} 
    If $\mu$ is non-trivial, then its non-empty support is a closed, embedded and orientable Lipschitz-surface in $\sphere^3$ of genus either $0$ or $1$, and moreover we have in this case: 
    \begin{equation} \label{Hausdorff.converg} 
	F_{t_{j_l}}(\Sigma) \longrightarrow \textnormal{spt}(\mu) 
	\quad \textnormal{as subsets of $\rel^4$ in 
	Hausdorff-distance}, 
	\end{equation}
	as $l \to \infty$. 
    \item[2)] In the ``non-degenerate case'' in which 
    there holds $\mu\not =0$ and $\textnormal{genus}(\textnormal{spt}(\mu))=1$ 
    for a limit varifold $\mu$ in \eqref{weak.convergence.mu}, 
    appropriate reparametrizations of possibly another  
    subsequence of the immersions $F_{t_{j_l}}$ from
    \eqref{weak.convergence.mu} converge 
    weakly in $W^{2,2}(\Sigma,\rel^4)$ and 
    weakly* in $W^{1,\infty}(\Sigma,\rel^4)$ to 
    a $(W^{2,2}\cap W^{1,\infty})$-parametrization $f$ of $\textnormal{spt}(\mu)$, which is a 
    bi-Lipschitz homeomorphism as well:
    \begin{equation}   \label{parametrization.mu}
    f:\Sigma \stackrel{\cong}\longrightarrow \textnormal{spt}(\mu) \subset \rel^4,
    \end{equation} 
    and $f$ is ``uniformly conformal with respect to $g_{\textnormal{poin}}$ on $\Sigma$'' in the sense that 
    $f^*g_{\textnormal{euc}}=e^{2u} \,g_{\textnormal{poin}}$ holds on $\Sigma$, for some smooth zero scalar curvature and unit volume metric $g_{\textnormal{poin}}$ on $\Sigma$ and for some real-valued function $u\in L^{\infty}(\Sigma)$ with 
    $\parallel u \parallel_{L^{\infty}(\Sigma)}\leq 
    \Lambda=\Lambda(\{F_{t_{j_l}}\},\mu)<\infty$. 
    Moreover, in this case we can reinterpret the integral varifold $\mu$ in \eqref{weak.convergence.mu} in two ways, namely there holds: 
    \begin{equation}  \label{spt.mu.f} 
    \mu_f = \mu = \Hn\lfloor_{\textnormal{spt}(\mu)} \quad \textnormal{on} \quad \rel^4,
    \end{equation}  
    where $\mu_f:=f(\mu_{f^*g_{\textnormal{euc}}})$ is the canonical surface measure of $f(\Sigma)=\textnormal{spt}(\mu)$ in $\rel^4$. 
    The coinciding varifolds $\mu$ and $\mu_f$ have weak mean curvature vectors $\vec H_{\mu}$, $\vec H_{\mu_f}$ 
    in $L^2(\Sigma,\mu)$, and they satisfy exactly: 
    \begin{equation}  \label{Willmore.energy}
    4\,\Will(\mu):=\int_{\rel^4} |\vec H_{\mu}|^2 \, d\mu 
    =\int_{\rel^4} |\vec H_{\mu_f}|^2 \, d\mu_f
    = \int_{\Sigma} |\vec H_{f,\rel^{4}}|^2 \, d\mu_{f^*{g_{\textnormal{euc}}}} \equiv 4\,\Will(f).  
    \end{equation}  
    \item[3)] If a limit varifold $\mu$ in \eqref{weak.convergence.mu} satisfies $\mu\not =0$ and 
    $\textnormal{genus}(\textnormal{spt}(\mu))=1$ - 
    as above in the second part of the theorem - 
    and additionally 
    $\Will(\mu)=\lim_{l\to \infty} \Will(F_{t_{j_l}})$
    for the subsequence $\{F_{t_{j_l}}\}$  
    from \eqref{weak.convergence.mu},   
    then there is another subsequence $\{F_{t_{j_k}}\}$ 
    of $\{F_{t_{j_l}}\}$ and some appropriate family of smooth diffeomorphisms 
    $\Theta_k:\Sigma \stackrel{\cong}\longrightarrow \Sigma$ 
    such that: 
    \begin{equation}   \label{W.2.2.convergence}
    \tilde F_{t_{j_k}}:=F_{t_{j_k}}\circ \Theta_k \longrightarrow f \quad \textnormal{in} \,\,\,W^{2,2}(\Sigma,g_{\textnormal{poin}}), \quad \textnormal{as} \,\, k\to \infty,   
    \end{equation}
    where $f$ is the uniformly conformal 
    bi-Lipschitz-parametrization of $\textnormal{spt}(\mu)$ from \eqref{parametrization.mu},
    and moreover each immersion $\tilde F_{t_{j_k}}$ in \eqref{W.2.2.convergence} is a uniformly bi-Lipschitz homeomorphism of $(\Sigma,g_{\textnormal{poin}})$ onto its image in $(\sphere^3,g_{\textnormal{euc}})$. Furthermore, there is for every fixed $x \in \sphere^3$ 
    some further subsequence 
    $\{\tilde F_{t_{j_{k_m}}}\}$ of the sequence   
    $\{\tilde F_{t_{j_k}}\}$, such that for any 
    $\varepsilon>0$ there is some sufficiently small 
    $\eta>0$ satisfying:      
    \begin{equation}  \label{no.A.concentration} 
    \int_{(\tilde F_{t_{j_{k_m}}})^{-1}
    (B^4_{\eta}(x)\cap \sphere^3)} 
    |A_{\tilde F_{t_{j_{k_m}}}}|^2 \,
    d\mu_{\tilde F_{t_{j_{k_m}}}^*g_{\textnormal{euc}}}
    <\varepsilon, \quad \forall \,m \in \nat.    
    \end{equation}  
    In particular, the measures 
    \begin{equation}  \label{no.concentration.Mill} 
    \Mill_l(\Omega):=
    \inf \Big{\{}\int_{F_{t_{j_l}}^{-1}(B\cap \sphere^3)} 
    |A_{F_{t_{j_l}}}|^2 \,d\mu_{F_{t_{j_l}}^*g_{\textnormal{euc}}}\, | \,B \supseteq \Omega \, \textnormal{and} \, 
    B \, \textnormal{is a Borel subset of} \,\,\rel^4 \Big{\}}
    \end{equation} 
    on $\rel^4$ {\bf do not concentrate at 
    any point of the ambient space $\rel^4$ as} 
    $l \to \infty$. 
    \item[4)] As in the second and third part of the 
    theorem we consider a limit varifold $\mu$ in \eqref{weak.convergence.mu} satisfying $\mu\not =0$ and 
    $\textnormal{genus}(\textnormal{spt}(\mu))=1$, 
    and we assume here again that the sequence $\{F_{t_{j_l}}\}$ in \eqref{weak.convergence.mu} satisfies \,   
    $\Will(\mu)=\lim_{l\to \infty} \Will(F_{t_{j_l}})$. 
    Suppose that $\{F_{t_{j_l}}\}$ also satisfies
    $\parallel |A^0_{F_{t_{j_l}}}|^2 \parallel_{L^{\infty}(\Sigma)} 
    \leq K$ and $|\frac{d}{dt}\Will(F_{t})|\lfloor_{t=t_{j_l}} 
    \leq K$ for all $l\in \nat$ and for some sufficiently 
    large number $K>1$, then the limit parametrization 
    $f$ of $\textnormal{spt}(\mu)$ in \eqref{parametrization.mu} is of class $W^{4,2}((\Sigma,g_{\textnormal{poin}}),\rel^4)$.        
\end{itemize}   
\end{theorem} 
\noindent 
Moreover, if we let the MIWF start moving in a smooth parametrization of a Hopf-torus - see Definition 
\ref{Hopf.torus.immersion} below - then we can derive more precise information from evolution equation \eqref{Moebius.flow}. For such flow lines of the MIWF we can actually rule out both degenerate cases $\mu=0$ and $\mu\not=0 \wedge \textnormal{spt}(\mu) \cong \sphere^2$ for singularities respectively for limit surfaces, as they appear in the first part of Theorem \ref{limit.MIWF}. Precisely, we prove the following theorem.    
\begin{theorem} \label{singular.time.MIWF.Hopf.tori}
Let $\{F_t\}$ be some flow line of the MIWF starting in a smooth parametrization $F_0:\Sigma \longrightarrow \sphere^3$ of a smooth Hopf-torus in $\sphere^3$ with 
$\Will(F_0)\leq 8\pi$. Moreover, as in Theorem \ref{limit.MIWF} we fix some sequence $t_j \nearrow T_{\textnormal{max}}(F_0)$ arbitrarily.  
\begin{itemize} 
\item[1)] For any subsequence 
$\{F_{t_{j_l}}\}$ as in \eqref{weak.convergence.mu} 
the embedded Hopf-tori $F_{t_{j_l}}(\Sigma)$ converge in Hausdorff-distance to some embedded $C^1$-Hopf-torus in $\sphere^3$, which is the support of the limit varifold 
$\mu$ appearing in \eqref{weak.convergence.mu}. In particular, for any such subsequence $\{F_{t_{j_l}}\}$ all statements of the second part of Theorem \ref{limit.MIWF} about $\{F_{t_{j_l}}\}$ and its limit varifold $\mu$ hold here. Hence, such an embedded limit Hopf-torus possesses a uniformly conformal bi-Lipschitz and $W^{2,2}$-parametrization
$f:(\Sigma,g_{\textnormal{poin}}) \stackrel{\cong}\longrightarrow \textnormal{spt}(\mu)$ from statement \eqref{parametrization.mu}, with 
$f^*(g_{\textnormal{euc}}) = e^{2u} \,g_{\textnormal{poin}}$ 
and $\parallel u \parallel_{L^{\infty}(\Sigma)} \leq \Lambda=\Lambda(\{F_{t_{j_l}}\},\mu)<\infty$, and with $\Will(f)<8\pi$, where $g_{\textnormal{poin}}$ is some suitable zero scalar curvature and unit volume metric on $\Sigma$. 
\item[2)] We suppose that there is some large constant $K>0$, such that $\parallel \vec H_{F_{t_j},\sphere^3} \parallel_{L^{\infty}(\Sigma)}$ remains uniformly bounded by $K$ for all $j \in \nat$, and we consider a subsequence $\{F_{t_{j_l}}\}$ of $\{F_{t_{j}}\}$ meeting property \eqref{weak.convergence.mu}, as in the first part of this theorem. Then any weakly/weakly* convergent sequence 
$\{\tilde F_{t_{j_k}}\}$ as in \eqref{weak.convergence.2} and \eqref{W.1.infty.convergence}, which we have obtained from $\{F_{t_{j_l}}\}$ in the second part of Theorem \ref{limit.MIWF}, can be uniformly estimated:
\begin{eqnarray}  \label{L.infty.W.4.2.estimate} 
\parallel \nabla^{g_{\textnormal{poin}}}(\tilde F_{t_{j_k}}) \parallel_{W^{3,2}(\Sigma,g_{\textnormal{poin}})}^2 \leq \\
\leq C(\Will(F_0),K,\Sigma,
g_{\textnormal{poin}},\Lambda)\cdot
\Big{(} \int_{\Sigma} |\nabla_{L^2} \Will(\tilde F_{t_{j_k}})|^2 \, d\mu_{(\tilde F_{t_{j_k}})^*(g_{\textnormal{euc}})} \,+1 \,\Big{)},  \nonumber
\end{eqnarray}
for every $k\in \nat$, where $g_{\textnormal{poin}}$ and 
$\Lambda$ are as in the first part of this theorem. 
\item[3)] We assume the same requirements as in part (2) 
of this theorem and additionally that the speed of 
``energy decrease'' $|\frac{d}{dt}\Will(F_{t})|$ shall remain uniformly bounded at the prescribed points of time
$t=t_{j} \nearrow T_{\textnormal{max}}(F_0)$.
Then, any weakly/weakly* convergent sequence 
$\{\tilde F_{t_{j_k}}\}$ as in \eqref{weak.convergence.2} and 
\eqref{W.1.infty.convergence} - as considered above 
in the second part - converges also weakly in $W^{4,2}((\Sigma,g_{\textnormal{poin}}),\rel^4)$,
strongly in $W^{3,2}((\Sigma,g_{\textnormal{poin}}),\rel^4)$ 
and in $C^{2,\alpha}((\Sigma,g_{\textnormal{poin}}),\rel^4)$, 
for any fixed $\alpha \in (0,1)$,   
to the uniformly conformal parametrization 
$f:\Sigma \stackrel{\cong} 
\longrightarrow \textnormal{spt}(\mu)$ of the corresponding 
limit Hopf-torus $\textnormal{spt}(\mu)$ from the first part of this theorem, and statement \eqref{no.A.concentration} holds here for the sequence $\{\tilde F_{t_{j_k}}\}$ respectively for the original sequence $\{F_{t_{j_k}}\}$ 
\footnote{This result confirms our geometric intuition, 
because every Hopf-torus consists of great circles in $\sphere^3$.}. 
In particular, the bi-Lipschitz parametrization 
$f$ of the limit Hopf-torus $\textnormal{spt}(\mu)$ 
is not only of class $W^{2,2}(\Sigma,\rel^4)$ but of class $W^{4,2}(\Sigma,\rel^4)$.  
\end{itemize}   
\end{theorem} 
\noindent 
Finally, we infer from Theorem \ref{singular.time.MIWF.Hopf.tori} above, from Proposition 
6 in \cite{Ruben.MIWF.II} and from Theorem 1 
in \cite{Ruben.MIWF.IV} the following dichotomy, which
contains a criterion for {\bf full convergence of 
global flow lines of the MIWF} starting in Hopf-tori in $\sphere^3$, only requiring the two additional conditions that they shall start moving below Willmore energy   
$4\pi^2$ and that the mean curvature vectors along those flow lines shall remain uniformly bounded in $L^{\infty}$.     
\begin{theorem} \label{limit.at.infinity}
Let $\{F_t\}_{t\geq 0}$ be a global flow line of the MIWF starting in some smooth parametrization $F_0:\Sigma \longrightarrow \sphere^3$ of a smooth Hopf-torus in $\sphere^3$ with $\Will(F_0)< 4\pi^2$. 
If there is a constant $K>0$, such that   
$\parallel \vec H_{F_{t},\sphere^3} \parallel_{L^{\infty}(\Sigma)}$ remains uniformly 
bounded by $K$ for all $t\in [0,\infty)$, 
then one and only one of the following two possibilities
will hold: 
\begin{itemize} 
\item[1)] The Willmore energy $\Will(F_t)$ strictly monotonically decreases to some value $v\geq 8 \pi$, 
and the flow line $\{F_t\}_{t\geq 0}$ \emph{diverges} 
as $t\nearrow \infty$ in the precise sense, 
that no smooth parametrization of the given flow line $\{F_t\}$ can fully converge in $C^4(\Sigma,\rel^4)$ to some 
$C^4$-immersion $F^*:\Sigma \longrightarrow \sphere^3$, 
as $t\nearrow \infty$.    
\item[2)] There is some smooth family of smooth diffeomorphisms $\Theta_t:\Sigma \stackrel{\cong} \longrightarrow \Sigma$, such that for each $m\in \nat$ the reparametrized flow line $\{F_t\circ \Theta_t\}_{t\geq 0}$ converges fully in $C^m(\Sigma,\rel^4)$ to a smooth and 
diffeomorphic parametrization of some torus in $\sphere^3$ which is conformally equivalent to the standard Clifford 
torus in $\sphere^3$, and this convergence takes place at an exponential rate as $t\nearrow \infty$.   
\end{itemize}      
\end{theorem} 
\noindent 
One should compare the first alternative in 
Theorem \ref{limit.at.infinity} with Theorem 
\ref{Singularities.do.exist} below, demonstrating  
that the first alternative in Theorem \ref{limit.at.infinity} 
might actually occur ! \\
On account of Propositions \ref{Hopf.Willmore.prop} and \ref{correspond.flows} the proof of Theorem \ref{limit.at.infinity} can easily be modified 
into a proof of the following counterpart of 
Theorem \ref{limit.at.infinity} for the degenerate flow \eqref{elastic.energy.flow} below.
\begin{corollary} \label{limit.at.infinity.1}
Let $\{\gamma_t\}_{t\geq 0}$ be a global flow line of the degenerate elastic energy flow \eqref{elastic.energy.flow}, starting in some smooth, closed and regular path 
$\gamma_0:\sphere^1 \longrightarrow \sphere^2$ 
with elastic energy 
$\Wil(\gamma_0):= \int_{\sphere^1} 
1 + |\vec \kappa_{\gamma_0}|^2 \,d\mu_{\gamma_0} 
< 4\pi$; see here Proposition 
\ref{compute.the.operator} below. 
If there is a constant $K>0$, such that   
the maximal curvature $\parallel \vec \kappa_{\gamma_{t}} \parallel_{L^{\infty}(\sphere^1)}$ remains uniformly 
bounded by $K$ for all $t\in [0,\infty)$, 
then one and only one of the following two possibilities
will hold: 
\begin{itemize} 
\item[1)] The elastic energy $\Wil(\gamma_t)$ strictly monotonically decreases to some value $v\geq 8$ 
and the flow line $\{\gamma_t\}_{t\geq 0}$ \emph{diverges} 
as $t\nearrow \infty$, in the precise sense that no smooth parametrization of the given flow line $\{\gamma_t\}$ can fully converge in $C^4(\sphere^1,\rel^3)$ to some 
regular curve $\gamma^*:\sphere^1 \longrightarrow \sphere^2$
of class $C^4$, as $t\nearrow \infty$.      
\item[2)] There is some smooth family of smooth diffeomorphisms $\sigma_t:\sphere^1 \stackrel{\cong} \longrightarrow \sphere^1$ such that for each $m\in \nat$ the reparametrized flow line $\{\gamma_t\circ \sigma_t\}_{t\geq 0}$ converges fully in $C^m(\sphere^1,\rel^3)$ to a smooth embedding of some great circle in $\sphere^2$,
and this convergence takes place at an exponential rate,
as $t\nearrow \infty$.       
\end{itemize}    
\end{corollary}
\noindent 
Similarly to Theorem \ref{Singularities.do.exist} below, 
Corollary \ref{limit.at.infinity.1} highlights the stark 
qualitative difference between the degenerate elastic 
energy flow \eqref{elastic.energy.flow} on $\sphere^2$ 
and its classical counterpart \eqref{classical.elastic.energy.flow} - a uniformly parabolic and subcritical flow of $4$th order whose flow lines 
have to fully converge - up to smooth reparametrization - to  
elastic curves in $\sphere^2$; compare here with the main result in \cite{Dall.Acqua.Pozzi.2018} and with the proof of Part III of Theorem 1 in \cite{Ruben.MIWF.II}, pp. 24--27.   

\section{Definitions and preparatory remarks}
 
First of all, we recall the following fundamental definitions, where parts (b) and (c) of Definition \ref{maximal.interval} are motivated by the classical terminology of Section 8.2 in \cite{Andrews.Hopper.2011}, examining the Ricci flow.
\begin{definition}  \label{maximal.interval} 
Let $\Sigma$ be a smooth compact torus and $n \geq 3$ an integer. 
\begin{itemize} 
	\item[a)] A \emph{flow line} of the MIWF \eqref{Moebius.flow} in the ambient manifold $M=\rel^n$ or $M=\sphere^n$ is a smooth family $\{F_t\}_{t\in [0,T)}$ of smooth immersions of  $\Sigma$ into $M$ such that the resulting smooth function 
	$F:\Sigma \times [0,T)\longrightarrow M$ satisfies equation \eqref{Moebius.flow} classically on $\Sigma \times [0,T)$. 	
	\item[b)] Let $F_0:\Sigma \longrightarrow M$ be a  
	smooth and umbilic-free immersion and $\{F_t\}_{t\in [0,T)}$ a smooth flow line of the MIWF starting in $F_0$. We call $[0,T)$ the \emph{interval of maximal existence} of the MIWF starting in $F_0$, if either $T=\infty$, or if there holds $T<\infty$ and for every $\varepsilon>0$ 
	there does not exist any smooth solution 
	$\{\tilde F_t\}_{t\in [0,T+\varepsilon)}$ 
	of the MIWF satisfying $\tilde F_t= F_t$ on $\Sigma$ for 
	$t\in [0,T)$. In both cases the element 
	$T \in \rel \cup \{\infty\}$ 
	is uniquely determined by the initial immersion $F_0$, and we call it the \emph{maximal time of existence} of the MIWF starting in $F_0$, in symbols: ``$T_{\textnormal{max}}(F_0)$''. 
	\item[c)] If $T_{\textnormal{max}}(F_0)$ is finite, 
	then we call $T_{\textnormal{max}}(F_0)$   
	\emph{the singular time} of the flow line $\{F_t\}$ 
	of the MIWF starting in $F_0$. In this case we will also 
	say that the flow line $\{F_t\}$ is ``singular'' or  
	``forms a singularity in finite time''.
	\end{itemize}  
\end{definition}  
\noindent 
\begin{remark} \label{divergent.flow}
	\begin{itemize} 
		\item[1)] Concerning Part (c) of the above definition 
		we should explain here that also a global flow line 
		of the MIWF might possibly ``diverge as $t \nearrow \infty$'' or ``form a singularity in infinite time'', if it either does not fully converge in $C^4(\Sigma)$ as $t \nearrow \infty$ or if it does not smoothly 
		subconverge to a Willmore immersion. This is obviously 
		not a consistent definition, but still both these 
		- slightly distinct - notions convey 
		the intuitive idea of a ``singularity being formed as 
		$t \nearrow \infty$'', and actually the first notion already appeared above in Theorem \ref{limit.at.infinity}, whereas the second notion will only appear below in Theorem \ref{Singularities.do.exist}.         
		\item[2)] We should also note here that all three parts of Definition \ref{maximal.interval} make sense because of Theorem 1 in \cite{Jakob_Moebius_2016} respectively Theorems 2 and 3 in \cite{Ruben.MIWF.III}, 
		proving existence and uniqueness of smooth 
		short-time solutions of the MIWF with 
		$C^{\infty}$-smooth, umbilic-free initial immersions 
		of a smooth torus into $\rel^n$ respectively $\sphere^n$.
	\end{itemize} 
\end{remark} 
\noindent 
Now we turn our attention to some basic differential geometric 
terms. As mentioned already in the introduction, we endow the unit $n$-sphere with the Euclidean scalar product of $\rel^{n+1}$, i.e. we set $g_{\sphere^n}:=\langle \,\cdot\,,\,\cdot\,\rangle_{\rel^{n+1}}$.
In Definition \ref{definitions} and Remark \ref{complex.conformal} we shall closely follow the classical book \cite{Jost.2006} about \emph{compact Riemann surfaces}, 
Rivi\`ere's lecture notes \cite{Riviere.Park.City.2013} and 
Tromba's introduction \cite{Tromba.Teichmueller.1992} into
\emph{Teichm\"uller Theory}, in order to recall here some
standard terminology and some basic facts. 
\begin{definition}  \label{definitions}
\begin{itemize}
\item[1)] We term two conformal atlases on a compact and 
orientable smooth surface $\Sigma$ equivalent, if 
their union is again a conformal atlas of $\Sigma$, 
and an equivalence class of conformal atlases
on $\Sigma$ is called a \emph{complex structure on $\Sigma$}. Having endowed $\Sigma$ with one of its complex structures, 
we call $\Sigma$ a \emph{compact Riemann surface}. 
\item[2)] Let $\Sigma_1$ and $\Sigma_2$ be two compact 
Riemann surfaces of the same genus $p\geq 0$. 
We term their complex structures $S_1(\Sigma_1)$ and $S_2(\Sigma_2)$ \emph{isomorphic}, if there is a 
\emph{conformal diffeomorphism}
$f:\Sigma_1 \stackrel{\cong}\longrightarrow \Sigma_2$, 
i.e. a bijective holomorphic map $f$ between $\Sigma_1$ and $\Sigma_2$ whose differential $D_xf:T_x\Sigma_1 \stackrel{\cong}\longrightarrow T_{f(x)}\Sigma_2$ is an invertible complex linear map between corresponding tangent spaces. 
\item[3)] The set of all isomorphism classes of complex structures on compact Riemann surfaces of some fixed genus $p\geq 0$ is called \emph{the moduli space of compact Riemann surfaces of genus $p$}, and we will denote this set by $\Mod_p$.
\item[4)] A \emph{conformal class} respectively 
\emph{conformal structure} on a compact and 
orientable smooth surface $\Sigma$ is a set $[g_0]$ 
of Riemannian metrics $g=e^{2u} \,g_0$ on $\Sigma$, where $g_0$ is a fixed Riemannian metric and $u\in C^{\infty}(\Sigma)$ an arbitrary smooth function on $\Sigma$, the \emph{conformal factor} of $g$ with respect to $g_0$.     
\end{itemize} 
\end{definition}
\noindent 
\begin{remark}  \label{complex.conformal}
\begin{itemize}
\item[1)]  
Any complex structure $S(\Sigma)$ on some  
compact and orientable smooth surface $\Sigma$ 
automatically yields a \emph{conformal class} $[g]$ of Riemannian metrics on $\Sigma$ which are compatible with any atlas $\mathcal{A}$ of $\Sigma$ representing the given complex structure $S(\Sigma)$, in the sense that there holds
$\psi_i^*(g_{\textnormal{euc}})=e^{-2v_i}\,g$ on 
$\Omega_i$, with $v_i\in L^{\infty}(\Omega_i)$, for any chart 
$\psi_i:\Omega_i \stackrel{\cong}\longrightarrow B_1^2(0)$
of the conformal atlas $\mathcal{A}$; 
see here Lemma 2.3.3 in \cite{Jost.2006}. 
Hence, a complex structure on any compact and orientable smooth surface $\Sigma$ automatically yields a conformal structure on $\Sigma$ in a canonical way, and actually also the converse holds on account of Theorem 3.11.1 in \cite{Jost.2006} respectively Theorems 2.9 and 2.13 in \cite{Riviere.Park.City.2013}, i.e. every Riemannian metric $g$ on a compact and orientable smooth surface $\Sigma$ yields a certain \emph{conformal atlas} $\mathcal{B}$ on $\Sigma$ 
which has the special property to only consist of \underline{isothermal} charts $\psi_i:\Omega_i \stackrel{\cong}\longrightarrow B_1^2(0)$ 
with respect to the prescribed metric $g$. 
Therefore, the two seemingly unrelated concepts of a compact and orientable smooth surface endowed with a complex structure 
respectively with a conformal class of Riemannian metrics 
coincide and thus will be used synonymously. 
\item[2)] It moreover follows from ``Poincar\'e's Theorem''  
- see Section 1.5 in \cite{Tromba.Teichmueller.1992} 
for an elegant proof treating the case ``$\textnormal{genus}(\Sigma)>1$'' -  
that every prescribed conformal class $[g_0]$ of Riemannian metrics on a compact and orientable smooth surface $\Sigma$ 
of genus $p\geq 1$ contains a \underline{unique smooth metric} $g_{\textnormal{poin}}$ of \underline{constant scalar curvature} $K_{g_{\textnormal{poin}}}$ and 
\underline{unit volume}, 
i.e. such that $K_{g_{\textnormal{poin}}} \equiv \textnormal{const}(p)\in \rel$ on $\Sigma$ and with $\mu_{g_{\textnormal{poin}}}(\Sigma)=1$. 
This can be easily seen, if we choose some 
$g$ from the prescribed conformal class $[g_0]$   
and consider the canonical \emph{Ansatz} $g_{\textnormal{poin}}:=e^{-2u}\,g$ for 
the unknown \emph{conformal factor} $u$ on $\Sigma$. 
Poincar\'e's Theorem now follows from the fact that 
Liouville's elliptic PDE 
\begin{equation}  \label{Liouville}  
-\triangle_{g}(u) + K_{g_{\textnormal{poin}}} \, e^{-2u} 
= K_g \quad \textnormal{on} \,\,\, \Sigma 
\end{equation}  
possesses for every prescribed 
negative number $K_{g_{\textnormal{poin}}}$ a unique 
smooth solution $u$ and in the special case 
$K_{g_{\textnormal{poin}}}\equiv 0$ a 
unique one-parameter family $\{u+r\}_{r\in \rel}$ 
of smooth solutions. Integrating equation \eqref{Liouville} 
over $\Sigma$ with respect to $\mu_g$, one infers from Gauss-Bonnet's theorem that the constant $K_{g_{\textnormal{poin}}}$ is actually determined by: 
\begin{equation}  \label{Gauss.Bonnet.K}
K_{g_{\textnormal{poin}}} = 
\frac{2\pi \,\chi(\Sigma)}{\mu_{g_{\textnormal{poin}}}(\Sigma)} 
= \frac{4\pi \,(1-p)}{\mu_{g_{\textnormal{poin}}}(\Sigma)},   
\end{equation}
see also formula (3.5) in \cite{Kuwert.Schaetzle.2012}. 
\end{itemize}    	        
\end{remark}
\noindent 
On account of this last remark we know in particular 
that every smooth immersion $f$ of a compact surface 
$\Sigma$ into some $\rel^n$, $n \geq 3$, yields a Riemannian 
metric $g_f:=f^*g_{\textnormal{euc}}$ which is 
conformally equivalent to a uniquely determined
smooth metric $g_{\textnormal{poin}}$ on $\Sigma$ 
of constant scalar curvature and unit volume, 
i.e. such that  $f^*g_{\textnormal{euc}}=e^{2u}\,g_{\textnormal{poin}}$
holds for some function $u\in C^{\infty}(\Sigma)$. 
This basic result will play a central role in the 
proofs of our first two main results, Theorems \ref{limit.MIWF} and 
\ref{singular.time.MIWF.Hopf.tori} below. 
Now by Theorem 4.3 and Corollary 4.4 in \cite{Riviere.Park.City.2013} respectively Theorem VII.12 
in \cite{Riviere.2011} or also by Theorem 1.4 in 
\cite{Riviere.Var.principle.2014} - building on work by 
M\"uller and Sver\'ak \cite{Mueller.Sverak.1995} and 
later H\'elein \cite{Helein.2004} about a geometric application of ``compensated compactness''-estimates and 
quasiconformal mapping theory - this result still holds 
more generally for immersions $f:\Sigma \longrightarrow \rel^n$ of class $W^{1,\infty}$ with second fundamental 
form $A_f$ in $L^2(\Sigma)$ in the precise sense,  
that such an immersion $f$ still gives rise to a certain  
complex structure on $\Sigma$ uniquely determining 
some smooth metric $g_{\textnormal{poin}}$ of constant scalar curvature and unit volume on $\Sigma$, such that at least up to reparametrization by some Lipschitz-homeomorphism $\Psi:\Sigma \stackrel{\cong}\longrightarrow \Sigma$ there holds $(f\circ\Psi)^*g_{\textnormal{euc}}=
e^{2u}\,g_{\textnormal{poin}}$ on $\Sigma$, 
for some $u \in L^{\infty}(\Sigma)$.  
We shall therefore introduce the concept of 
\emph{uniformly conformal Lipschitz immersions} as in Section 2 of \cite{Kuwert.Schaetzle.conf.class.2013}.    
\begin{definition} \label{conformal.immersion} 
Let $g_{\textnormal{poin}}$ be a smooth constant 
scalar curvature and unit-volume metric on some smooth 
compact, orientable surface $\Sigma$.  
\begin{itemize} 
\item[1)] We call a Lipschitz immersion $f:\Sigma \longrightarrow \sphere^3$ ``uniformly conformal  
to $g_{\textnormal{poin}}$'', if there is some 
function $u\in L^{\infty}(\Sigma)$ such that 
$f^*g_{\textnormal{euc}}=e^{2u}\,g_{\textnormal{poin}}$. 
If the Lipschitz immersion $f$ is here additionally 
of class $W^{2,2}(\Sigma,\rel^4)$, then we call $f$ 
a $(W^{1,\infty} \cap W^{2,2})$-immersion being
uniformly conformal to $g_{\textnormal{poin}}$. 
\item[2)] Similarly, we define the above two notions 
for Lipschitz- respectively 
$(W^{1,\infty} \cap W^{2,2})$-immersions 
$f:\Sigma \longrightarrow \rel^3$, 
replacing $(\sphere^3,g_{\textnormal{euc}})$ by 
$(\rel^3,g_{\textnormal{euc}})$ and $\rel^4$ by $\rel^3$.
\end{itemize}   
\end{definition}
\noindent 
Motivated by Rivi\`ere's \cite{Riviere.2008}, 
\cite{Riviere.2011} and Bernard's \cite{Bernard.2016} work 
on the divergence form of the Euler-Lagrange-operator 
$\nabla_{L^2}\Will$ of the Willmore functional 
we shall follow here exactly Section 7 in \cite{Riviere.Park.City.2013}, 
Sections VII.5.2 and VII.6.5 in \cite{Riviere.2011}, or 
Section 1.5 in \cite{Palmurella.2022} and define for any 
fixed Lipschitz-immersion 
$F:\Sigma \longrightarrow \rel^3$ with second fundamental 
form $A_F$ in $L^2(\Sigma)$ its \emph{weak Willmore operator} below in \eqref{distributional.Willmore} as a distribution of second order on $\Sigma$, thus generalizing its classical meaning in \eqref{first_variation} for smooth immersions $f$ of $\Sigma$ into $\rel^3$. 
The differential geometric terminology appearing below  
in Definition \ref{definitions.3} and in 
Remark \ref{Immersion.in.R4} is fairly standard 
and can be looked up in Definition 4 of \cite{Ruben.MIWF.II} and in Section 2 of both \cite{Jakob_Moebius_2016} and 
\cite{Ruben.MIWF.III}. 
\begin{definition}   \label{definitions.3}
Let $\Sigma$ be a compact Riemann surface, 
and let $F:\Sigma \longrightarrow \rel^3$ be 
a Lipschitz immersion with second fundamental form 
$A_F$ in $L^2(\Sigma)$ and with induced metric  $g_{F}:=F^*{g_{\textnormal{euc}}}$.
\begin{itemize}  
\item[1)] We define the ``weak Willmore operator'' $\nabla_{L^2}\Will(F)$ as the following distribution on $\Sigma$ of second order:
\begin{eqnarray}  \label{distributional.Willmore} 
\langle \nabla_{L^2}\Will(F),\varphi \rangle_{\Dom'(\Sigma)} 	 
:= \frac{1}{2} \, \int_{\Sigma} 
\langle \vec H_{F}, \triangle_F \varphi \rangle_{g_{\textnormal{euc}}} 
-  \, g_F^{\nu\alpha} \, g_F^{\mu\xi}\,\langle (A_F)_{\xi\nu}, 
\vec H_{F} \rangle_{g_{\textnormal{euc}}} \,\langle \partial_{\mu}F, \partial_{\alpha}\varphi \rangle_{g_{\textnormal{euc}}} \nonumber \\
- \, g_F^{\nu\alpha} \, g_F^{\mu\xi}\,\langle (A_F^0)_{\xi\nu}, \vec H_{F} \rangle_{g_{\textnormal{euc}}} \, \langle \partial_{\mu}F, \partial_{\alpha}\varphi \rangle_{g_{\textnormal{euc}}}
\, d\mu_{g_F}  \quad 
\end{eqnarray}  
$\forall \, \varphi \in C^{\infty}(\Sigma,\rel^3)$. 
\item[2)] Moreover, as in \cite{Riviere.2011}, Definition VII.3, we shall define the differential $1$-form 
$$
w_F:\partial_{\nu} \mapsto \frac{1}{2} \,
\Big{(} \nabla^F_{\nu}(\vec H_F) + \, g_F^{\mu\xi}\,
\langle (A_F)_{\xi\nu}, \vec H_{F} \rangle_{g_{\textnormal{euc}}} \,\partial_{\mu}F 
+ \, g_F^{\mu\xi}\,\langle (A_F^0)_{\xi\nu}, \vec H_{F} \rangle_{g_{\textnormal{euc}}} \,\partial_{\mu}F \Big{)}
$$ 
mapping smooth vector fields on $\Sigma$ into 
distributions of first order on $\Sigma$, 
acting on sections into $F^*T\rel^3$ concretely by:
\begin{eqnarray}  \label{distributional.w.form} 
\langle w_F(\partial_{\nu}),\varphi \rangle_{\Dom'(\Sigma)} 	 
:= \frac{1}{2} \,\int_{\Sigma} 
-\langle \vec H_{F}, \nabla^F_{\nu} \varphi \rangle_{g_{\textnormal{euc}}} 
+   \, g_F^{\mu\xi}\,\langle (A_F)_{\xi\nu}, 
\vec H_{F} \rangle_{g_{\textnormal{euc}}} \,
\langle \partial_{\mu}F, \varphi \rangle_{g_{\textnormal{euc}}} 
\nonumber \\
+ \, g_F^{\mu\xi}\,\langle (A_F^0)_{\xi\nu}, \vec H_{F} \rangle_{g_{\textnormal{euc}}} \, \langle \partial_{\mu}F, \varphi \rangle_{g_{\textnormal{euc}}}
\, d\mu_{g_F},  \quad 
\end{eqnarray}   
$\forall \, \varphi \in C^{\infty}(\Sigma,\rel^3)$, for $\nu=1,2$. As in \cite{Palmurella.2022}, formula (2.5), 
we will also use the short notation:
\begin{equation}  \label{distributional.w.form.2} 
w_F = \frac{1}{2}\Big{(}
\nabla^F(\vec H_F) + \,
\langle (A_F), \vec H_{F} \rangle_{g_{\textnormal{euc}}}^{\sharp_{g_F}} 
+ \,\langle (A_F^0), \vec H_{F} \rangle_{g_{\textnormal{euc}}}^{\sharp_{g_F}}\Big{)},
\end{equation}  
for the distributional $1$-form in \eqref{distributional.w.form}. 
\end{itemize}
\end{definition}
\noindent 
\begin{remark}  \label{Immersion.in.R4}
Concerning the definition of the ``weak Willmore operator'' in the first part of Definition \ref{definitions.3} we should mention here, that the definition in line \eqref{distributional.Willmore}    
becomes an identity for smooth immersions 
$F:\Sigma \longrightarrow \rel^3$ on account of 
Theorem VII.7 in \cite{Riviere.2011}, respectively 
on account of the main theorem in \cite{Bernard.2016}, originating from Theorem 1.1 in \cite{Riviere.2008}. 
We only have to recall the fact, that for smooth $F$ the $L^2$-gradient  $\nabla_{L^2}\Will(F)$ of $\Will$
can be directly computed, yielding the differential operator
in \eqref{Euler.Lagrange.operator}, which coincides 
with the expression in \eqref{first_variation} respectively \eqref{distributional.Willmore} via integration by parts and the main theorem in \cite{Bernard.2016}:
\begin{equation}  \label{Riviere.Bernard.Willmore} 
\langle \nabla_{L^2}\Will(F),\varphi 
\rangle_{L^2(\Sigma,\mu_{g_{F}})} 
=\langle \nabla_{L^2}\Will(F),\varphi \rangle_{\Dom'(\Sigma)} 
\qquad \forall \, \varphi \in C^{\infty}(\Sigma,\rel^3).   
\end{equation}    
For smooth immersions into $\sphere^3$ this is no longer true. In order to see this, but also in view of the proof of 
Theorem \ref{limit.MIWF}, we should point out here some precise, basic facts from Differential Geometry:   
On account of Lemma 2.1 in \cite{Ndiaye.Schaetzle.2014} 
the Willmore functional of a $C^{\infty}$-smooth immersion 
$F:\Sigma \longrightarrow \sphere^n$, $n \geq 3$, 
the pullback metric induced by $F$, the tracefree part 
$A^{0}_{F}$ of the second fundamental form of $F$, its 
squared length $|A^{0}_{F}|^2$,  
and the classical Willmore Lagrange operator
\begin{equation}   \label{Euler.Lagrange.operator} 
\nabla_{L^2}\Will(F) =
\frac{1}{2} \,\Big{(} \triangle_{F}^{\perp} \vec H_{F} + Q(A^{0}_{F})(\vec H_{F}) \Big{)} 
\end{equation}
remain unchanged, if the immersion 
$F:\Sigma \longrightarrow \sphere^{n}\subset \rel^{n+1}$ 
is interpreted as an immersion of $\Sigma$ into $\rel^{n+1}$, 
although this distinction has a certain effect on its 
entire second fundamental form and on its mean curvature vector. Denoting by $A_{F,{\rel^{n+1}}}$ and $\vec H_{F,{\rel^{n+1}}}$ the second fundamental form and the mean curvature vector of $F:\Sigma \longrightarrow \sphere^{n}$ 
interpreted as an immersion into $\rel^{n+1}$ and 
similarly by $A_{F,{\sphere^{n}}}$ and 
$\vec H_{F,{\sphere^{n}}}$ the second fundamental form 
and the mean curvature vector of $F:\Sigma \longrightarrow \sphere^{n}$ interpreted as an immersion into $\sphere^{n}$,   
we have by formulae (2.1)--(2.4) in 
\cite{Ndiaye.Schaetzle.2014}:
\begin{eqnarray}
g_{F,\rel^{n+1}}=g_{F,\sphere^{n}}	          \label{no1}   \\
A_{F,\rel^{n+1}}=A_{F,\sphere^{n}} - F\,g_F   \label{no2}   \\
A^0_{F,\rel^{n+1}}=A^0_{F,\sphere^{n}}        \label{no3}   \\
\vec H_{F,\rel^{n+1}}= \vec H_{F,\sphere^{n}} - 2F       \label{no4}   \\ 
\triangle_{F,\rel^{n+1}}^{\perp} \vec H_{F,\rel^{n+1}}     
+ Q(A^{0}_{F,\rel^{n+1}})(\vec H_{F,\rel^{n+1}})
=\triangle_{F,\sphere^{n}}^{\perp} \vec H_{F,\sphere^{n}} 
+ Q(A^{0}_{F,\sphere^{n}})(\vec H_{F,\sphere^{n}}).     \label{no5} 
\end{eqnarray}  
\end{remark}
\noindent \\  
Finally, we will also need:
\begin{definition}  \label{Umbilic.free}
 Let $\Sigma$ be a compact smooth torus and 
$n \geq 3$ an integer. We denote by $\textnormal{Imm}_{\textnormal{uf}}(\Sigma,\rel^n)$
the subset of $C^2(\Sigma,\rel^n)$ consisting of 
\underline{umbilic-free immersions}, i.e.:
$$ 	
\textnormal{Imm}_{\textnormal{uf}}(\Sigma,\rel^n):=
\{\,f \in C^{2}(\Sigma,\rel^n)\,|\, f \,\,\textnormal{is an immersion satisfying}  \mid A^0_{f} \mid^{2}>0 \,\, \textnormal{on} \,\,\Sigma \,\}.
$$
\end{definition} 
\noindent \\
At this point we also mention and apply the first part of Theorem 4.1 in Palmurella's and Rivi\`ere's paper \cite{Palmurella.2022}, a tricky regularity theorem for conformal $(W^{2,2}\cap W^{1,\infty})$-immersions 
$f:B^2_1(0) \longrightarrow \rel^3$ with distributional Willmore operator $\nabla_{L^2}\Will(f)$ from 
\eqref{distributional.Willmore} of class $L^2(B^2_1(0),\rel^3)$, in order to achieve the analogous statement for uniformly conformal 
$(W^{2,2}\cap W^{1,\infty})$-immersions 
$F:\Sigma \longrightarrow \rel^3$ whose distributional 
Willmore operator $\nabla_{L^2}\Will(F)$ from 
\eqref{distributional.Willmore} is of class $L^2((\Sigma,g_{\textnormal{poin}}),\rel^3)$.
We shall not repeat here the proof of Theorem 4.1 in \cite{Palmurella.2022} whose mathematical background 
can be found in \cite{Riviere.2011}, \cite{Riviere.Park.City.2013}, \cite{Bernard.2016} and  \cite{Palmurella.2022} itself.  
\begin{theorem}  \label{Regul.theorem} 
Let $\Sigma$ be a compact Riemann surface of genus 
$p \geq 1$ and $g_{\textnormal{poin}}$ 
the smooth metric of constant scalar curvature and 
unit volume representing the given conformal class of $\Sigma$. Then any uniformly conformal 
$(W^{2,2}\cap W^{1,\infty})$-immersion 
$F:\Sigma \longrightarrow \rel^3$, whose distributional 
Willmore operator $\nabla_{L^2}\Will(F)$ from 
\eqref{distributional.Willmore} can be identified 
with a function of class $L^2((\Sigma,g_{\textnormal{poin}}),\rel^3)$, is of class $W^{4,2}((\Sigma,g_{\textnormal{poin}}),\rel^3)$. 
\end{theorem}  
\emph{Sketch of a proof:\,}  
Since the given Riemann surface $\Sigma$ is compact, 
we have to show the asserted regularity of $F$ 
only locally. We therefore choose some coordinate patch 
$\Omega \subset \Sigma$ and some isothermal chart 
$\psi:B^2_1(0) \longrightarrow \Omega$, i.e. some chart 
$\psi$ satisfying 
$\psi^*(g_{\textnormal{poin}})=e^{2v} \,g_{\textnormal{euc}}$ 
on the open unit disc $B^2_1(0)$ for some smooth, bounded 
function $v$ on $B^2_1(0)$; see here the first part of 
Remark \ref{complex.conformal}.
Now, since there holds $F^*(g_{\textnormal{euc}})
=e^{2u}\, g_{\textnormal{poin}}$ on $\Sigma$ for some 
function $u\in L^{\infty}(\Sigma)$ - see here Definition \ref{conformal.immersion} - the composition 
$f:=F\circ \psi:B^2_1(0) \longrightarrow \rel^3$ is 
of class $(W^{2,2}\cap W^{1,\infty})(B_1^2(0),\rel^3)$ 
and satisfies:
\begin{equation}  \label{isothermal.coordinates} 
f^*(g_{\textnormal{euc}})=(F\circ \psi)^*(g_{\textnormal{euc}})
=\psi^*(e^{2u}\, g_{\textnormal{poin}}) 
= e^{2u\circ \psi+2v}\, g_{\textnormal{euc}}
\quad \textnormal{on} \,\,B^2_1(0). 
\end{equation} 
Moreover, the second requirement regarding the considered immersion $F$ implies that the restriction of $\nabla_{L^2}\Will(F)$ from 
\eqref{distributional.Willmore} to our coordinate 
patch $\Omega$ has to be a function of class 
$L^2((\Omega,g_{\textnormal{poin}}),\rel^3)$, hence  
its pullback $\nabla_{L^2}\Will(F)\circ \psi = 
\nabla_{L^2}\Will(f)$ via our isothermal chart $\psi$  
has to be a function of class $L^2(B^2_1(0),\rel^3)$.
Furthermore, we can easily see by means of formulae 
\eqref{distributional.Willmore} and \eqref{distributional.w.form} 
that the distributional covariant divergence of the distributional $1$-form $w_F$ in \eqref{distributional.w.form} and \eqref{distributional.w.form.2} is the distributional 
Willmore operator $\nabla_{L^2}\Will(F)$ in \eqref{distributional.Willmore}, i.e.: 
\begin{eqnarray*} 
\langle \textnormal{div}_F(w_F),\varphi \rangle_{\Dom'(\Sigma)} = 
\int_{\Sigma} \big{\langle} g_F^{\nu\alpha} \,\nabla^F_{\alpha}(w_F(\partial_{\nu})), \varphi \big{\rangle} \,d\mu_{g_F}                        
=  \langle \nabla_{L^2}\Will(F),\varphi \rangle_{\Dom'(\Sigma)} 
\end{eqnarray*}      
$\forall \, \varphi \in C^{\infty}(\Sigma,\rel^3)$.
We can therefore write here:
\begin{equation}   \label{d.star.wF}
\textnormal{div}_F(w_F) = \nabla_{L^2}\Will(F)
\qquad \textnormal{in} \,\, \Dom'(\Omega),
\end{equation}
as if these were smooth functions on the coordinate patch $\Omega$. Equation \eqref{d.star.wF} 
translates via our isothermal chart $\psi$ and equation 
\eqref{isothermal.coordinates} into the equation
\begin{equation}   \label{d.star.wf}
d^*(w_f) = e^{2u\circ\psi+2v}\, d^{*_{g_f}}(w_{f})   
=e^{2u\circ\psi+2v}\, \nabla_{L^2}\Will(f)  \qquad 
\textnormal{in} \,\, \Dom'(B^2_1(0)),   
\end{equation}
where we have also used the two adjoint operators 
$$
d^*=-\star \,d \,\star :\Omega^{1}(B^2_1(0))
\longrightarrow \Omega^{0}(B^2_1(0)) 
\,\, \textnormal{and} \,\,  
d^{*_{g_f}}:= -\star_{g_f} \,d \,\star_{g_f}: 
\Omega^{1}(B^2_1(0))\longrightarrow \Omega^{0}(B^2_1(0))
$$
of the exterior derivative  
$d:\Omega^0(B^2_1(0)) \longrightarrow \Omega^{1}(B^2_1(0))$
both with respect to the Euclidean metric and with respect to the pullback metric $g_f \equiv f^*(g_{\textnormal{euc}}) 
= e^{2u\circ \psi+2v}\, g_{\textnormal{euc}}$ 
on the open unit disc $B^2_1(0)$;
see here also Section 3.3 in \cite{Jost.2017}. 
Equation \eqref{d.star.wf} together with the 
$L^2$-regularity of the right hand side 
of \eqref{d.star.wf} are all required ingredients 
at the beginning of the proof of Proposition 4.6  
in \cite{Palmurella.2022}. Hence, we can apply here 
Propositions 4.6 and 4.3 in \cite{Palmurella.2022}  
to the conformal $(W^{2,2}\cap W^{1,\infty})(B^2_1(0))$-immersion $f= F \circ \psi$, 
in order to conclude firstly that the mean curvature vector $\vec H_f$ of $f$ is of class $L^p_{\textnormal{loc}}(B^2_1(0),\rel^3)$, 
for every $p\geq 1$, and then that $f$ is indeed of class  $W^{4,2}_{\textnormal{loc}}(B^2_1(0),\rel^3)$.

\section{Proof of Theorem  \ref{limit.MIWF}}
In parts (1)--(4) of Theorem \ref{limit.MIWF}  
every immersion $F$ of $\Sigma$ into $\sphere^3$ will be  
interpreted as an immersion of $\Sigma$ into $\rel^4$.	 
The corresponding effect of this choice on geometric 
tensors and scalars of such immersions is summarized in 
Remark \ref{Immersion.in.R4}. Hence, we shall simply  
write $F^*g_{\textnormal{euc}}$ instead of
$F^*\langle \cdot,\cdot \rangle_{\rel^4}$,   
$A_{F}$ instead of $A_{F,\rel^{4}}$, 
$A^0_{F}$ instead of $A^0_{F,\rel^{4}}$,
$\vec H_{F}$ instead of $\vec H_{F,\rel^{4}}$, etc..
First of all, we infer from the requirement
$\Will(F_0)\leq 8\pi$ and from the fact that 
\begin{equation} \label{Willmore.monotonicity}  
\frac{d}{dt}\Will(F_t) = 
\langle \partial_t F_t, \nabla_{L^2} \Will(F_t) 
\rangle_{L^2(\Sigma,\mu_{g_{F_t}})} 
= -\frac{1}{2} \frac{1}{|A^0_{F_t}|^4} \,
\,\int_{\Sigma} |\nabla_{L^2} \Will(F_t)|^2 \,d\mu_{g_{F_t}} 
\leq 0 
\end{equation}
holds for every $t\in [0,T_{\textnormal{max}}(F_0))$, i.e. 
from the well-known monotonicity of $t\mapsto \Will(F_t)$ 
along the considered flow line $\{F_t\}_{t\in [0,T_{\textnormal{max}}(F_0))}$ of the MIWF, 
that we should distinguish two different possibilities: 
(a) There holds $\Will(F_t)=8 \pi$ on some 
arbitrarily short, but non-empty time interval $[0,\varepsilon)$, or (b) there holds 
$\Will(F_t)<8 \pi$ for every 
$t\in (0,T_{\textnormal{max}}(F_0))$.
In the special case (a) we infer again from the computation 
in equation \eqref{Willmore.monotonicity}:
\begin{equation} \label{Willmore.stationary}
0 =\frac{d}{dt}\Will(F_t) 
= -\frac{1}{2} \frac{1}{|A^0_{F_t}|^4} \,
\,\int_{\Sigma} |\nabla_{L^2} \Will(F_t)|^2 \,d\mu_{g_{F_t}}, 
\,\,\, \textnormal{for every}\,\, t\in [0,\varepsilon),
\end{equation}
implying that $\nabla_{L^2} \Will(F_t)$ and  
therefore also $\partial_t F_t$ vanishes on $\Sigma$ 
for every $t\in [0,\varepsilon)$. This means that in case (a) 
$\{F_t\}$ has to be a stationary flow line of the MIWF 
consisting only of the Willmore-torus $F_0$ 
with Willmore energy $8\pi$, at least for some short time
$\varepsilon>0$. Now, recalling the regularity 
result of the second part of Theorem 3 in \cite{Ruben.MIWF.III} we know that the Willmore energy 
$t \mapsto \Will(F_t)$ is a real-analytic function 
along the entire flow line $\{F_t\}$, more precisely 
for every $t\in (0,T_{\textnormal{max}}(F_0))$.  
Hence, combining this result with the above conclusion 
we infer that in case (a) actually  
$\partial_t F_t=0$ has to hold for every 
$t\in [0,T_{\textnormal{max}}(F_0))$ and that moreover 
$T_{\textnormal{max}}(F_0)$ cannot be finite, i.e. that 
in case (a) $\{F_t\}_{t \in [0,T_{\textnormal{max}}(F_0))}$ 
is a global flow line of the MIWF which only 
consists of the Willmore torus $F_0$ having
Willmore energy $8\pi$. In this special case 
$\{F_t\}_{t \in [0,T_{\textnormal{max}}(F_0))}$ obviously  
does not produce any singularity as 
$t \nearrow T_{\textnormal{max}}(F_0)=\infty$, and 
all statements of this theorem are trivially true.  
We may therefore assume in the sequel of this proof
that there holds:  
\begin{equation} \label{Energy.smaller.8pi}
\Will(F_t)<8 \pi  \quad \forall \, 
t\in (0,T_{\textnormal{max}}(F_0)).
\end{equation}  
\begin{itemize} 
\item[1)] First of all, on account of assumption \eqref{Energy.smaller.8pi} we have here:     
\begin{equation}  \label{monotonic.Will}
\Will(F_{t_{j+1}})\leq \Will(F_{t_{j}})<8 \pi \quad 
\forall\, j \in \nat, 
\end{equation} 
which together with formula \eqref{no4} 
yields inequality A.20 in \cite{Kuwert.Schaetzle.Annals}
for each integral $2$-varifold $\mu_j:= \Hn\lfloor_{F_{t_j}(\Sigma)}$, and hence formulae 
A.17 and A.21 in \cite{Kuwert.Schaetzle.Annals} show 
that each immersion $F_{t_j}$ is a smooth embedding
of $\Sigma$ into $\sphere^3$. 
Moreover, we infer from the compactness of each embedded surface $F_{t_j}(\Sigma)$ and from formulae A.6 and A.16 in \cite{Kuwert.Schaetzle.Annals} applied to each integral $2$-varifold $\mu_j = \Hn\lfloor_{F_{t_j}(\Sigma)}$ that 
$$ 
\varrho^{-2} \, \Hn(F_{t_j}(\Sigma)\cap B_{\varrho}(x))
\leq \textnormal{const.} \,\Will(F_{t_{j}})  
\quad \forall \,j \in \nat,
$$ 
for every fixed $x\in \rel^4$ and every $\varrho>0$. 
Hence, combining this again with inequality \eqref{monotonic.Will} and formula \eqref{no4} we can infer from Allard's Compactness Theorem that there is some subsequence $\{F_{t_{j_l}}\}$ of $\{F_{t_j}\}$ such that  
\begin{equation}  \label{mu.in.the.limit} 
\Hn\lfloor_{F_{t_{j_l}}(\Sigma)} \longrightarrow 
\mu \quad \textnormal{weakly as Radon measures on} 
\,\,\, \rel^4,  
\end{equation} 
as $l\to \infty$, for some integral 
$2$-varifold $\mu$ on $\rel^4$, which is the asserted 
convergence \eqref{weak.convergence.mu}.   
Now we immediately suppose that the limit measure $\mu$ in  \eqref{mu.in.the.limit} is non-trivial, i.e. the case  in which $\textnormal{spt}(\mu)\not = \emptyset$.
\footnote{The special case in which the limit measure $\mu$ in \eqref{mu.in.the.limit} is trivial cannot be ruled out a-priori and does not allow any topological conclusion 
about the set of limit points of the sequence of embedded surfaces $F_{t_{j_l}}(\Sigma)$ in $\sphere^3$. 
One can quickly understand on account of the proofs of 
Proposition 4.1 and Theorem 5.2 in \cite{Kuwert.Li.2012} or on account of Theorem 5.3 in \cite{Riviere.Park.City.2013} that we cannot prove any reasonable type of convergence of the immersions $F_{t_{j_l}}$ in this special case, since we do not 
have the freedom to arbitrarily apply 
M\"obius transformations 
$\Psi_l: \sphere^3 \longrightarrow \sphere^3$ to the 
embedded surfaces $F_{t_{j_l}}(\Sigma)$ for infinitely 
many $l$.} 
Here we can combine \eqref{monotonic.Will} and \eqref{mu.in.the.limit} with the fact that there holds $F_{t_{j_l}}(\Sigma) \subset \sphere^3$, $\forall \, l\in \nat$, and we can argue as in 
the proof of Proposition 2.1 in \cite{Schaetzle.Conf.factor.2013}, respectively as in the proof of Theorem 4.2 in \cite{Kuwert.Schaetzle.2012},
using Simon's monotonicity formula (1.2) in \cite{Simon.1993} 
for integral $2$-rectifiable varifolds in $\rel^4$ as in the proof of Theorem 3.1 in \cite{Simon.1993}, pp. 310--311:
\begin{equation}  \label{Hausdorff.converge} 
F_{t_{j_l}}(\Sigma) \longrightarrow \textnormal{spt}(\mu) 
\,\,\, \textnormal{as subsets of $\rel^4$ in 
Hausdorff-distance, as} \,\,\, l\to \infty, 
\end{equation}
which is just the asserted convergence \eqref{Hausdorff.converg} and shows particularly 
that $\textnormal{spt}(\mu)$ has to be contained in $\sphere^3$, whenever $\mu \not =0$. 
Again because of $\mu\not =0$ we can apply here Proposition 2.1 in \cite{Schaetzle.Conf.factor.2013} and infer that $\textnormal{spt}(\mu)$ is an orientable and 
embedded Lipschitz-surface, having genus either $0$ or 
$1$, and that $\mu$ has constant Hausdorff-$2$-density $1$.	
\item[2)]
Now we suppose that the limit varifold $\mu$ in \eqref{weak.convergence.mu} respectively in \eqref{mu.in.the.limit} is non-trivial and that its 
support $\textnormal{spt}(\mu)$ is an embedded torus in $\sphere^3$. For ease of notation and in view of the 
aims of this part of the theorem we may relabel henceforth 
the subsequence $\{F_{t_{j_l}}\}$ satisfying  \eqref{weak.convergence.mu} and \eqref{Hausdorff.converg} 
into $\{F_{t_{j}}\}$. 
As explained in Remark \ref{complex.conformal} we obtain unique smooth metrics $g_{\textnormal{poin},j}$ of zero scalar curvature and unit volume on $\Sigma$, such that each immersion $F_{t_{j}}$ is uniformly conformal with respect to $g_{\textnormal{poin},j}$ in the sense of Definition \ref{conformal.immersion}, i.e. such that there holds:
\begin{equation}   \label{conformal.reparametrization}  
(F_{t_{j}})^*(g_{\textnormal{euc}}) 
= e^{2u_{j}} \,g_{\textnormal{poin},j} \quad 
\textnormal{on} \,\, \Sigma, \quad \forall \,j \in \nat,
\end{equation}
for unique functions $u_{j} \in C^{\infty}(\Sigma)$.
Now we shall try to prove that the smooth conformal 
factors $u_j$ appearing in \eqref{conformal.reparametrization} remain uniformly bounded on $\Sigma$ for all $j \in \nat$. 
First of all, on account of the assumptions in 
this part of the theorem we may apply 
Propositions 2.1 and 7.2 in \cite{Schaetzle.Conf.factor.2013} 
in order to see that there is a smooth compact torus 
$\tilde \Sigma$ and a homeomorphic   
$(W^{2,2}\cap W^{1,\infty})$-parametrization 
$F:\tilde \Sigma \stackrel{\cong}\longrightarrow \textnormal{spt}(\mu)$ of $\textnormal{spt}(\mu)$,
being uniformly conformal with respect to some smooth metric 
$\tilde g$ on $\tilde \Sigma$ and whose ``pushforward''-measure  
$F(\mu_{F^*(g_{\textnormal{euc}})})=:\mu_F$ 
coincides with $\mu$ on $\rel^4$ by 
formula (2.5) in \cite{Schaetzle.Conf.factor.2013}
and has square integrable weak mean curvature vector 
$\vec H_{\mu_F}$ with $\Will(\mu_F)=\Will(F)\in [4\pi,\infty)$ by formula (7.8) in \cite{Schaetzle.Conf.factor.2013}.
Moreover, on account of the proven Willmore conjecture, 
Theorem A in \cite{Marques.Neves.2014}, combined with 
Theorem 1.7 in \cite{Riviere.Var.principle.2014}, 
and by means of stereographic projection $\Stereo$ from 
$\sphere^3\setminus \{(0,0,0,1)\}$ into $\rel^3$ - 
assuming here that $\textnormal{spt}(\mu)$ is 
contained in $\sphere^3\setminus \{(0,0,0,1)\}$ 
without loss of generality on account of convergence 
\eqref{Hausdorff.converg} - we conclude from 
the above properties of $F$: 
\begin{eqnarray*}
\Will(F) = \Will(\Stereo \circ F) \geq \qquad  \\
\geq \min \{ \Will(f) \,|\, 
f \in (W^{2,2}\cap W^{1,\infty})(\tilde \Sigma,\rel^3), \,\,
f^*(g_{\textnormal{euc}}) \geq \varepsilon \,\tilde g, \,\, 
\textnormal{for some} \, \varepsilon>0 \, \} 
\geq 2\pi^2. 
\end{eqnarray*} 
Combining this with $\Will(\mu)=\Will(\mu_F)=\Will(F)$ 
and with estimate \eqref{monotonic.Will}
we obtain: 
\begin{eqnarray}   \label{Abschaetz.nerd}
\Will(\mu) + \frac{8\pi}{3}    
= \Will(F) + \frac{8\pi}{3} \geq 2\pi^2 + \frac{8\pi}{3} 
> 8\pi \geq \Will(F_0) \geq 
\limsup_{j\to \infty} \Will(F_{t_{j}}),  \nonumber
\end{eqnarray}
obtaining exactly condition (2.91) in \cite{Schaetzle.Conf.factor.2013} for $n=4$ and  
$e_4:= \frac{8\pi}{3}$ - as introduced in equation (1.2) 
of \cite{Schaetzle.Conf.factor.2013} - and we may therefore 
apply Proposition 2.4 in \cite{Schaetzle.Conf.factor.2013}
\footnote{Proposition 2.4 in \cite{Schaetzle.Conf.factor.2013} builds directly on the proof of Proposition 2.1 in \cite{Schaetzle.Conf.factor.2013}, i.e. on Lemma 2.1 in \cite{Kuwert.Schaetzle.2012}, on Theorem 3.1 in \cite{Simon.1993}, and on the fundamental results in 
Theorems 3.1 and 6.1 of \cite{Kuwert.Schaetzle.2012}. 
It is worth mentioning here that the latter theorem 	together with its improved versions, 
Theorem 5.1 and Proposition 5.1 
in \cite{Schaetzle.Conf.factor.2013}, will again play 
a central role in the proof of the second part of 
Theorem \ref{singular.time.MIWF.Hopf.tori} below.}.
Hence, recalling here the conformality relations in \eqref{conformal.reparametrization} of the 
embeddings $F_{t_{j}}$ we conclude from Proposition 2.4 in \cite{Schaetzle.Conf.factor.2013} 
that the resulting conformal factors $u_{j} \in C^{\infty}(\Sigma)$ of the pullback metrics 
$(F_{t_{j}})^*(g_{\textnormal{euc}})$
in \eqref{conformal.reparametrization}
are uniformly bounded, i.e. that there holds: 
\begin{equation}   \label{bounded.factor} 
(F_{t_{j}})^*(g_{\textnormal{euc}}) 
= e^{2u_{j}} \,g_{\textnormal{poin},j} \quad 
\textnormal{with} \quad 
\parallel u_{j} \parallel_{L^{\infty}(\Sigma)}\leq\Lambda, 
\quad \forall\, j \in \nat,  
\end{equation}
just as desired. Here, the constant $\Lambda$ does not only depend on measure theoretic properties of the varifold 
$\mu$ and on the limit of the growing gaps 
$\Delta_j:=8 \pi-\Will(F_j)$, 
but also on the sizes of local integrals of the squared fundamental forms $|A_{F_{t_{j_l}}}|^2$ in $\rel^4$   
of the embedded surfaces $F_{t_{j_l}}(\Sigma)\subset \sphere^3$ satisfying \eqref{mu.in.the.limit} and \eqref{Hausdorff.converge}, as $l \to \infty$; see here formulae (2.20)--(2.27) in the proof of Proposition 2.1 in \cite{Schaetzle.Conf.factor.2013}, the proof of Proposition 2.4 and Theorem 5.2 in \cite{Schaetzle.Conf.factor.2013}. 
Having obtained estimate \eqref{bounded.factor}
we can dive into the proof of Proposition 6.1 in \cite{Schaetzle.Conf.factor.2013}. 
As in that proof we infer from \eqref{monotonic.Will} 
and \eqref{bounded.factor} and from Lemma 5.2 in \cite{Kuwert.Schaetzle.2012} that the conformal structures 
$S(F_{t_{j}}(\Sigma))$ corresponding to the conformal classes of the metrics $(F_{t_{j}})^*(g_{\textnormal{euc}})$     
respectively $g_{\textnormal{poin},j}$ in \eqref{conformal.reparametrization} are compactly contained 
in the moduli space $\Mod_1$
\footnote{Alternatively one could also apply here 
Theorem 5.2 in \cite{Kuwert.Li.2012} respectively Theorem 1.1 in \cite{Riviere.Lip.conf.imm.2013}.}.
Hence, up to extraction of a subsequence - which we shall relabel again - there exist diffeomorphisms 
$\Phi_j:\Sigma \stackrel{\cong}\longrightarrow \Sigma$ 
such that 
\begin{equation} \label{smooth.converg.metrics} 
\Phi_{j}^*g_{\textnormal{poin},j}
\longrightarrow g_{\textnormal{poin}} \quad 
\textnormal{smoothly as} \,\,\,j\to \infty
\end{equation}  
for some zero scalar curvature and unit volume metric $g_{\textnormal{poin}}$ on $\Sigma$. In view of 
statements \eqref{bounded.factor} and \eqref{smooth.converg.metrics}, in view of the desired 
assertion of this part of the theorem and since the 
Willmore functional is invariant with respect to smooth 
reparametrization, we shall replace the embeddings 
$F_{t_{j}}$ by their reparametrizations 
$\tilde F_{t_{j}}:= F_{t_{j}} \circ \Phi_j$, 
for every $j\in \nat$, and continue this proof  
using the modified equations
\begin{equation}   \label{bounded.factor.mod}
(\tilde F_{t_{j}})^*(g_{\textnormal{euc}}) 
= e^{2u_{j}\circ \Phi_j} \,\Phi_j^*g_{\textnormal{poin},j} \quad \textnormal{with} \quad 
\parallel u_{j}\circ \Phi_j \parallel_{L^{\infty}(\Sigma)} \leq \Lambda, \quad \forall\, j \in \nat, 
\end{equation} 	
instead of equations \eqref{bounded.factor}.
For ease of notation we shall rename 
$\Phi_j^*g_{\textnormal{poin},j}$ again into 
$g_{\textnormal{poin},j}$ and $u_{j}\circ \Phi_j$ 
into $u_{j}$, $\forall \,j\in \nat$, so that we can assume 
the convergence 
\begin{equation}  \label{smooth.converg.metrics.mod} 
g_{\textnormal{poin},j}
\longrightarrow g_{\textnormal{poin}} \quad 
\textnormal{smoothly as} \,\,\,j\to \infty, 
\end{equation} 
instead of the convergence in \eqref{smooth.converg.metrics}. 
Moreover, from equations \eqref{bounded.factor.mod}, together with formulae \eqref{Liouville} and \eqref{Gauss.Bonnet.K} above we infer, just as in formula (6.4) of \cite{Schaetzle.Conf.factor.2013}: 
$$
-\triangle_{\tilde F_{t_{j}}^*(g_{\textnormal{euc}})}(u_{j})    
= K_{\tilde F_{t_{j}}^*(g_{\textnormal{euc}})} 
\quad \textnormal{on} \,\,\, \Sigma,
$$
or expressed equivalently: 
\begin{equation} \label{Gauss.curvature}
-\triangle_{g_{\textnormal{poin},j}}(u_{j})    
= e^{2 u_{j}} \, K_{\tilde F_{t_{j}}^*(g_{\textnormal{euc}})} 
\quad \textnormal{on} \,\,\, \Sigma.
\end{equation}
Following now exactly the lines of the proof of 
Proposition 6.1 in \cite{Schaetzle.Conf.factor.2013}, 
the equations in \eqref{Gauss.curvature} yield together 
with the differential equations
\begin{equation} \label{mean.curvature} 
\triangle^{\rel^4}_{g_{\textnormal{poin},j}}(\tilde F_{t_{j}})  = e^{2 u_{j}} \, 
\vec H_{\tilde F_{t_{j}},\rel^4}   \quad 
\textnormal{on} \,\,\, \Sigma,
\end{equation}  
and together with statements \eqref{monotonic.Will}, \eqref{bounded.factor.mod} and \eqref{smooth.converg.metrics.mod}, and also using the fact that each $\tilde F_{t_{j}}$ maps $\Sigma$ into the compact $3$-sphere, the estimates: 
\begin{eqnarray}  \label{valuable.bounds} 
\parallel \nabla u_{j} \parallel_{L^{2}(\Sigma,g_{\textnormal{poin}})} 
\leq C(\Lambda)  \quad \textnormal{and} \quad  
\parallel \tilde F_{t_{j}} \parallel_{W^{2,2}(\Sigma,g_{\textnormal{poin}})} 
\leq C(\Lambda)   	
\end{eqnarray} 
for every $j\in \nat$. Hence, we obtain convergent 
subsequences $\{u_{j_k}\}$ and $\{\tilde F_{t_{j_k}}\}$ of 
$\{u_{j}\}$ and $\{\tilde F_{t_{j}}\}$: 
\begin{eqnarray}
	u_{j_k} \longrightarrow u \qquad \textnormal{weakly in} \,\, W^{1,2}(\Sigma,g_{\textnormal{poin}})     \\ 
	u_{j_k} \longrightarrow u \qquad \textnormal{weakly* in} \,\, L^{\infty}(\Sigma,g_{\textnormal{poin}})     \label{weak.convergence}     \\
	u_{j_k} \longrightarrow u \qquad 
	\textnormal{pointwise a.e. in} \,\, \Sigma              \label{pointwise.convergence}     \\
	\textnormal{and} \quad 	
	\tilde F_{t_{j_k}}\longrightarrow f \qquad 
	\textnormal{weakly in} \,\, W^{2,2}(\Sigma,g_{\textnormal{poin}}) 
	\label{weak.convergence.2}
\end{eqnarray} 
as $k\to \infty$, for appropriate functions 
$u \in W^{1,2}(\Sigma,g_{\textnormal{poin}}) \cap 
L^{\infty}(\Sigma,g_{\textnormal{poin}})$ and 
$f\in W^{2,2}(\Sigma,g_{\textnormal{poin}})$. 
Moreover, we infer from \eqref{bounded.factor.mod} and 
\eqref{smooth.converg.metrics.mod} that 
$\nabla^{g_{\textnormal{poin}}} \tilde F_{t_{j}}$ is 
uniformly bounded in $L^{\infty}(\Sigma,g_{\textnormal{poin}})$. 
Hence, again recalling the fact that $\tilde F_{t_{j}}$ 
map $\Sigma$ into the $3$-sphere, we obtain:   
\begin{equation}   \label{bounded.factor.2} 
\parallel \tilde F_{t_{j}} 
\parallel_{W^{1,\infty}(\Sigma,g_{\textnormal{poin}})}
\leq \textnormal{Const}(\Lambda).        
\end{equation}
We thus infer from Theorem 8.5 in \cite{Alt.2015} 
and from Theorem 4.12 in \cite{Alt.2015}, i.e.  
from Arzela's and Ascoli's Theorem,  
that the convergent subsequence $\{\tilde F_{t_{j_k}}\}$ in 
\eqref{weak.convergence.2} also converges in the senses: 
\begin{eqnarray}   \label{W.1.infty.convergence} 
\tilde F_{t_{j_k}} \longrightarrow f  
\qquad \textnormal{weakly* in} \,\,\, W^{1,\infty}(\Sigma,g_{\textnormal{poin}})       \\ 
\label{C0.convergence} \textnormal{and} \quad 
\tilde F_{t_{j_k}} \longrightarrow f  
\qquad \textnormal{in} \,\,\,C^{0}(\Sigma,g_{\textnormal{poin}})  
\end{eqnarray}   
as $k\to \infty$. Now we can conclude from \eqref{bounded.factor.mod} and from the above convergences \eqref{smooth.converg.metrics.mod}, \eqref{weak.convergence}, \eqref{pointwise.convergence} and
\eqref{W.1.infty.convergence}, that in the limit there holds
actually:
\begin{equation}   \label{conformal.limit} 
f^*(g_{\textnormal{euc}}) = e^{2u} \,g_{\textnormal{poin}}, 
\end{equation}  
showing that $f$ is a uniformly conformal 
$(W^{2,2}\cap W^{1,\infty})$-immersion with respect to 
$g_{\textnormal{poin}}$ on $\Sigma$, with 
$\parallel u \parallel_{L^{\infty}(\Sigma)} \leq \Lambda$ 
for the same constant $\Lambda$ as in \eqref{bounded.factor}, 
depending on the sequence $\{F_{t_{j_l}}\}$ from \eqref{weak.convergence.mu} and on the limit varifold $\mu$. 
Moreover, we obtain as in the proof of Proposition 2.1 in \cite{Schaetzle.Conf.factor.2013}, line (2.30),
that $\mu_f:=f(\mu_{f^*(g_{\textnormal{euc}})})$ coincides  
with the varifold $\mu$ in \eqref{mu.in.the.limit}. 
Similarly to the argument above before 
estimate \eqref{Abschaetz.nerd}, we now continue as in the proof of Proposition 2.1 in \cite{Schaetzle.Conf.factor.2013} and conclude from the facts, that $f\in W^{2,2}(\Sigma,\rel^4)$
and that $\Will(f) \leq 
\liminf_{k\to \infty} \Will(\tilde F_{t_{j_k}})<8 \pi$ and from equation \eqref{conformal.limit}  via Proposition 7.2 in \cite{Schaetzle.Conf.factor.2013} that $f$ is injective, and that not only $\mu_f=\mu$ holds, but also 
$\mu_f = \Hn\lfloor_{f(\Sigma)}$, and moreover that   
$$
\textnormal{spt}(\mu_f) = f(\Sigma), \,\,\, 
\textnormal{thus also that} 
\,\,\, \textnormal{spt}(\mu)= f(\Sigma),
$$
and that therefore 
\begin{equation}  \label{f.Lipschitz.homeo} 
f: \Sigma \stackrel{\cong}\longrightarrow \textnormal{spt}(\mu)
\end{equation} 
is a bi-Lipschitz continuous homeomorphism. Hence, 
assertions \eqref{parametrization.mu} and \eqref{spt.mu.f} 
are already proven. Continuing here as in the proof of Proposition 7.2 of \cite{Schaetzle.Conf.factor.2013}, we infer that the coinciding integral varifolds $\mu$, $\mu_f$ have weak mean curvature vectors $\vec H_{\mu}=\vec H_{\mu_f}$ 
in $L^2(\mu_f)$ and that formula 
\eqref{Willmore.energy} holds here.         
\item[3)] If we suppose again that the support of the 
limit varifold $\mu$ is an embedded torus and that additionally its Willmore energy $\Will(\mu)$ 
from \eqref{Willmore.energy} equals the limit 
of the Willmore energies 
$\Will(F_{t_{j_l}})$ of the embeddings 
$F_{t_{j_l}}$ appearing in \eqref{weak.convergence.mu} and 
\eqref{Hausdorff.converg}, 
then we can follow the argument at the end of the proof of Proposition 5.3 in \cite{Kuwert.Schaetzle.conf.class.2013}.
First of all, similarly to the proof of the second part 
of the theorem, we ease our notation and relabel the given subsequence $\{F_{t_{j_l}}\}$ again into $\{F_{t_{j}}\}$,
and then we replace them by appropriate reparametrizations 
$\tilde F_{t_j}:=F_{t_{j}}\circ \Phi_j$ in 
order to have statements \eqref{bounded.factor.mod}--\eqref{valuable.bounds} and \eqref{bounded.factor.2} 
at our disposal. Hence, we obtain in the 
weak limits \eqref{weak.convergence.2} and 
\eqref{W.1.infty.convergence} a homeomorphic parametrization 
$f\in (W^{2,2}\cap W^{1,\infty})(\Sigma,g_{\textnormal{poin}})$ 
of $\textnormal{spt}(\mu)$ which is a uniformly conformal 
$(W^{2,2}\cap W^{1,\infty})$-immersion with respect to 
$g_{\textnormal{poin}}$ on account of formula \eqref{conformal.limit}. This formula 
implies particularly: 
\begin{equation} \label{mean.curvature.limit} 
\triangle^{\rel^4}_{g_{\textnormal{poin}}}(f)    
= e^{2 u} \, \vec H_{f,\rel^4}   \,\,\, \textnormal{on} \,\,\, \Sigma,
\end{equation} 
similarly to equations \eqref{mean.curvature}.  
Now, from our additional requirement that 
$\Will(\mu) = \lim_{j \to \infty} \Will(F_{t_j})$, 
and from statements \eqref{Willmore.energy}, 
\eqref{bounded.factor.mod}, \eqref{conformal.limit} and 
\eqref{mean.curvature.limit} we infer that 
\begin{eqnarray}  \label{converg.L2.norms}
\int_{\Sigma} |\triangle^{\rel^4}_{g_{\textnormal{poin},j}}(\tilde F_{t_{j}})|^2 
\,e^{-2u_{j}} \,  d\mu_{g_{\textnormal{poin,j}}}   
= \int_{\Sigma} |\vec H_{\tilde F_{t_{j}}}|^2 \,e^{2u_{j}} \, d\mu_{g_{\textnormal{poin,j}}}             
= \int_{\Sigma} |\vec H_{\tilde F_{t_{j}}}|^2 \,
d\mu_{\tilde F_{t_{j}}^*{g_{\textnormal{euc}}}}  \quad  \\
\longrightarrow 4\, \Will(\mu) =
\int_{\Sigma} |\vec H_{f}|^2 \,d\mu_{f^*{g_{\textnormal{euc}}}} 
= \int_{\Sigma} |\vec H_{f}|^2 \,e^{2u} \, d\mu_{g_{\textnormal{poin}}}                  
= \int_{\Sigma} |\triangle^{\rel^4}_{g_{\textnormal{poin}}}(f)|^2 
\,e^{-2u} \,  d\mu_{g_{\textnormal{poin}}}  \nonumber    
\end{eqnarray} 
as $j\to \infty$. 
Now we observe that estimates \eqref{bounded.factor.mod} 
and \eqref{valuable.bounds} imply the uniform bound 
\begin{equation}   \label{bounded.laplace} 
\parallel \triangle^{\rel^4}_{g_{\textnormal{poin},j}}
(\tilde F_{t_{j}}) \, e^{-u_{j}}  \parallel_{L^2(\Sigma,g_{\textnormal{poin}})} 
\leq C(\Lambda) 
\end{equation}  
for every $j\in \nat$. Hence, combining convergence  \eqref{converg.L2.norms} with estimate \eqref{bounded.laplace} and with the smooth convergence of metrics in \eqref{smooth.converg.metrics.mod} 
we obtain elementarily, that also 
\begin{eqnarray}  \label{convergence.L2.norms}
\int_{\Sigma} |\triangle^{\rel^4}_{g_{\textnormal{poin},j}}(\tilde F_{t_{j}})|^2 
\,e^{-2u_{j}} \,  d\mu_{g_{\textnormal{poin}}} 
\longrightarrow 
\int_{\Sigma} |\triangle^{\rel^4}_{g_{\textnormal{poin}}}(f)|^2 
\,e^{-2u} \,  d\mu_{g_{\textnormal{poin}}}    
\end{eqnarray} 
as $j\to \infty$. In order to finally prove our assertion in   
\eqref{W.2.2.convergence}, we shall rather continue 
working with the subsequence $\{t_{j_k}\}$ of $\{t_j\}$
appearing in convergences \eqref{weak.convergence}--\eqref{weak.convergence.2} 
instead of $\{t_j\}$ itself. 
One can combine estimate \eqref{bounded.laplace}
with convergences \eqref{pointwise.convergence} and \eqref{weak.convergence.2} and estimate \eqref{bounded.factor.mod} in order to prove that 
\begin{equation}  \label{weak.convergence.in.L2.true} 
\triangle^{\rel^4}_{g_{\textnormal{poin},j_k}}
(\tilde F_{t_{j_k}}) \, e^{-u_{j_k}} \longrightarrow \triangle^{\rel^4}_{g_{\textnormal{poin}}}(f) \, e^{-u} 
\quad  \textnormal{weakly in} \,\,\,L^{2}(\Sigma,g_{\textnormal{poin}}) 
\end{equation}  
as $k\to \infty$. Combining this with convergence 
\eqref{convergence.L2.norms} we finally obtain: 
$$ 
\triangle^{\rel^4}_{g_{\textnormal{poin},j_k}}(\tilde F_{t_{j_k}}) \, e^{- u_{j_k}} \longrightarrow \triangle^{\rel^4}_{g_{\textnormal{poin}}}(f) \, e^{-u} \quad  \textnormal{strongly in} \,\,\,L^{2}(\Sigma,g_{\textnormal{poin}}) 
$$  
for the subsequence appearing in 
\eqref{weak.convergence.in.L2.true}. 
Combining this again with estimate \eqref{bounded.factor.mod}, with convergence \eqref{pointwise.convergence} and with
Vitali's convergence theorem, we finally obtain:  
\begin{equation}  \label{strong.convergence.in.L2} 
\triangle^{\rel^4}_{g_{\textnormal{poin},j_k}}(\tilde F_{t_{j_k}}) \longrightarrow \triangle^{\rel^4}_{g_{\textnormal{poin}}}(f)  
\quad  \textnormal{strongly in} \,\,\,L^{2}(\Sigma,g_{\textnormal{poin}}) 
\end{equation}  
as $k\to \infty$. Just as in the end of the 
proof of Proposition 5.3 in \cite{Kuwert.Schaetzle.conf.class.2013}, p. 506, we can 
infer now from convergence \eqref{strong.convergence.in.L2} and again from \eqref{smooth.converg.metrics.mod} and
\eqref{valuable.bounds} that 
\begin{eqnarray*} 
\partial_u\Big{(} g^{uv}_{\textnormal{poin}}
\sqrt{\det g_{\textnormal{poin}}} \,\,
\partial_v(\tilde F_{t_{j_k}}-f) \Big{)}              
\longrightarrow 0    \quad  
\textnormal{strongly in} \,\,\,L^{2}(\Sigma,g_{\textnormal{poin}}),  
\end{eqnarray*} 	
as $k\to \infty$. Hence, by standard $L^2$-estimates for 
uniformly elliptic partial differential equations on 
$(\Sigma,g_{\textnormal{poin}})$ together with   	
convergence \eqref{C0.convergence}, we conclude
that the subsequence $\{\tilde F_{t_{j_k}}\}$ in 
\eqref{strong.convergence.in.L2} of the reparametrized 
sequence of embeddings $\{\tilde F_{t_{j}}\}$ 
converges strongly to $f$ in 
$W^{2,2}(\Sigma,g_{\textnormal{poin}})$, just as asserted  
in \eqref{W.2.2.convergence}.\\
\noindent 
Now we finally prove that the subsequence  
$\{\tilde F_{t_{j_k}}\}$ from \eqref{W.2.2.convergence} 
is a sequence of uniformly bi-Lipschitz continuous 
homeomorphisms between  
$(\Sigma,g_{\textnormal{poin}})$ and their images in 
$\sphere^3$. On account of estimate \eqref{bounded.factor.2} the sequence $\{\tilde F_{t_{j_k}}\}$ is certainly uniformly Lipschitz continuous. Hence, if the asserted uniform bi-Lipschitz property of 
$\{\tilde F_{t_{j_k}}\}$ would not hold here, then there existed another subsequence 
$f_m:=\tilde F_{t_{j_{k_m}}}$ of $\tilde F_{t_{j_k}}$ and certain points $p_m, q_m \in \Sigma$ such that 
\begin{equation} \label{f.k.contradiction}      
0\leq |f_m(p_m)-f_m(q_m)| <\frac{1}{m}\,
\textnormal{dist}_{(\Sigma,g_{\textnormal{poin}})}(p_m,q_m)   
\qquad \forall \, m\in \nat
\end{equation}
would have to hold. On account of the compactness of 
$(\Sigma,g_{\textnormal{poin}})$, we first 
extract convergent subsequences  
$p_{m'} \longrightarrow p$ and $q_{m'} \longrightarrow q$ in 
$(\Sigma,g_{\textnormal{poin}})$, which we relabel again. Inserting this into  
\eqref{f.k.contradiction} and using the uniform convergence 
of $\{f_m\}$ to $f$ from \eqref{C0.convergence}, we obtain
in the limit as $m \to \infty$: $|f(p)-f(q)|=0$ and 
thus that $p=q$ by injectivity of $f$. Hence, similarly 
to the proof of Proposition 7.2 in 
\cite{Schaetzle.Conf.factor.2013} we now introduce local 
conformal coordinates about this limit point $p$ with respect to 
the zero scalar curvature metrics $g_{\textnormal{poin},j_{k_m}}$. 
By Definitions 2.3.1, 2.3.3 and 2.3.4 in \cite{Jost.2006}
this means, that we consider 
open neighbourhoods $U_m(p)$ of $p$ in $\Sigma$ and charts
$\varphi_m:B^2_1(0) \stackrel{\cong}\longrightarrow U_m(p)$ 
mapping $0$ to $p$, such that 
\begin{eqnarray} \label{phi.k} 
\varphi_m^*g_{\textnormal{poin},j_{k_m}}= 
e^{2v_m} \,g_{\textnormal{euc}} \,\,\, \textnormal{with} \,\,\,
\triangle_{g_{\textnormal{euc}}}(v_m) = 0 \,\,\,
\textnormal{on} \,\,\,B^2_1(0), \nonumber \\ 
\textnormal{and with} \,\,\, v_m \in L^{\infty}(B^2_1(0))      
\,\,\, \textnormal{and} \,\,\,
\parallel v_m \parallel_{L^{\infty}(B^2_1(0))} \leq 
C(g_{\textnormal{poin}}) \quad \forall \, m\in \nat, \quad
\end{eqnarray}
where we have used convergence \eqref{smooth.converg.metrics.mod}.
Using Cauchy estimates we thus also have: 
\begin{equation}  \label{phi.k.2} 
\parallel \nabla^s(v_m) \parallel_{L^{\infty}(B^2_{3/4}(0))} 
\leq C(g_{\textnormal{poin}},s) \quad \forall \, m\in \nat
\end{equation} 
and for each fixed $s\in \nat$. Hence, we obtain:
\begin{eqnarray*}  
(f_m \circ \varphi_m)^*(g_{\textnormal{euc}})
= \varphi_m^*(e^{2u_{j_{k_m}}} \,g_{\textnormal{poin},j_{k_m}})= \\ 
= e^{2u_{j_{k_m}}\circ \varphi_m}\, \varphi_m^*g_{\textnormal{poin},j_{k_m}} 
= e^{2u_{j_{k_m}}\circ \varphi_m+2v_m}\, g_{\textnormal{euc}}
\end{eqnarray*}  
showing that $f_m \circ \varphi_m:B^2_1(0) \longrightarrow \sphere^3$ are smooth conformal embeddings with respect to the Euclidean metric on the unit disc $B^2_1(0)$. We set 
$$
M:=\sup_{m\in \nat}\big{(} 
\parallel u_{j_{k_m}} \parallel_{L^{\infty}(\Sigma)} 
+ \parallel v_m \parallel_{L^{\infty}(B^2_1(0))}\big{)}<\infty,
$$   
and on account of the strong $W^{2,2}$-convergence 
\eqref{W.2.2.convergence} of the embeddings $f_m$ and 
on account of estimates \eqref{phi.k} and \eqref{phi.k.2} we 
may choose $\varrho\in (0,1)$ that small, such that 
\begin{equation}  \label{little.integral} 
\int_{B^2_{\varrho}(0)} |D^2(f_m \circ \varphi_m)|^2 \,d\Lno
< \frac{\pi}{2} \, \tanh(\pi) \, e^{-6M},  \quad \forall \,
m\in \nat.  
\end{equation}  
On account of $p_m \rightarrow p=q \leftarrow q_m$ in 
$(\Sigma,g_{\textnormal{poin}})$ and due to estimate 
\eqref{phi.k} we know that  
$p_m, q_m \in \varphi_m(B^2_{\varrho}(0))\subset \subset U_m(p)$ for sufficiently large $m$, 
and thus we obtain from equation \eqref{bounded.factor.mod}, 
convergence \eqref{smooth.converg.metrics.mod}, estimate \eqref{little.integral} and 
from Lemmata 4.2.7 and 4.2.8 in \cite{Mueller.Sverak.1995}, similarly to the end of the proof of Proposition 7.2 in \cite{Schaetzle.Conf.factor.2013}: 
$$     
\textnormal{dist}_{(\Sigma,g_{\textnormal{poin}})}(p_m,q_m)  
\leq 2\,e^{2M} \, \textnormal{dist}_{(\Sigma,f_m^*g_{\textnormal{euc}})}(p_m,q_m)
\leq 2\,\sqrt{2}\,\, e^{2M} \, |f_m(p_m)-f_m(q_m)|
$$
for sufficiently large $m$, which contradicts \eqref{f.k.contradiction} for very large $m$. 
Hence, the sequence of embeddings 
$\tilde F_{t_{j_k}}:\Sigma \stackrel{\cong}
\longrightarrow \tilde F_{t_{j_k}}(\Sigma) \subset \sphere^3$
is indeed uniformly bi-Lipschitz continuous. \\
Now we choose some point $x\in \sphere^3$ arbitrarily.
On account of the uniform bi-Lipschitz property of the 
embeddings $\{\tilde F_{t_{j_k}}\}$, we can find for any 
small $r>0$ some sufficiently small $\eta>0$, depending 
on $x$ and $r$, such that 
the preimages $\tilde F_{t_{j_k}}^{-1}(B^4_{\eta}(x)
\cap \sphere^3)$ are contained in open geodesic discs $B^{g_{\textnormal{poin}}}_{r}(p_k)$ about certain points 
$p_k$ in $\Sigma$ with respect to the metric $g_{\textnormal{poin}}$, 
for every $k\in \nat$. On account of the compactness of $(\Sigma,g_{\textnormal{poin}})$ we can extract some 
convergent subsequence $p_{k_m}\to p$ in $\Sigma$,
depending on $x\in \sphere^3$.  
Setting now again $f_m:=\tilde F_{t_{j_{k_m}}}$,     	
we thus obtain the existence of some large 
$K=K(x) \in \nat$, such that $f_m^{-1}(B^4_{\eta}(x)\cap \sphere^3)$ is contained in $B^{g_{\textnormal{poin}}}_{2r}(p)$ 
for every $m \geq K$. Combining this again with the strong $W^{2,2}$-convergence \eqref{W.2.2.convergence} of the 
embeddings $f_m$, we infer that for every $\varepsilon>0$ 
there is some sufficiently small $\eta>0$, depending 
on $x$ and $\varepsilon$, such that: 
$$
	\int_{(f_{m})^{-1}(B^4_{\eta}(x)\cap \sphere^3)} 
	|A_{f_m}|^2 \,d\mu_{f_m^*g_{\textnormal{euc}}}
	\leq \int_{B^{g_{\textnormal{poin}}}_{2r}(p)} 
	|D^2f_m|^2 \,d\mu_{f_m^*g_{\textnormal{euc}}} 
	<\varepsilon, \quad \forall \,m \geq K,   
$$ 
which has already proved assertion \eqref{no.A.concentration}. 
This assertion automatically implies  
the last assertion of the third part of this theorem, 
because if the measures $\Mill_l$ from 
\eqref{no.concentration.Mill} would concentrate at 
some point $x\in \sphere^3$, then there was some 
sufficiently small $\varepsilon >0$ and some 
``bad subsequence'' $\{F_{t_{j_{l_b}}}\}$ 
of the original sequence
$\{F_{t_{j_l}}\}$ from \eqref{weak.convergence.mu} 
such that   
\begin{eqnarray} \label{contradiction.concentration} 
\sup \{\rho >0 \,\big{|}\, \Mill_{l_b}(B^4_{\rho}(x))<\varepsilon \, \} 
= \sup \Big{\{}\rho >0 \,\big{|}\, 
\int_{F_{t_{j_{l_b}}}^{-1}(B^4_{\rho}(x)\cap \sphere^3)} 
|A_{F_{t_{j_{l_b}}}}|^2\,
d\mu_{F_{t_{j_{l_b}}}^*g_{\textnormal{euc}}}<\varepsilon
\, \Big{\}} \nonumber \\
\longrightarrow 0 \quad \textnormal{as}\quad b \to \infty.
\qquad 
\end{eqnarray} 
But we can exchange the original weakly convergent sequence $\{F_{t_{j_l}}\}$ from \eqref{weak.convergence.mu} by any of its subsequences, e.g. by $\{F_{t_{j_{l_b}}}\}$, and the entire above argument yielding the strong $W^{2,2}$-convergence \eqref{W.2.2.convergence} again applies to $\{F_{t_{j_{l_b}}}\}$ and thus guarantees the existence of some subsequence 
$\{F_{t_{j_{l_{b_k}}}} \}$ and of smooth diffeomorphisms  
$\Theta_k:\Sigma \stackrel{\cong}\longrightarrow \Sigma$ such 
that $\tilde F_{t_{j_{l_{b_k}}}}:=F_{t_{j_{l_{b_k}}}}\circ \Theta_k$ converge as in \eqref{W.2.2.convergence} 
to the same bi-Lipschitz-parametrization   
$f:\Sigma \stackrel{\cong}\longrightarrow \textnormal{spt}(\mu)$. 
Moreover, each immersion $\tilde F_{t_{j_{l_{b_k}}}}$ 
has to be a uniformly bi-Lipschitz homeomorphism of 
$(\Sigma,g_{\textnormal{poin}})$ onto its image in 
$(\sphere^3,g_{\textnormal{euc}})$, just on account of 
the above reasoning. Hence, adopting the above
argument leading finally to assertion \eqref{no.A.concentration} we again obtain some 
further subsequence 
$\{f_n\}$ of $\tilde F_{t_{j_{l_{b_k}}}}$       	
and some large $K=K(x) \in \nat$ such that for every $\varepsilon>0$ there is some sufficiently small $\eta>0$, depending on $x$ and $\varepsilon$, such that: 
$$
\int_{(f_{n})^{-1}(B^4_{\eta}(x)\cap \sphere^3)} 
|A_{f_n}|^2 \,d\mu_{f_n^*g_{\textnormal{euc}}}
<\varepsilon, \quad \forall \,n \geq K.   
$$ 
But this result shows the existence of some particular subsequence of the sequence  
$\{\sup \{\rho >0 \,\big{|}\,\Mill_{l_b}(B_{\rho}(x))
<\varepsilon \,\}\}_{b\in \nat}$
appearing in \eqref{contradiction.concentration}   
which remains bounded from below by $\eta=\eta(x,\varepsilon)>0$. Hence, hypothesis \eqref{contradiction.concentration} 
indeed turns out to be wrong, and the measures $\Mill_l$ from 
\eqref{no.concentration.Mill} cannot concentrate at 
any $x\in \sphere^3$ respectively $x\in \rel^4$.       
\item[4)] As in the third part of this theorem we 
consider a limit varifold $\mu$ in \eqref{weak.convergence.mu}, 
whose support is a compact and embedded torus in 
$\sphere^3$, we assume again that 
$\lim_{l\to \infty} \Will(F_{t_{j_l}})=\Will(\mu)$ 
holds for the sequence $\{F_{t_{j_l}}\}$ from \eqref{weak.convergence.mu}, and   
additionally we require here: 
\begin{eqnarray}  \label{boendae} 
\parallel |A^0_{F_{t_{j_l}},\sphere^3}|^2 \parallel_{L^{\infty}(\Sigma)} \leq K \quad 
\textnormal{and} \quad \big{|}\frac{d}{dt}\Will(F_{t})\big{|}\lfloor_{t=t_{j_l}} 
\leq K, \quad \textnormal{for all}\,\, k\in \nat,
\end{eqnarray}
for some large $K>1$.  
We know already from \eqref{W.2.2.convergence}, \eqref{W.1.infty.convergence} and \eqref{C0.convergence}, 
that some particular subsequence $\{F_{t_{j_k}}\}$ 
of the sequence $\{F_{t_{j_l}}\}$ can be reparametrized into embeddings $\tilde F_{t_{j_k}}$ which converge strongly 
in $W^{2,2}(\Sigma,g_{\textnormal{poin}})$, weakly* 
in $W^{1,\infty}(\Sigma,g_{\textnormal{poin}})$ 
and also uniformly:  
\begin{equation}  \label{C0.convergence.2}
\tilde F_{t_{j_k}} \longrightarrow f  
\qquad \textnormal{in} \,\,\,C^{0}(\Sigma,g_{\textnormal{poin}}),  
\end{equation}
as $k\to \infty$, to the uniformly conformal 
$(W^{2,2}\cap W^{1,\infty})$-parametrization $f$ of $\textnormal{spt}(\mu)$ from 
\eqref{f.Lipschitz.homeo}. Since the MIWF \eqref{Moebius.flow} is conformally invariant, we can 
therefore assume that the images $F_{t_{j_k}}(\Sigma)$ 
of the entire subsequence $\{F_{t_{j_k}}\}$ satisfy:
\begin{equation}  \label{containoe}
F_{t_{j_k}}(\Sigma) \subset 
\sphere^3 \setminus B^4_{\delta}((0,0,0,1)),	\quad 
\forall \,k\in \nat,
\end{equation}
for some sufficiently small $\delta>0$.    
We therefore project the entire flow line $\{F_{t}\}$ 
of the MIWF stereographically from 
$\sphere^3 \setminus \{(0,0,0,1)\}$ into $\rel^3$, 
and we shall consider henceforth the sequence of 
stereographically projected embeddings 
$\Stereo \circ F_{t_{j_k}}:\Sigma \longrightarrow \rel^3$ instead of the original sequence $\{F_{t_{j_k}}\}$. 
Using the facts that $\{F_{t}\}$ solves flow equation \eqref{Moebius.flow} and that both the MIWF and 
the Willmore functional itself are 
\underline{conformally invariant}, 
we can compute by means of the chain 
rule and the additional requirements in \eqref{boendae} 
for the subsequence $\{F_{t_{j_k}}\}$:
\begin{eqnarray*} 
\int_{\Sigma} 
\frac{1}{|A^0_{\Stereo \circ F_{t_{j_k}}}|^4} \,
|\nabla_{L^2}\Will(\Stereo \circ F_{t_{j_k}})|^2 
\,d\mu_{(\Stereo \circ F_{t_{j_k}})^*(g_{\textnormal{euc}})} \qquad\qquad \\
= 2 \,\Big{|}\frac{d}{dt}\Will(\Stereo \circ
F_{t})\Big{|}\lfloor_{t=t_{j_k}} 
=2\, \Big{|}\frac{d}{dt}\Will(F_{t})\Big{|}\lfloor_{t=t_{j_k}} 
\leq 2K, \quad \forall \,k\in \nat. 
\end{eqnarray*} 
Moreover, we notice that by \eqref{containoe} our first 
condition in \eqref{boendae} also implies: 
$$ 
|A^0_{\Stereo \circ F_{t_{j_k}}}|^2 \leq K'(K,\delta)<\infty
\quad \textnormal{on} \,\, \Sigma, \quad \forall \,k\in \nat,
$$ 
and we therefore obtain the estimate:
\begin{equation} \label{first.bound} 
\int_{\Sigma} |\nabla_{L^2}\Will(\Stereo \circ F_{t_{j_k}})|^2 
\,d\mu_{(\Stereo \circ F_{t_{j_k}})^*(g_{\textnormal{euc}})} 
\leq 2K \,(K')^2  \quad  \textnormal{for every} \,\, k \in \nat.    
\end{equation}
Now, we again replace the embeddings 
$\Stereo \circ F_{t_{j_k}}$ by their reparametrizations 
$\tilde f_k:=\Stereo \circ F_{t_{j_k}} \circ \Phi_{j_k}
\equiv \Stereo\circ \tilde F_{t_{j_k}}$ in view of \eqref{bounded.factor.mod} and \eqref{smooth.converg.metrics.mod}, 
and we obtain from \eqref{first.bound}, \eqref{bounded.factor.mod} and \eqref{smooth.converg.metrics.mod} together with the area formula:     
\begin{eqnarray}     \label{bounded.Willmore.operator} 
\parallel \nabla_{L^2}\Will(\tilde f_k) \parallel_{L^2(\Sigma,g_{\textnormal{poin}})}^2 
\equiv \int_{\Sigma} |\nabla_{L^2}\Will(\tilde f_k)|^2 \,  
d\mu_{g_{\textnormal{poin}}}    
\leq 2 \,e^{2\tilde \Lambda} \,
\int_{\Sigma} |\nabla_{L^2}\Will(\tilde f_k)|^2 \,  
d\mu_{\tilde f_k^*(g_{\textnormal{euc}})}        \nonumber      \\
= 2\,e^{2\tilde \Lambda} \,
\int_{\Sigma} |\nabla_{L^2}\Will(\Stereo \circ F_{t_{j_k}})|^2 \, d\mu_{(\Stereo \circ F_{t_{j_k}})^*(g_{\textnormal{euc}})} 
\leq 4K \,(K')^2  \,e^{2\tilde \Lambda} \qquad 
\end{eqnarray}
for every $k \in \nat$, where $\tilde \Lambda$ denotes an upper bound for the new conformal factor $\tilde u_{j_k}$ in 
$L^{\infty}(\Sigma)$, appearing in: 
\begin{equation}  \label{uniform.conformal.R3}
(\tilde f_k)^*(g_{\textnormal{euc}}) 
=e^{2\tilde u_{j_k}} \,g_{\textnormal{poin},j_k} \quad \textnormal{on} \,\, \Sigma
\end{equation} 
for every $k\in \nat$, which holds on account of 
\eqref{bounded.factor.mod} and \eqref{containoe} and due to 
the conformality of the stereographic projection. 
On account of \eqref{bounded.Willmore.operator} 
there is some subsequence of the sequence 
$\{\tilde f_k\}$, which we relabel
into $\{\tilde f_k\}$ again, and some function 
$q\in L^2(\Sigma,g_{\textnormal{poin}})$ with values in 
$\rel^3$, such that
\begin{equation}  \label{weaken}
\nabla_{L^2}\Will(\tilde f_k)  \longrightarrow  q  \quad 
\textnormal{weakly in} \,\,L^2(\Sigma,g_{\textnormal{poin}})
\end{equation}
as $k \to \infty$. Concerning the homeomorphic parametrization
$\tilde f :=\Stereo \circ f$ of the limit torus $\Stereo(\textnormal{spt}(\mu))\subset \rel^3$  
from \eqref{f.Lipschitz.homeo} projected stereographically 
from $\sphere^3\setminus \{(0,0,0,1)\}$ into $\rel^3$, 
we derive from \eqref{conformal.limit}, \eqref{C0.convergence.2} and \eqref{containoe} 
immediately:    
\begin{equation}   \label{conformal.limit.2} 
(\tilde f)^*(g_{\textnormal{euc}})
=(\Stereo \circ f)^*(g_{\textnormal{euc}}) 
= e^{2\tilde u} \,g_{\textnormal{poin}} \quad 
\textnormal{on} \,\,\, \Sigma, 
\end{equation} 
for some function $\tilde u \in L^{\infty}(\Sigma)$, 
which has to additionally satisfy:   
\begin{equation}  \label{pointwise.convergence.2}
\tilde u_{j_k} \longrightarrow \tilde u \quad \textnormal{pointwise \,a.e. \,in}\,\,\, \Sigma, 
\end{equation} 
on account of \eqref{smooth.converg.metrics.mod}, \eqref{uniform.conformal.R3} and 
on account of the strong $W^{2,2}(\Sigma,g_{\textnormal{poin}})$-convergence
of the sequence $\{\tilde f_k\}$ to $\tilde f$ due to 
\eqref{W.2.2.convergence}. Obviously, 
$\tilde f$ is a uniformly conformal 
$(W^{2,2}\cap W^{1,\infty})$-immersion  
with respect to $g_{\textnormal{poin}}$ on account of 
\eqref{conformal.limit.2} - just as $f$ is by \eqref{conformal.limit} - and statements \eqref{weaken} and \eqref{conformal.limit.2} also imply that
\begin{equation}     \label{weakens}
\nabla_{L^2}\Will(\tilde f_k)  \longrightarrow  q \quad \textnormal{weakly in} \,\,
L^2(\Sigma,\tilde f^*{g_{\textnormal{euc}}}),
\end{equation}
as $k \to \infty$, just on account of the definition of 
weak $L^2$-convergence and $\tilde u \in L^{\infty}(\Sigma)$ 
in \eqref{conformal.limit.2}, where   
$q$ is here the same $L^2$-function as in \eqref{weaken}. 
On the other hand, we can immediately 
infer from \eqref{W.2.2.convergence}, \eqref{W.1.infty.convergence}, \eqref{C0.convergence} 
respectively \eqref{C0.convergence.2}, 
combined with \eqref{containoe}, 
that the embeddings 
$\tilde f_k =\Stereo \circ F_{t_{j_k}}\circ \Phi_{j_k}$ converge strongly 
in $W^{2,2}(\Sigma,g_{\textnormal{poin}})$, weakly* 
in $W^{1,\infty}(\Sigma,g_{\textnormal{poin}})$ 
and in $C^{0}(\Sigma,g_{\textnormal{poin}})$ to 
the parametrization $\tilde f =\Stereo \circ f$ 
of the projected torus 
$\Stereo(\textnormal{spt}(\mu))\subset \rel^3$. 
Using particularly the fact that the $\tilde f_k$  
are uniformly bounded in $W^{1,\infty}(\Sigma,g_{\textnormal{poin}})$ 
by \eqref{bounded.factor.2} and \eqref{containoe}, we 
can combine the above mentioned convergences of the sequences $\{\tilde f_k\}$ and $\{\tilde u_{j_k}\}$, equations \eqref{uniform.conformal.R3} and \eqref{conformal.limit.2} and convergence \eqref{smooth.converg.metrics.mod} 
with H\"older's inequality and Vitali's convergence theorem and with formulae \eqref{distributional.Willmore} and \eqref{Riviere.Bernard.Willmore}, in order to conclude that  
\begin{eqnarray}  \label{distributional.Willmore.limit} 
\langle \nabla_{L^2}\Will(\tilde f_k),\varphi 
\rangle_{L^2(\Sigma,\tilde f_k^*(g_{\textnormal{euc}}))} = 	
\langle \nabla_{L^2}\Will(\tilde f_k),\varphi \rangle_{\Dom'(\Sigma)} =	\nonumber \\
= \int_{\Sigma} 
\langle \vec H_{\tilde f_k}, \triangle_{\tilde f_k} \varphi \rangle_{g_{\textnormal{euc}}} 
-  \, g_{\tilde f_k}^{\nu\alpha} \, g_{\tilde f_k}^{\mu\xi}\,\langle 
(A_{\tilde f_k})_{\xi\nu}, \vec H_{\tilde f_k} \rangle_{g_{\textnormal{euc}}} \,\langle \partial_{\mu}\tilde f_k, \partial_{\alpha}\varphi \rangle_{g_{\textnormal{euc}}} \nonumber \\
- \, g_{\tilde f_k}^{\nu\alpha} \, 
g_{\tilde f_k}^{\mu\xi}\,\langle 
(A_{\tilde f_k}^0)_{\xi\nu}, \vec H_{\tilde f_k} \rangle_{g_{\textnormal{euc}}} \, \langle \partial_{\mu}\tilde f_k, \partial_{\alpha}\varphi \rangle_{g_{\textnormal{euc}}}\, d\mu_{g_{\tilde f_k}}   \nonumber  \\
\longrightarrow 
\int_{\Sigma} \Big{(} \langle \vec H_{\tilde f}, \triangle_{\tilde f} \varphi \rangle_{g_{\textnormal{euc}}} 
-  \, g_{\tilde f}^{\nu\alpha} \, g_{\tilde f}^{\mu\xi}\,\langle 
(A_{\tilde f})_{\xi\nu}, \vec H_{\tilde f} \rangle_{g_{\textnormal{euc}}} \,\langle \partial_{\mu}\tilde f, \partial_{\alpha}\varphi \rangle_{g_{\textnormal{euc}}} \nonumber \\
- \, g_{\tilde f}^{\nu\alpha} \, g_{\tilde f}^{\mu\xi}\,\langle 
(A_{\tilde f}^0)_{\xi\nu}, \vec H_{\tilde f} \rangle_{g_{\textnormal{euc}}} \, \langle \partial_{\mu}\tilde f, \partial_{\alpha}\varphi \rangle_{g_{\textnormal{euc}}}\, \Big{)} \,
d\mu_{g_{\tilde f}}    =   \nonumber  \\                 
= \langle \nabla_{L^2}\Will(\tilde f),\varphi \rangle_{\Dom'(\Sigma)} 	 
\end{eqnarray}  
as $k \to \infty$, for every fixed 
$\varphi \in C^{\infty}(\Sigma,\rel^3)$. 
In the above argument one has to recall, that the 
$(W^{2,2}\cap W^{1,\infty})$-immersions 
$\tilde f_k:\Sigma \longrightarrow \rel^3$
and $\tilde f:\Sigma \longrightarrow \rel^3$ map into $\rel^3$, but not into $\sphere^3$,  
such that we can actually apply our special version \eqref{distributional.Willmore} of the distributional 
Willmore operator. Now, combining the weak convergence 
\eqref{weaken} with both the pointwise convergence  \eqref{pointwise.convergence.2} and the uniform 
boundedness of the conformal factors $\tilde u_{j_k}$ 
of the $\tilde f_k$ in \eqref{uniform.conformal.R3}
- on account of \eqref{bounded.factor.mod} and \eqref{containoe} - 
and again with convergence \eqref{smooth.converg.metrics.mod}, we can infer from E8.3 in \cite{Alt.2015} that
$$
\lim_{k\to \infty} 
\langle \nabla_{L^2}\Will(\tilde f_k),\varphi 
\rangle_{L^2(\Sigma,\tilde f_k^*(g_{\textnormal{euc}}))} 
=\lim_{k\to \infty} 
\langle \nabla_{L^2}\Will(\tilde f_k),\varphi 
\rangle_{L^2(\Sigma,\tilde f^*(g_{\textnormal{euc}}))}.
$$    
Combining this equation again with  \eqref{distributional.Willmore.limit} and with 
convergence \eqref{weakens}, then we obtain:
\begin{eqnarray*}    
\langle \nabla_{L^2}\Will(\tilde f),\varphi \rangle_{\Dom'(\Sigma)} 
= \lim_{k\to \infty} 
\langle \nabla_{L^2}\Will(\tilde f_k),\varphi 
\rangle_{L^2(\Sigma,\tilde f_k^*(g_{\textnormal{euc}}))} 
\qquad \\
=\lim_{k\to \infty} 
\langle \nabla_{L^2}\Will(\tilde f_k),\varphi 
\rangle_{L^2(\Sigma,\tilde f^*(g_{\textnormal{euc}}))} 
= \langle\, q, \varphi \,
\rangle_{L^2(\Sigma,\tilde f^*(g_{\textnormal{euc}}))}     \qquad 
\end{eqnarray*}  
$\forall  \,\varphi \in C^{\infty}(\Sigma,\rel^3)$, 
where the function $q$ is of class $L^2((\Sigma,g_{\textnormal{poin}}),\rel^3)$ by \eqref{weakens}. 
This shows that $\nabla_{L^2}\Will(\tilde f)$ is here 
not only a distribution of second order acting on 
$C^{\infty}(\Sigma,\rel^3)$, but it can be identified 
with an $\rel^3$-valued function of class $L^2(\Sigma,g_{\textnormal{poin}})$.  
We can therefore apply here Theorem \ref{Regul.theorem}
to the uniformly conformal 
$(W^{2,2}\cap W^{1,\infty})$-immersion 
$\tilde f$ and conclude that $\tilde f$ is actually of class  $W^{4,2}((\Sigma,g_{\textnormal{poin}}),\rel^3)$, 
whence that $f=\Stereo^{-1} \circ \tilde f$ is of 
class $W^{4,2}((\Sigma,g_{\textnormal{poin}}),\rel^4)$.  	 
\end{itemize}

\section{Dimension-reduction of the M\"obius-invariant 
Willmore flow} 
\label{dimension.reduction} 
The basic ingredient of this approach to the proofs of 
Theorems \ref{singular.time.MIWF.Hopf.tori} and 
\ref{limit.at.infinity} is the Hopf-fibration
$$
\sphere^1 \hookrightarrow \sphere^3
\stackrel{\pi}{\longrightarrow} \sphere^2
$$
and its equivariance with respect to rotations
on $\sphere^3$ and $\sphere^2$ - 
see Lemma \ref{Hopf.properties} below - 
and also with respect to the first variation of the Willmore-energy along Hopf-tori in $\sphere^3$ and closed curves in $\sphere^2$, see formula \eqref{DHopf.first.Variation} below.
In order to work with the most effective formulation of the
Hopf-fibration, we consider $\sphere^3$ as the subset of the
four-dimensional $\rel$-vector space $\quat$ of quaternions, whose elements have length $1$, i.e.
$$
\sphere^3:=\{q \in \quat \,|\, \bar q\cdot q=1 \}.
$$
We shall use the usual notation for the generators of
the division algebra $\quat$, i.e. 1,$i$,$j$,$k$.
We therefore decompose every quaternion in the way
$$
q=q_1+i\, q_2 +j\, q_3 + k\,q_4,
$$
for unique ``coordinates'' $q_1, q_2, q_3, q_4 \in \rel$,
such that in particular there holds\\
$\bar q = q_1 - i\, q_2 - j\, q_3 - k\,q_4$.
Moreover, we identify
$$
\sphere^2 =\{q \in \textnormal{span}\{1,j,k\}\, |\, 
\bar q\cdot q=1 \}
=\sphere^3 \cap \textnormal{span}\{1,j,k\},
$$
and use the particular involution $q \mapsto \tilde q$ on $\quat$, which fixes the generators $1$, $j$ and $k$, but sends $i$ to $-i$. Following \cite{Pinkall}, we employ this involution to write the Hopf-fibration in the elegant way
\begin{equation} \label{Hopf}
\pi :\quat \longrightarrow  \quat, \quad q \mapsto \tilde q \cdot q,
\end{equation}
for $q \in \quat$. We shall recall here its most important, algebraic properties from Lemma 1 in \cite{Ruben.MIWF.II} 
in the following lemma, without proof.
\begin{lemma}  \label{Hopf.properties}
	\begin{itemize}
	\item[1)] $\pi(\sphere^3)=\sphere^2$, and moreover 
	$\pi(e^{i\varphi} q )= \pi(q)$,  $\forall \,
	\varphi \in \rel$ and $\forall \,q\in \sphere^3$.
	\item[2)] There holds
	$$
	\pi(q \cdot r) = \tilde r\cdot \pi(q) \cdot r,  \quad \forall  \,q,r\in \sphere^3.
	$$
	This formula expresses the fact that right multiplication on $\sphere^3$ translates equivariantly via $\pi$ to rotation in $\sphere^2$,  
	because every $r \in \sphere^3$ induces the rotation
	$q \mapsto \tilde r \cdot q \cdot r$ on $\sphere^2$.
	\item[3)] The derivative of $\pi$ in any $q\in \quat$,
	applied to some $v \in \quat$, reads:
	$$
	D\pi_q(v) = \,\tilde v \cdot q + \tilde q \cdot v.
	$$
\end{itemize}
\end{lemma}
\noindent \\
By means of the Hopf-fibration we introduce Hopf-tori  
as in Definition 1 of \cite{Ruben.MIWF.II}. 
The only slight difference between our 
Definition \ref{Hopf.torus.immersion} and 
Definition 1 in \cite{Ruben.MIWF.II} is that we define here  ``$C^k$-Hopf-tori'' - for any $k\in \nat$ -
as preimages with respect to the Hopf-fibration of closed 
$C^k$-curves in $\sphere^2$ without any reference to  particularly useful or canonical parametrizations of 
these subsets of $\sphere^3$, whereas in the $C^{\infty}$-smooth case there is no need for such a technical distinction on account of Lemma \ref{closed.lifts} below    
\footnote{See here also \cite{Pinkall} for further explanations
and applications concerning Hopf-tori.}.
\begin{definition} \label{Hopf.torus.immersion}
Let $\gamma: [a,b] \longrightarrow \sphere^2$ be a regular and closed path in $\sphere^2$ of regularity class $C^k$, 
with $k\in \nat$ or $k=\infty$.
	\begin{itemize}
	\item[1)] We call the preimage 
	$\pi^{-1}(\textnormal{trace}(\gamma))$  
	the \emph{Hopf-torus in $\sphere^3$ corresponding 
	to $\gamma$}, or less precisely a ``$C^k$-Hopf-torus'' 
    in $\sphere^3$. 	
	\item[2)] In the smooth case ``$k=\infty$'', 
	we can consider a smooth lift $\eta: [a,b] \longrightarrow \sphere^3$ of $\gamma$ with respect to $\pi$ into $\sphere^3$, 
	i.e. a map from $[a,b]$ into $\sphere^3$ of class $C^{\infty}$ satisfying $\pi \circ \eta = \gamma$.
	\footnote{See here Lemma \ref{closed.lifts} below regarding existence and uniqueness of such smooth lifts.}
	We define
	\begin{equation}   \label{Hopf.Torus}
	X(s,\varphi):= e^{i\varphi}  \cdot \eta(s), \quad
	\forall\, (s,\varphi) \in [a,b] \times [0,2\pi],
	\end{equation}
	and we note that $(\pi \circ X)(s,\varphi) = \gamma(s)$ holds $\forall \,(s,\varphi)\in [a,b] \times [0,2\pi]$.
	\item[3)] In the smooth case ``$k=\infty$'', 
	we call the map $X$ appearing in \eqref{Hopf.Torus} the 
	\emph{standard parametrization} of  
	the smooth Hopf-torus $\pi^{-1}(\textnormal{trace}(\gamma))$, or less precisely a smooth ``Hopf-torus-immersion''. 
	\end{itemize}
\end{definition}
\noindent 
In order to compute the position of the projection 
of the conformal structure of a given Hopf-torus into the moduli space in terms of its profile curve $\gamma$,
we introduce \emph{abstract Hopf-tori}:
\begin{definition} \label{Hopf.curve}
Let $\gamma: [0,L/2] \to \sphere^2$ be a path
with constant speed $2$ which traverses an embedded, 
closed, smooth curve in $\sphere^2$ of length
$L>0$ and encloses the area $A$ of the domain on $\sphere^2$, ``which lies on the left hand side'' when performing one loop through $\textnormal{trace}(\gamma)$.
We assign to $\gamma$ the lattice $\Gamma_{\gamma}$, which
is generated by the vectors $(2\pi,0)$ and $(A/2,L/2)$.
We call the torus $M_{\gamma}:=\com/\Gamma_{\gamma}$
the \emph{abstract Hopf-torus} corresponding to $\textnormal{trace}(\gamma)$. 
\end{definition}
\noindent 
In this context we should quote the following result, which is Proposition 1 in \cite{Pinkall}. 
\begin{proposition} \label{Pinkall} 
Let $\gamma:\sphere^1 \longrightarrow \sphere^2$ parametrize 
an embedded, closed, smooth curve of length $L$, which encloses the area $A$ in the sense of Definition 
\ref{Hopf.curve}. Its associated embedded Hopf-torus $\pi^{-1}(\textnormal{trace}(\gamma)) \subset 
\sphere^3$ endowed with the Euclidean metric of the ambient space $\rel^4$ is conformally equivalent to its corresponding \emph{abstract Hopf-torus} $\com/\Gamma_{\gamma}$ in the sense of Definition \ref{Hopf.curve}. In particular, the projection of the point $(A/4\pi,L/4\pi) \in \quat$ 
into the moduli space $\quat/\textnormal{PSL}_2(\ganz)$   
yields exactly the isomorphism class of the conformal class of 
$(\pi^{-1}(\textnormal{trace}(\gamma)), g_{\textnormal{euc}})$, when interpreted as a Riemann surface.  	 
\end{proposition}

\begin{remark} \label{Surfaces.genus.1}
On account of Proposition \ref{Pinkall} 
the conformal structures $[\pi^{-1}(\textnormal{trace}(\gamma))]$ and $[M_{\gamma}]$ induced by $\pi^{-1}(\textnormal{trace}(\gamma))\subset \sphere^3$ and $M_{\gamma}\subset \com$ correspond to each other via some suitable conformal diffeomorphism between $M_{\gamma}$ and $\pi^{-1}(\textnormal{trace}(\gamma))$, and their common projection into moduli space $\Mill_1$ lies in some prescribed compact subset 
$K$ of $\Mill_1$, if and only if the pair of ``moduli'' $(A/2,L/2)$ of $M_{\gamma}$
is situated at some sufficiently large distance to the boundary of $\quat$, i.e. if and only if the length $L$ of $\gamma$ is bounded from above and from below by two appropriate positive numbers $R_1(K),R_2(K)$. We will use this fact below in the proof of Proposition \ref{L2.estimates.1}, Part (5).  
\end{remark}
In order to rule out unnecessarily complicated 
parametrizations of Hopf-tori, and especially in order to avoid technical mistakes regarding formulae \eqref{Will.Wil} and \eqref{First.Var.Will.Wil} below, we follow 
here the strategy in \cite{Ruben.MIWF.II} and 
introduce ``simple'' parametrizations of smoothly 
immersed tori - a modified, more effective version of Definition 3 in \cite{Ruben.MIWF.II}:
\begin{definition} \label{Simple.map}
Let $\Sigma$ be a compact smooth torus and 
$F:\Sigma \longrightarrow \sphere^3$ a smooth immersion. 
We call $F$ a \emph{simple parametrization} of the 
immersed torus $F(\Sigma) \subset \sphere^3$, if there holds 
$\sharp \{F^{-1}(z)\} =1$ in $\Hn$-a.e. $z \in F(\Sigma)$. 
\end{definition}

\begin{remark} \label{Simple.maps}
\begin{itemize} 
\item[1)] We should recall here that a regular closed path 
$\gamma:\sphere^1 \longrightarrow \rel^n$ cannot map 
$\sphere^1$ to a point and closes up after a finite number of loops through the trace of $\gamma$. This implies 
that the number of preimages $\sharp \{\gamma^{-1}(z)\}$
is a unique natural number $l=l(\gamma)$ for all 
but finitely many $z\in \textnormal{trace}(\gamma)$ - 
the \emph{self-intersections} of $\textnormal{trace}(\gamma)$. 
The number $l=l(\gamma)$ is exactly the above number of loops 
which the path $\gamma$ travels through $\textnormal{trace}(\gamma)$. Since every 
smooth compact torus $\Sigma$ is homeomorphic to the 
product $\sphere^1\times \sphere^1$, one can argue similarly that a smooth immersion $F:\Sigma \longrightarrow \sphere^3$ cannot map $\Sigma$ into any set of Hausdorff-dimension 
strictly less than $2$ and has to wrap the torus $\Sigma$ finitely often about the doubly periodic, immersed 
image $F(\Sigma)$. 
More precisely this means: for every smooth immersion 
$F:\Sigma \longrightarrow \sphere^3$ there is some natural number $k=k(F)$, such that $\sharp \{F^{-1}(z)\}=k$ in $\Hn$-a.e. $z\in F(\Sigma)$.
\item[2)] For any smooth immersion 
$F:\Sigma \longrightarrow \sphere^3$ we can estimate 
the number $k=k(F)$ of wraps about $F(\Sigma)$ 
by means of the proven Willmore conjecture,  
Theorem A in \cite{Marques.Neves.2014}, provided the 
Willmore energy of $F$ can be estimated from above. 
More precisely, from the inequality $\Will(F)< 2K \pi^2$ 
it follows that $1\leq k=k(F)< K$. For, if we had here $k=k(F)\geq K$, then: 
$$ 
\Will(F) = k \cdot \Will(F(\Sigma)) \geq k \cdot 2\pi^2 
\geq 2K \,\pi^2,  
$$ 
which contradicts our assumption on $\Will(F)$. 
In particular, for any smooth immersion 
$F:\Sigma \longrightarrow \sphere^3$ 
with $\Will(F)< 4 \pi^2$ we can determine that 
$k(F)=1$, i.e. that $F$ is a \emph{simple parametrization} 
of the immersed torus $F(\Sigma)$ in the sense of 
Definition \ref{Simple.map}.  
\item[3)] 
In the special case in which $F(\Sigma)$ is a compact manifold of genus $1$, i.e. a smooth torus in the sense of differential topology - the condition in Definition \ref{Simple.map} is satisfied, if and only if the induced homomorphism 
$$
(F_{*})_2: H_2(\Sigma,\ganz) \stackrel{\cong}\longrightarrow
H_{2}(F(\Sigma),\ganz)
$$
is an isomorphism between these two singular homology 
groups in degree $2$; see here Remark 4 in \cite{Ruben.MIWF.II}.
\end{itemize}
\end{remark}
\noindent 
In the next two propositions we recall some basic differential geometric formulae from Propositions 3 and 4 in \cite{Ruben.MIWF.II}, which particularly yield the useful correspondence between the MIWF and the degenerate variant \eqref{elastic.energy.flow} of the classical 
elastic energy flow \eqref{classical.elastic.energy.flow}; 
see Proposition \ref{correspond.flows} below. The proofs are either straight forward or can be found in \cite{Ruben.MIWF.II}.   
\begin{proposition} \label{compute.the.operator}
	The $L^2$-gradient of the elastic energy
	$\Wil(\gamma):= \int_{\sphere^1} 1 + |\vec{\kappa}_{\gamma}|^2 \, d\mu_{\gamma}$, with 
	$\vec{\kappa}_{\gamma}$ as in \eqref{Curva.vector} below, evaluated in an arbitrary closed curve $\gamma \in
	C^{\infty}_{\textnormal{reg}}(\sphere^1,\sphere^2)$, reads exactly:
\begin{equation}  \label{first.variation}
\nabla_{L^2} \Wil(\gamma)(x) = 2\,
\Big{(}\nabla_{\frac{\gamma'}{|\gamma'|}}^{\perp}\Big{)}^2
(\vec \kappa_{\gamma})(x)
+\, |\vec{\kappa}_{\gamma}|^2(x) \,\vec{\kappa}_{\gamma}(x)
+\, \vec{\kappa}_{\gamma}(x),  \qquad \textnormal{for} \,\, 
x\in \sphere^1,
\end{equation}
where we denote in \eqref{first.variation} and in the sequel by $\nabla_{\frac{\gamma'}{|\gamma'|}}(W)$ 
the classical covariant derivative of some 
smooth tangent vector field $W$ on $\sphere^2$ 
along the given curve $\gamma$ with respect to the unit tangent vector field 
$\frac{\partial_x\gamma}{|\partial_x\gamma|}$
along $\gamma$ and moreover by 
$\nabla^{\perp}_{\frac{\gamma'}{|\gamma'|}}(W)$ 
the orthogonal projection of the tangent vector
field $\nabla_{\frac{\gamma'}{|\gamma'|}}(W)$ into the 
normal bundle of the given curve $\gamma$ within $T\sphere^2$. 
Abbreviating furthermore $\partial_s\gamma:=\frac{\partial_x\gamma}{|\partial_x\gamma|}$\, we can reformulate the leading term of the right hand side of equation \eqref{first.variation} as: 
\begin{eqnarray}  \label{split.operator}
\Big{(}\nabla_{\frac{\gamma'}{|\gamma'|}}^{\perp}\Big{)}^2
(\vec{\kappa}_{\gamma})(x)
=\Big{(}\nabla_{\partial_s\gamma}^{\perp}\Big{)}^2
\big{(}(\partial_{ss}\gamma)(x) - \langle \gamma(x), \partial_{ss}\gamma(x) \rangle \,\gamma(x) \big{)}     \nonumber \\
= (\partial_s)^{4}(\gamma)(x)
- \langle (\partial_s)^{4}(\gamma)(x),\gamma(x) \rangle \,\gamma(x)         \\
- \,\langle (\partial_s)^{4}(\gamma)(x),\partial_s\gamma(x) \rangle \,\partial_s\gamma(x)
+ \,|(\nabla_{\partial_s\gamma})^2(\gamma)(x)|^2 \,\partial_{ss}\gamma(x)
\quad \textnormal{for} \,\, x\in \sphere^1.   \nonumber
\end{eqnarray}
The fourth normalized derivative $(\partial_s)^{4}(\gamma)\equiv
\Big{(} \frac{\partial_x}{|\partial_x \gamma|}\Big{)}^4(\gamma)$ is non-linear with respect to $\gamma$, and at least its leading term can be computed in terms of ordinary partial derivatives of $\gamma$:
	\begin{eqnarray}  \label{split.operator.x}
	(\partial_s)^4(\gamma)
	=\frac{(\partial_x)^4(\gamma)}{|\partial_x\gamma|^4}
	- \frac{1}{|\partial_x\gamma|^4}\,
	\Big{\langle} (\partial_x)^4(\gamma), \frac{\partial_x\gamma}{|\partial_x\gamma|} \Big{\rangle}\, \frac{\partial_x\gamma}{|\partial_x\gamma|}  \nonumber \\
	 + \,C((\partial_x)^2(\gamma),\partial_x(\gamma))  \cdot  (\partial_x)^3(\gamma)    \\
	+ \textnormal{rational expressions which only involve} \,\,(\partial_x)^2(\gamma)
	\,\, \textnormal{and} \,\, \partial_x(\gamma),
	\nonumber
	\end{eqnarray}
	where $C:\rel^6 \longrightarrow \textnormal{Mat}_{3,3}(\rel)$ is a $\textnormal{Mat}_{3,3}(\rel)$-valued function whose components are rational functions in $(y_1,\ldots,y_6) \in \rel^6$.      	
\end{proposition}

\begin{proposition}  \label{Hopf.Willmore.prop}
Let $F:\Sigma \longrightarrow \sphere^3$ be a smooth immersion which maps the compact torus $\Sigma$ simply 
onto some smooth Hopf-torus in $\sphere^3$ in the sense of Definition \ref{Simple.map}, and let $\gamma:\sphere^1 \longrightarrow \sphere^2$ be a smooth regular parametrization of the closed curve $\textnormal{trace}(\pi \circ F)$ 
which performs exactly one loop through its trace. 
Let's moreover use quaternionic notation - 
as explained at the beginning of this section - 
in order to introduce the curvature vector
\begin{equation}  \label{Curva.vector}
\vec{\kappa}_{\gamma} \equiv \vec{\kappa}^{\sphere^2}_{\gamma}
:= - \frac{1}{|\gamma'|^2} \,
(\overline{\gamma'} \cdot \nu'_{\gamma}) \cdot \nu_{\gamma}
\end{equation}
and the signed curvature $\kappa_{\gamma}:=\langle \vec{\kappa}_{\gamma},\nu_{\gamma} \rangle_{\rel^3}$  
along the given path $\gamma$, where 
$\nu_{\gamma} \in \Gamma(\gamma^*T\sphere^2)$ denotes some fixed unit normal field along $\gamma$.
Then there is some $\varepsilon =\varepsilon(F,\gamma)>0$
such that for an arbitrarily fixed point $x^* \in \sphere^1$ the following differential-geometric formulae hold for the immersion $F$:
	\begin{equation}  \label{second.fundam.}
	A_{F,\sphere^3}(\eta_F(x))  = N_F(\eta_F(x))
	\begin{pmatrix}
		2 \kappa_{\gamma}(x)  &  1  \\
		1  &  0
	\end{pmatrix}
	\end{equation}
	where $\eta_F:\sphere^1\cap B_{\varepsilon}(x^*) \longrightarrow \Sigma$ denotes an arbitrary horizontal smooth lift of
	$\gamma\lfloor_{\sphere^1\cap B_{\varepsilon}(x^*)}$ with respect to
	the fibration $\pi \circ F$, as introduced in Lemma \ref{closed.lifts} below, and $N_F$ denotes a fixed unit normal field along the immersion $F$. This implies
	\begin{equation} \label{H}
	\vec H_{F,\sphere^3}(\eta_F(x)) = \textnormal{trace} A_{F,\sphere^3}(\eta_F(x))= 
	2 \kappa_{\gamma}(x)\, N_F(\eta_F(x))
	\end{equation}
	$\forall \, x  \in \sphere^1\cap B_{\varepsilon}(x^*)$,
	for the mean curvature vector of $F$ and also
	\begin{equation}  \label{second.fundam.trace.free}
	A^0_{F,\sphere^3}(\eta_F(x)) = N_F(\eta_F(x))
	\begin{pmatrix}
		\kappa_{\gamma}(x)  &  1  \\
		1  &   -\kappa_{\gamma}(x)
	\end{pmatrix}
	\end{equation}
	and $|A^0_{F,\sphere^3}|^2(\eta_F(x)) 
	= 2 (\kappa_{\gamma}(x)^2+1)$,
	\begin{eqnarray}  \label{Q.trace.free.H}
	Q(A^0_{F,\sphere^3})(\vec H_{F,\sphere^3})(\eta_F(x))  
	= 4 \,(\kappa_{\gamma}^3(x) + \kappa_{\gamma}(x)) \, N_F(\eta_F(x))
	\end{eqnarray}
	and finally
	\begin{equation}  \label{normal.laplace.H}
	\triangle^{\perp}_F (\vec H_{F,\sphere^3})(\eta_F(x))
	=  8\,  \Big{(} \nabla_{\frac{\gamma'}{|\gamma'|}}\Big{)}^2
	(\kappa_{\gamma})(x) \, N_F(\eta_F(x))
	\end{equation}
	and for the traced sum of all derivatives of $A_F$ 
	of order $k \in \nat$:
	\begin{equation}  \label{all.derivatives}
	|(\nabla^{\perp_F})^k(A_{F,\sphere^3})(\eta_F(x))|^2
	=  2^{2+2k} \, \Big{|} \Big{(} \nabla^{\perp}_{\frac{\gamma'}{|\gamma'|}}\Big{)}^k
	(\vec \kappa_{\gamma})(x)\Big{|}^2
	\end{equation}
	$\forall \,x \in \sphere^1\cap B_{\varepsilon}(x^*)$.
	In particular, we derive
	\begin{eqnarray} \label{MIWF.Hopf.Tori}
	\nabla_{L^2} \Will(F)(\eta_F(x)) 
	= 2 \,\Big{(} 2\,\Big{(} \nabla_{\frac{\gamma'}{|\gamma'|}}\Big{)}^2
	(\kappa_{\gamma})(x) + \kappa_{\gamma}^3(x) +
	\kappa_{\gamma}(x) \Big{)} \, N_F(\eta_F(x)),
	\end{eqnarray}
	and the \emph{Hopf-Willmore-identity}:
	\begin{eqnarray} \label{DHopf.first.Variation}
	D\pi_{F(\eta_F(x))}.\Big{(} 
	\nabla_{L^2} \Will(F)(\eta_F(x)) \Big{)}   
	= 4 \, \nabla_{L^2} \Wil(\gamma)(x)   
	\end{eqnarray}
	$\forall \,x \in \sphere^1\cap B_{\varepsilon}(x^*)$,
	where there holds $\pi \circ F \circ \eta_F = \gamma$ on
	$\sphere^1\cap B_{\varepsilon}(x^*)$, as in Lemma \ref{closed.lifts} below. Finally, we have
	\begin{equation}  \label{Will.Wil}
	\Will(F) \equiv \int_{\Sigma} 1 + \frac{1}{4} \,
	|\vec H_{F,\sphere^3}|^2 \, d\mu_{F}
	=\pi  \,\int_{\sphere^1} 1 + \,|\kappa_{\gamma}|^2 \, d\mu_{\gamma} = \pi  \,\Wil(\gamma), 
	\end{equation}
	and 
	\begin{eqnarray}  \label{First.Var.Will.Wil}
	\int_{\Sigma} \frac{1}{|A^0_{F,\sphere^3}|^4} \, |\nabla_{L^2}\Will(F)|^2 \, d\mu_{F} 
	= \pi  \,\int_{\sphere^1} 
	\frac{1}{(\kappa_{\gamma}^2+1)^2} 
	\Big{|}2\,\Big{(}\nabla^{\perp}_{\frac{\gamma'}{|\gamma'|}}\Big{)}^2(\vec{\kappa}_{\gamma})
	+ |\vec{\kappa}_{\gamma}|^2 \vec{\kappa}_{\gamma}
	+ \vec{\kappa}_{\gamma} \Big{|}^2 \, d\mu_{\gamma} \quad \\
	\equiv \pi  \,\int_{\sphere^1} \frac{1}{(\kappa_{\gamma}^2+1)^2} \,  |\nabla_{L^2}\Wil(\gamma)|^2 d\mu_{\gamma}.\quad  \nonumber
	\end{eqnarray} 
\end{proposition}
\noindent \\ 
Now we arrive at the main result of this section, 
a reduction of the MIWF to some \emph{degenerate variant} 
of the classical elastic energy flow on $\sphere^2$ by means
of the Hopf-fibration: 
\begin{proposition}  \label{correspond.flows}
	Let $[0,T] \subset \rel$ be a non-void compact interval, and let $\gamma_t: \sphere^1 \to \sphere^2$ be a smooth family of closed, smooth and regular paths, for $t \in [0,T]$. 
	Moreover, let $F_t:\Sigma \longrightarrow \sphere^3$
	be an arbitrary smooth family of smooth immersions, 
	parametrizing the Hopf-tori $\pi^{-1}(\textnormal{trace}(\gamma_t))\subset \sphere^3$, 
    for every $t\in [0,T]$. Then the following statement holds:\\
	The family of immersions $\{F_t\}$ moves according to the
	MIWF-equation \eqref{Moebius.flow} on 
	$[0,T] \times \Sigma$ - up to smooth, time-dependent reparametrizations $\Phi_t$ with $\Phi_0=\textnormal{id}_{\Sigma}$ - if and only if there 
	is a smooth family $\sigma_t:\sphere^1 \to \sphere^1$ of reparametrizations with $\sigma_0=\textnormal{id}_{\sphere^1}$, such that
	the family $\{\gamma_t\circ \sigma_t\}$ satisfies the
	``degenerate elastic energy evolution equation''
	\begin{eqnarray}  \label{elastic.energy.flow}
	\partial_t \tilde \gamma_t =
	- \,\frac{1}{(\kappa_{\tilde \gamma_t}^2+1)^2}\,
	\Big{(} 2 \, \Big{(}\nabla^{\perp}_{\frac{\tilde \gamma_t'}{|\tilde \gamma_t'|}} \Big{)}^2(\vec{\kappa}_{\tilde \gamma_t})
	+ |\vec{\kappa}_{\tilde \gamma_t}|^2 \vec{\kappa}_{\tilde \gamma_t} + \vec{\kappa}_{\tilde \gamma_t} \Big{)}        \equiv - \,\frac{1}{(\kappa_{\tilde \gamma_t}^2+1)^2}\,
	\nabla_{L^2}\Wil(\tilde \gamma_t)
	\end{eqnarray}
	on $[0,T] \times \sphere^1$, where $\nabla_{L^2} \Wil$ denotes the $L^2$-gradient of the elastic energy $\Wil$, as above in Proposition \ref{compute.the.operator}.
\end{proposition}
\emph{Proof:\,} The proof is essentially an adaption of the proof of the corresponding Proposition 5 in \cite{Ruben.MIWF.II} - up to only minor modifications - and uses only the 
formulae of Proposition \ref{Hopf.Willmore.prop}.  \\\\
\noindent 
Moreover, we will need the following short-time existence and uniqueness result. 
\begin{proposition}  \label{short.time.existence.elastic} 
Let $\gamma_0:\sphere^1 \longrightarrow \sphere^2$ be a 
$C^{\infty}$-smooth, closed and regular curve. Then 
there is some small $T>0$ and a $C^{\infty}$-smooth solution 
$\{\gamma_t\}$ of the degenerate elastic energy flow \eqref{elastic.energy.flow} on $\sphere^1\times [0,T]$, 
starting in $\gamma_0$ at $t=0$.   	
\end{proposition}
\emph{Proof:\,} The proof works exactly as the proof of Theorem 3.1 in \cite{Dall.Acqua.Pozzi.2018}, where short-time existence and uniqueness of the classical elastic energy flow 
\begin{eqnarray}  \label{classical.elastic.energy.flow}
	\partial_t \gamma_t =
	- \,\Big{(} 2 \, \Big{(}\nabla^{\perp}_{\frac{\gamma_t'}{|\gamma_t'|}} \Big{)}^2(\vec{\kappa}_{\gamma_t})
	+ |\vec{\kappa}_{\gamma_t}|^2 \vec{\kappa}_{\gamma_t}
	+ \vec{\kappa}_{\gamma_t} \Big{)}                     
	\equiv - \,\nabla_{L^2}\Wil(\gamma_t)
\end{eqnarray}
for smooth curves $\gamma_t:\sphere^1 \longrightarrow \sphere^2$ is proved by means of a concrete, stereographic chart from $\rel^2$ onto $\sphere^2\setminus \{(0,0,1)\}$ 
and by means of normal representation 
$\hat \gamma_t:= \hat \gamma_0 +u_t \, N_{\hat \gamma_0}$ of 
the projected plane curves $\hat \gamma_t$ 
with respect to the projected plane initial curve  
$\hat \gamma_0:\sphere^1 \longrightarrow \rel^2$. 
The only difference here is the additional factor  
$\frac{1}{(\kappa_{\gamma_t}^2+1)^2}$ arising in 
\eqref{elastic.energy.flow} in front of the right hand side 
of \eqref{classical.elastic.energy.flow}. Since the 
curvature vector $\vec \kappa_{\hat \gamma_t}$ of the projected curve $\hat \gamma_t = 
\hat \gamma_0 +u_t \, N_{\hat \gamma_0}$ 
can be explicitly computed here just as in the proof of Theorem 3.1 in \cite{Dall.Acqua.Pozzi.2018} 
and since the factor 
$\frac{1}{((\kappa_{\hat \gamma_0 + u_t \, 
N_{\hat \gamma_0}})^2+1)^2}$
- now to be multiplied with the right hand side of 
equation (3.1) in \cite{Dall.Acqua.Pozzi.2018} - 
is bounded from above by $1$ and from below by 
$\frac{1}{2} \frac{1}{(\max_{\sphere^1}(\kappa_{\hat \gamma_0})^2+1)^2}$ on $\sphere^1\times [0,T]$ for every perturbation $\{u_t\}$ which is sufficiently small in 
$C^{2+\alpha,\frac{2+\alpha}{4}}(\sphere^1 \times [0,T],\rel^2)$, hence the decisive arguments of the proof of Theorem 3.1 in \cite{Dall.Acqua.Pozzi.2018} - employing linearization of the quasilinear parabolic differential equation (3.1) in \cite{Dall.Acqua.Pozzi.2018} and parabolic 
Schauder theory - can be adopted here without any changes. \\\\
\noindent 
In the following proposition we collect the most 
fundamental information about flow lines of evolution 
equation \eqref{elastic.energy.flow} - our degenerate 
variant of the elastic energy flow.
\begin{proposition} \label{L2.estimates.1}
Let $\{\gamma_t\}_{t\in [0,T)}$, with either $T>0$ or $T=\infty$, be a flow line of evolution equation \eqref{elastic.energy.flow}, starting in a 
smooth closed path $\gamma_0:\sphere^1\longrightarrow \sphere^2$ with elastic energy $\Wil_0:=\Wil(\gamma_0)\leq 8$. Then the following statements hold: 
\begin{itemize} 
\item[1)] The elastic energy $\Wil(\gamma_t)$ along the flow line $\{\gamma_t\}$ is monotonically decreasing and stays strictly smaller than $8$ for $t \in (0,T)$,
\item[2)] the curves $\gamma_t$ are smooth embeddings of $\sphere^1$ into $\sphere^2$ for $t \in (0,T)$, 
\item[3)] the lengths of the curves $\gamma_t$ 
are uniformly bounded from above by $\Wil_0$ and 
from below by $\pi$ for $t \in [0,T)$, 
\item[4)] the areas of the domains $\Omega_t$ lying on the left hand sides of the embedded curves $\gamma_t$ in $\sphere^2$ - see Definition \ref{Hopf.curve} - remain 
strictly bigger than $2\,(\pi - 2)$ and strictly smaller 
than $2\,(\pi + 2)$ for $t \in (0,T)$,  
\item[5)] the conformal structures 
induced by the Euclidean metric of $\sphere^3$ restricted to the embedded Hopf-tori $\pi^{-1}(\textnormal{trace}(\gamma_t))$ lie in a compact subset of the moduli space 
$\Mill_1\cong \quat/\textnormal{PSL}_2(\ganz)$, 
for every $t\in (0,T)$. 
\end{itemize} 
\end{proposition}
\emph{Proof:\,}   
First of all, just as in \eqref{Willmore.monotonicity} 
we argue that any flow line $\{\gamma_t\}$ of equation \eqref{elastic.energy.flow} satisfies:   
\begin{eqnarray} \label{Pseudo.Gradient.flow}
	\frac{d}{dt} \Wil(\gamma_t)
	= \langle \nabla_{L^2} \Wil(\gamma_t), \partial_t \gamma_t
	\rangle_{L^2(\mu_{\gamma_{t}})} = \\
	= - \int_{\sphere^1}  \frac{1}{(\kappa_{\gamma_t}^2+1)^2}\,
	\Big{|} 2 \,\Big{(}\nabla^{\perp}_{\gamma_t'/|\gamma_t'|}\Big{)}^2
	(\vec{\kappa}_{\gamma_t})
	+ | \vec{\kappa}_{\gamma_t}|^2\, \vec\kappa_{\gamma_t} +
	\,\vec{\kappa}_{\gamma_t} \Big{|}^2 \, d\mu_{\gamma_t}  \leq 0    \nonumber
\end{eqnarray}
for every $t\in [0,T)$. Now we recall that 
$\Wil(\gamma_0) \leq 8$. As in the proof of Theorem 
\ref{limit.MIWF} we shall distinguish  
the two cases in which there either 
holds (a) $\Wil(\gamma_t)=8$ on some 
arbitrarily short, but non-empty time interval $[0,\varepsilon)$, or (b) $\Wil(\gamma_t)<8$ for every 
$t\in (0,T)$. 
On account of inequality \eqref{Pseudo.Gradient.flow} 
we can argue exactly as in the proof of Theorem \ref{limit.MIWF} and infer that case (a) actually means that 
$\{\gamma_t\}$ can be smoothly continued as  
a global flow line of flow \eqref{elastic.energy.flow} 
only consisting of the initial curve $\gamma_0$,  
which would have to be here a closed elastic curve in 
$\sphere^2$ with elastic energy $8$. But on account of Proposition 6 in \cite{Ruben.MIWF.II} there are no closed elastic curves in $\sphere^2$ whose elastic energy 
lies in the interval $(2\pi,4 \pi)$
\footnote{In Proposition 6 of \cite{Ruben.MIWF.II} the 
elastic energy was restricted to closed paths in $\sphere^2$ which traverse their traces only once, leading to  
the large energy gap 
$\big{(}2\pi,\frac{8 \pi}{\sqrt{2}}\big{]}$. 
But in the situations of Proposition \ref{L2.estimates.1}, 
Theorem \ref{limit.at.infinity} or Corollary \ref{limit.at.infinity.1} we also have to take elastic curves of higher multiplicities into account, i.e. 
stationary closed paths of $\Wil$ which might perform 
several loops. Indeed, in this broader sense the 
``double loop great circle'' is an elastic curve of 
energy $4 \pi$, whereas the 
``triple loop great circle'' already has energy 
$6\pi>\frac{8 \pi}{\sqrt{2}}$.} - 
obviously containing the value $8$.     
Therefore case (a) is automatically excluded for 
the degenerate elastic energy flow  \eqref{elastic.energy.flow}, and we must have 
$\Wil(\gamma_t)<8$ for every $t\in (0,T)$, whenever 
$\Wil(\gamma_0)\leq 8$ is required.     
This proves already the first assertion of the proposition.
\\ 
Combining this with formula \eqref{Will.Wil} and with the Li-Yau inequality, we infer that the corresponding Hopf-tori $\pi^{-1}(\textnormal{trace}(\gamma_t))$ are embedded surfaces in $\sphere^3$ for positive times $t$, implying that their profile curves $\gamma_t$ have to map $\sphere^1$ injectively into $\sphere^2$, i.e. that $\gamma_t$ have to be
smooth embeddings, at least for positive times. \\
As for the third assertion, we infer from inequality \eqref{Pseudo.Gradient.flow} 
in particular the two inequalities:
\begin{eqnarray}  \label{Length.total}
\textnormal{length}(\gamma_t) \leq \Wil(\gamma_0),      \\
\int_{\sphere^1} |\vec{\kappa}_{\gamma_t}^{\sphere^2}|^2 \, 
d\mu_{\gamma_t} \leq \Wil(\gamma_0),    \label{L2.bounded}
\end{eqnarray}
for every $t \in [0,T)$. Applying the elementary inequality
$$
\Big{(} \int_{\sphere^1} |\vec{\kappa}_{\gamma}^{\sphere^2}| \,
d\mu_{\gamma}\Big{)}^2 \geq  4\pi^2 - \textnormal{length}(\gamma)^2,
$$
which holds for closed smooth regular paths
$\gamma:\sphere^1 \longrightarrow \sphere^2$,  
see here \cite{Teufel}, one can easily derive the 
lower bound
\begin{equation} \label{lower.bound.length}
\textnormal{length}(\gamma) \geq
\min\Big{\{} \pi, \frac{3\pi^2}{\Wil(\gamma)} \Big{\}},
\end{equation}	
see Lemma 2.9 in \cite{Dall.Acqua.Pozzi.2018}, for any
closed smooth regular path 
$\gamma:\sphere^1 \longrightarrow \sphere^2$.  
In combination with the monotonicity of $\Wil(\gamma_t)$,  i.e. with statement \eqref{Pseudo.Gradient.flow}, 
and with the requirement $\Wil(\gamma_0)\leq 8$ we obtain:  
\begin{equation}  \label{bound.below}
\textnormal{length}(\gamma_t) \geq
\min\Big{\{} \pi, \frac{3\pi^2}{\Wil(\gamma_0)} \Big{\}}=\pi,
\end{equation} 
for every $t \in [0,T)$, which proves the third assertion. \\
The fourth assertion of the proposition now follows from 
Gauss-Bonnet's Theorem for simply connected subdomains 
$\Omega$ of $\sphere^2$ with smooth boundary 
$\partial \Omega$: 
$$ 
\Hn(\Omega)+ \int_{\partial \Omega} 
\kappa_{\partial \Omega}^{\sphere^2} \,d\Hnm =2\pi, 
$$ 
yielding here together with Cauchy-Schwarz' inequality,
with the first statement of the proposition, especially 
with $\Wil(\gamma_t) < 8$ for positive times $t$: 
$$ 
2\,(\pi +2) > 2\pi + \int_{\sphere^1} |\kappa_{\gamma_t}^{\sphere^2}| \, d\mu_{\gamma_t} \geq \Hn(\Omega_t) \geq 2\pi - \int_{\sphere^1} 
|\kappa_{\gamma_t}^{\sphere^2}| \, d\mu_{\gamma_t} 
> 2\,(\pi-2),
$$
for every for every $t \in (0,T)$. 
The fifth assertion of the proposition follows immediately from the second and third statement of this proposition combined with Remark \ref{Surfaces.genus.1}.

\section{Proofs of Theorems 
\ref{singular.time.MIWF.Hopf.tori} and 
\ref{limit.at.infinity}} 
\noindent 
\underline{{\bf Proof of Theorem \ref{singular.time.MIWF.Hopf.tori}}}:
\begin{itemize} 
\item[1)] 
First of all, on account of Remark \ref{Simple.maps} (2) 
the fact that the considered flow line $\{F_t\}$ satisfies 
$\Will(F_t)\leq \Will(F_0) \leq 8 \pi$, 
for every $t\in [0,T_{\textnormal{max}}(F_0))$, 
implies that each immersion $F_t$ belonging to the 
flow line $\{F_t\}$ is a \emph{simple} parametrization of its image $F_t(\Sigma)$, in the sense of Definition \ref{Simple.map}. Moreover, by means of Proposition \ref{correspond.flows} and on account of the uniqueness of classical flow lines of both the MIWF in $\sphere^3$ and the degenerate elastic energy flow 
\eqref{elastic.energy.flow} in $\sphere^2$ 
one can easily show - as in the proof of the 
first part of Theorem 1 in \cite{Ruben.MIWF.II} treating 
the classical Willmore flow in $\sphere^3$ - that a flow line $\{F_t\}$ of the MIWF has to consist of smooth parametrizations of Hopf-tori in $\sphere^3$, whenever it starts moving in a smooth parametrization $F_0$ of a Hopf-torus in $\sphere^3$, and that such a flow line $\{F_t\}$ is projected by the Hopf-fibration onto a smooth flow line $\{\gamma_t\}_{t\in [0,T_{\textnormal{max}}(F_0))}$ 
of the degenerate elastic energy flow \eqref{elastic.energy.flow}, which satisfies $\textnormal{trace}(\gamma_t)=\pi(F_t(\Sigma))$ 
for every $t\in [0,T_{\textnormal{max}}(F_0))$
and exactly $T_{\textnormal{max}}(\gamma_0) 
=T_{\textnormal{max}}(F_0)$. 
Combining now this particular 
\emph{correspondence of flows lines} with identity \eqref{Will.Wil} and with the first part of Proposition \ref{L2.estimates.1} we conclude that already 
$\Will(F_t)<8 \pi$ must hold for every 
$t\in (0,T_{\textnormal{max}}(F_0))$,
implying that the simple Hopf-torus-immersions $F_t$ 
have to be embeddings, at least for positive $t$.   
Now, in order to prove the first statement of Theorem 
\ref{singular.time.MIWF.Hopf.tori} we consider the 
sequence $t_j \nearrow T_{\textnormal{max}}(F_0)$ 
appearing in the statement of the theorem,   
and we recall \eqref{weak.convergence.mu}, i.e. that 
there is some suitable subsequence $\{F_{t_{j_l}}\}$ of  
$\{F_{t_{j}}\}$ and some integral, $2$-rectifiable varifold $\mu$ in $\rel^4$ such that
\begin{equation}  \label{weak.convergence.F.j}  
\Hn\lfloor_{F_{t_{j_l}}(\Sigma)} \longrightarrow \mu 
\qquad \textnormal{weakly as Radon measures on} \,\,\,\rel^4, 
\end{equation}
as $l \to \infty$.  
As mentioned above, we know that for every fixed 
time $t\in [0,T_{\textnormal{max}}(F_0))$ the surface $F_{t}(\Sigma)$ is a Hopf-torus in $\sphere^3$ and that we can parametrize its projection $\pi(F_t(\Sigma))$ into $\sphere^2$ by means of a smooth closed path 
$\gamma_t:\sphere^1 \longrightarrow \sphere^2$, such that 
the resulting family of paths $\{\gamma_t\}_{t\in [0,T_{\textnormal{max}}(F_0))}$ constitutes a maximal flow line of flow \eqref{elastic.energy.flow}, 
starting with elastic energy $\Wil_0=\Wil(\gamma_0)\leq 8$. 
From Proposition \ref{L2.estimates.1} we infer 
that the curves $\gamma_t:\sphere^1\longrightarrow \sphere^2$ 
are smooth embeddings whose lengths are uniformly bounded 
from below and which additionally enclose simply connected domains $\Omega_{t} \subset \sphere^2$ 
- in the sense of Definition \ref{Hopf.curve} - 
whose $\Hn$-measures are bounded from below by the positive 
number $2(\pi-2)$, for all 
$t\in [0,T_{\textnormal{max}}(F_0))$. This implies first 
of all, that the diameters of the Hopf-tori 
$F_{t_{j_l}}(\Sigma)=\pi^{-1}(\textnormal{trace}
(\gamma_{t_{j_l}}))$ - interpreted as subsets of $\rel^4$ - are bounded from below for every $l\in \nat$. Hence, it follows as in the proof of Proposition 2.2 in \cite{Schaetzle.Conf.factor.2013} that the integral 
limit varifold $\mu$ in \eqref{weak.convergence.F.j}  
satisfies $\mu \not =0$. Therefore,     
we can infer from the first part of Theorem \ref{limit.MIWF} 
that $\textnormal{spt}(\mu)$ is an embedded, closed and orientable Lipschitz-surface in $\sphere^3$ either 
of genus $0$ or of genus $1$. Moreover, we know by weak 
lower semicontinuity of the Willmore energy $\Will$ with respect to the convergence in \eqref{weak.convergence.F.j} 
- see \cite{Schaetzle.Lower.semicont.2009} - that 
\begin{equation}  \label{Willmore.mu}
\Will(\mu):=\frac{1}{4} \,\int_{\rel^4} |\vec H_{\mu}|^2 \, d\mu \leq \liminf_{l\to \infty}  
\frac{1}{4} \,\int_{\rel^4} |H_{F_{t_{j_l}},\rel^4}|^2 \, 
d\mu_{F_{t_{j_l}}^*g_{\textnormal{euc}}}  
< \Will(F_0) \leq 8\pi,
\end{equation}  
and as in the proof of the first part of Theorem \ref{limit.MIWF} we can also infer from the fact that 
$F_{t_{j_l}}(\Sigma) \subset \sphere^3$\, 
$\forall \, l\in \nat$, from $\mu \not =0$ and from convergence \eqref{weak.convergence.F.j} the convergence 
of the embedded surfaces $F_{t_{j_l}}(\Sigma)$ to $\textnormal{spt}(\mu)$ in Hausdorff distance as 
$l\to \infty$, i.e. convergence \eqref{Hausdorff.converg}.  
Because of $F_{t_{j_l}}(\Sigma)=
\pi^{-1}(\textnormal{trace}(\gamma_{t_{j_l}}))$ we 
can apply the Hopf-fibration to the latter convergence 
and infer:
\begin{equation}  \label{Hausdorff.convergence.2} 
\textnormal{trace}(\gamma_{t_{j_l}}) \longrightarrow 
\pi(\textnormal{spt}(\mu))
\,\,\, \textnormal{as subsets of $\rel^3$ 
in Hausdorff distance, as} \,\,\, l\to \infty. 
\end{equation}
Moreover, we note that for any 
smooth curve $c:(-\varepsilon,\varepsilon) 
\longrightarrow \sphere^2$ we have the formula 
$$ 
\langle N^{\sphere^2}(c(t)),c''(t) \rangle 
= II^{\sphere^2}_{c(t)}(c'(t),c'(t)),  \quad  \forall \,
t\in  (-\varepsilon,\varepsilon),
$$
from elementary Differential Geometry, where 
$N^{\sphere^2}$ and $II^{\sphere^2}$ denote the 
Gauss-map and the second fundamental form of the 
standard embedding $\sphere^2 \hookrightarrow \rel^3$, 
respectively. Hence, requiring also that $c'(t)$ has 
length one for $t\in  (-\varepsilon,\varepsilon)$, 
we see that the ``normal component'' 
$\langle N^{\sphere^2}(c(t)),c''(t) \rangle$ 
of the curvature vector
$\vec \kappa^{\rel^3}_{c}(t) \equiv c''(t)$ 
along $c$ - when considered as a path in 
$\rel^3$ - is exactly given by $II^{\sphere^2}_{c(t)}(c'(t),c'(t))$ 
and thus equals $1$, the only principle 
curvature of the standard unit sphere $\sphere^2$. 
We can therefore reformulate 
the elastic energy $\Wil$ of any smooth closed curve 
$\gamma:\sphere^1\longrightarrow \sphere^2$ in the way: 
\begin{equation} \label{Wil.Bil} 
\Wil(\gamma) \equiv \int_{\sphere^1} 
1 + |\vec \kappa^{\sphere^2}_{\gamma}|^2 \,d\mu_{\gamma} 
= \int_{\sphere^1} |\vec \kappa^{\rel^3}_{\gamma}|^2 \,d\mu_{\gamma}, 
\end{equation}  
which is simply the standard elastic energy of a smooth 
closed curve in $\rel^3$.   
Combining now formula \eqref{Wil.Bil} with 
the bounds \eqref{Length.total} and
\eqref{L2.bounded}, again Allard's compactness 
theorem implies that the integral $1$-varifolds $\nu_l:=\Hnm\lfloor_{\textnormal{trace}(\gamma_{t_{j_l}})}$ 
converge weakly - up to extraction of another 
subsequence - to an integral $1$-varifold $\nu$ 
in $\rel^3$:    
\begin{equation} \label{weak.convergence.nu} 
\Hnm\lfloor_{\textnormal{trace}(\gamma_{t_{j_l}})} 
\longrightarrow \nu  \quad 
\textnormal{weakly as Radon measures on} \,\,\, \rel^3,  
\end{equation}  
as $l\to \infty$. Now we also know that the 
$1$-dimensional Hausdorff-densities 
$\theta^1(\nu_l)$ exist in every point of $\rel^3$ 
and satisfy $\theta^1(\nu_l)\geq 1$ pointwise on 
$\textnormal{spt}(\nu_l)=
\textnormal{trace}(\gamma_{t_{j_l}})$, 
since $\gamma_{t_{j_l}}$ are closed and smooth curves.
Hence, combining this with convergence \eqref{weak.convergence.nu}, formula \eqref{Wil.Bil}, 
and with estimates \eqref{Length.total}, \eqref{L2.bounded}, 
we can easily check that the conditions of Proposition \ref{14.7} below are satisfied by the integral $1$-varifolds  $\nu_l=\Hnm\lfloor_{\textnormal{trace}(\gamma_{t_{j_l}})}$  
with $n=1$, $m=2$, $\alpha=\frac{1}{2}$ and  
$\beta=\frac{1}{4}$. Combining the statement of 
Proposition \ref{14.7} with convergence \eqref{Hausdorff.convergence.2} we obtain first of all:
\begin{eqnarray} \label{convergence.supports}
\pi(\textnormal{spt}(\mu)) \longleftarrow 
\textnormal{trace}(\gamma_{t_{j_l}}) 
= \textnormal{spt}
(\Hnm\lfloor_{\textnormal{trace}(\gamma_{t_{j_l}})}) 
\longrightarrow \textnormal{spt}(\nu) \\
\textnormal{as subsets of $\rel^3$ 
in Hausdorff distance, as} \,\,\, l\to \infty, \nonumber
\end{eqnarray} 
which implies in particular:  
$\pi(\textnormal{spt}(\mu))=\textnormal{spt}(\nu)$, 
and we infer furthermore: 
\begin{equation}   \label{spt.nu.nu.l} 
\textnormal{spt}(\nu) = 
\{x \in \rel^3 \,|\,\forall\, l \in \nat \,\,\exists\, 
x_l \in \textnormal{spt}(\nu_l) \,\,
\textnormal{such that} \,\, x_l \longrightarrow x \, \}.
\end{equation}
Since $\textnormal{spt}(\mu)$ is already known 
to be either an embedded $2$-sphere or an embedded 
compact torus, the equation $\pi(\textnormal{spt}(\mu))=\textnormal{spt}(\nu)$
particularly shows that $\textnormal{spt}(\nu)$ 
is a compact and path-connected subset of $\sphere^2$. 
Moreover, since $\nu$ is a $1$-rectifiable varifold, 
one can easily derive from Theorem 3.2 in \cite{Simon.1984}
that $\nu$ coincides with the measure 
$\theta^{*1}(\nu)\cdot \Hnm\lfloor_{[\theta^{*1}(\nu)>0]}$
on entire $\rel^3$ and that $[\theta^{*1}(\nu)>0]$ is  
a countably $1$-rectifiable subset of $\rel^3$, 
where ``$\theta^{*1}(\nu)$'' denotes the upper 
$1$-dimensional Hausdorff-density of $\nu$; compare with 
Paragraph 3 in \cite{Simon.1984}. 
Since $\nu$ is here additionally integral, we 
therefore infer especially: 
\begin{eqnarray}    \label{Hnm.nu.A} 
\Hnm(A) = \int_{A} (\theta^{*1}(\nu))^{-1} 
\cdot \theta^{*1}(\nu) \, d\Hnm 
= \int_{A} (\theta^{*1}(\nu))^{-1} \, d\nu 
\leq \nu(A) <\infty                      \\ 
\textnormal{for all} \,\,\,
\Hnm\textnormal{-measurable subsets 
A of} \,\, [\theta^{*1}(\nu)>0],      \nonumber
\end{eqnarray} 
recalling that $\nu$ is especially 
a Radon measure on $\rel^3$, and that here the upper $1$-dimensional density $\theta^{*1}(\nu)$ has to satisfy
$\theta^{*1}(\nu)\geq 1$\, $\nu$-almost everywhere 
on $\rel^3$. Moreover, again on account of 
convergence \eqref{weak.convergence.nu}, 
formula \eqref{Wil.Bil}, and estimates \eqref{Length.total}, \eqref{L2.bounded}, the conditions of Proposition \ref{14.5} are satisfied by  $\nu_l=\Hnm\lfloor_{\textnormal{trace}(\gamma_{t_{j_l}})}$  
with $n=1$, $m=2$, $\alpha=\frac{1}{2}$ and  
$\beta=\frac{1}{4}$, and we can conclude 
together with equation \eqref{spt.nu.nu.l}: 
\begin{eqnarray} \label{upper.semicontinuity}
\theta^1(\nu,x) \,\, \textnormal{exists and} \,\, 
\theta^1(\nu,x) \geq \limsup_{l\to \infty} \theta^1(\nu_l,x_l) \geq 1 \,\,\, \textnormal{for every}\, x \in \textnormal{spt}(\nu),
\end{eqnarray}
where we have chosen an arbitrary point 
$x\in \textnormal{spt}(\nu)$ and an appropriate sequence 
$x_l \in \textnormal{spt}(\nu_l)$ with $x_l \to x$ in $\rel^3$, according to equation \eqref{spt.nu.nu.l}, and 
where we have again used the obvious fact that $\theta^1(\nu_l)\geq 1$ pointwise on 
$\textnormal{spt}(\nu_l)=\textnormal{trace}(\gamma_{t_{j_l}})$ for every $l \in \nat$. Combining the obvious general fact that $[\theta^{*1}(\xi)>0]$ is contained in $\textnormal{spt}(\xi)$, for every rectifiable $1$-varifold $\xi$, with statement \eqref{upper.semicontinuity}, 
we finally obtain:
\begin{eqnarray} \label{sets.included}
[\theta^{1}(\nu)\geq 1] \subseteq 
[\theta^{*1}(\nu)>0] \subseteq  
\textnormal{spt}(\nu)\subseteq [\theta^{1}(\nu)\geq 1],  
\end{eqnarray}
proving that these three subsets of $\sphere^2$ 
coincide with each other. Since the set $[\theta^{*1}(\nu)>0]$ is already known to be countably $1$-rectifiable, 
statement \eqref{sets.included} proves in particular that the 
compact set $\pi(\textnormal{spt}(\mu))=\textnormal{spt}(\nu)$ 
is actually a countably $1$-rectifiable 
subset of $\sphere^2$, which additionally has to have 
finite $\Hnm$-measure on account of \eqref{Hnm.nu.A},     
simply taking here $A=\textnormal{spt}(\nu)=[\theta^{*1}(\nu)>0]$.
This shows especially that $\textnormal{spt}(\nu)$ 
cannot be a dense subset of $\sphere^2$, because otherwise 
there would hold $\textnormal{spt}(\nu)=\overline{\textnormal{spt}(\nu)}
=\sphere^2$ and therefore $\Hn(\textnormal{spt}(\nu))=4\pi$, 
which obviously contradicts  $\Hnm(\textnormal{spt}(\nu))<\infty$. 
Hence, there has to be some point $x_0\in \sphere^2$ 
and some radius $\varrho>0$ such that
$\textnormal{spt}(\nu)\subset 
\sphere^2\setminus B_{2\varrho}^3(x_0)$.  
On account of the equivariance of the degenerate 
elastic energy flow \eqref{elastic.energy.flow} with respect to rotations of $\sphere^2$ we may assume that here $x_0$ is exactly the north pole $(0,0,1)$ on $\sphere^2$. 
Moreover, from convergence \eqref{convergence.supports} 
we can thus infer that the converging sets $\textnormal{trace}(\gamma_{t_{j_l}})$ still have to be 
contained in $\sphere^2\setminus B^3_{\varrho}(x_0)$ 
for sufficiently large $l \in \nat$, say 
for every $l \in \nat$ without loss of generality. 
Hence, we can apply here stereographic projection 
$\Stereo:\sphere^2\setminus \{x_0\} \stackrel{\cong}\longrightarrow \rel^2$, 
$(x,y,z) \mapsto \frac{1}{1-z}(x,y)$, 
and thus map all the sets $\textnormal{trace}(\gamma_{t_{j_l}})$   
stereographically onto closed planar curves 
with smooth and regular parametrizations 
$\hat \gamma_{t_{j_l}}:=\Stereo(\gamma_{t_{j_l}})$,
which have to be contained in some compact subset $K=K(\varrho)$ of $\rel^2$. 
Now, given any smooth closed curve 
$\gamma:\sphere^1 \longrightarrow \sphere^2 \setminus B^3_{\varrho}(x_0)$ we can compare the pointwise values 
of its Euclidean curvature 
$|\vec \kappa^{\rel^3}_{\gamma}|$ with the corresponding 
values of the Euclidean curvature 
$|\vec \kappa^{\rel^2}_{\Stereo(\gamma)}|$ of 
the stereographically projected curve 
$\Stereo(\gamma):\sphere^1 \longrightarrow \rel^2$. 
Hence, there is a constant $C=C(\varrho)>0$, being 
independent of $\gamma$, such that:
\begin{equation} \label{estimate.curv} 
|\vec \kappa^{\rel^2}_{\Stereo(\gamma)}| 
\leq C(\varrho) \, |\vec \kappa^{\rel^3}_{\gamma}| 
\quad \textnormal{pointwise on} \,\,\, \sphere^1, 
\end{equation}   
provided there holds $\gamma:\sphere^1 \longrightarrow \sphere^2 \setminus B^3_{\varrho}(x_0)$. 
Recalling now formula \eqref{Wil.Bil} and 
noting that there also holds  
\begin{equation} \label{estimate.tangent} 
|\partial_x(\Stereo(\gamma))| 
\leq \tilde C(\varrho) \, |\partial_x\gamma| 
\quad \textnormal{pointwise on} \,\,\, \sphere^1, 
\end{equation}  
for any smooth closed curve 
$\gamma:\sphere^1 \longrightarrow \sphere^2 \setminus B^3_{\varrho}(x_0)$, we finally arrive at the estimates: 
\begin{equation} \label{Wil.up.down} 
\int_{\sphere^1} |\vec \kappa^{\rel^2}_{\Stereo(\gamma)}|^2 \, 
d\mu_{\Stereo(\gamma)}              
\leq C^2(\varrho) \,\tilde C(\varrho)\,
\int_{\sphere^1} |\vec \kappa^{\rel^3}_{\gamma}|^2 \, 
d\mu_{\gamma}     
= C^2(\varrho) \,\tilde C(\varrho)\,\Wil(\gamma) 
\end{equation}  
and similarly 
\begin{equation} \label{Length.up.down} 
	\textnormal{length}(\Stereo(\gamma)) =   
	\int_{\sphere^1} \,1 \, d\mu_{\Stereo(\gamma)} 
	\leq \tilde C(\varrho)\,
	\int_{\sphere^1}  \,1 \, d\mu_{\gamma} 
	\leq \tilde C(\varrho)\, \Wil(\gamma), 
\end{equation}   
for any smooth closed curve 
$\gamma:\sphere^1 \longrightarrow \sphere^2 \setminus B^3_{\varrho}(x_0)$. 
Recalling now that we could guarantee above 
that $\textnormal{trace}(\gamma_{t_{j_l}})$ is actually contained in $\sphere^2 \setminus B^3_{\varrho}(x_0)$ for every $l \in \nat$ and that we derived above 
from formula \eqref{Will.Wil} and from Proposition \ref{L2.estimates.1} that the curves $\gamma_t$ - being driven by flow \eqref{elastic.energy.flow} - satisfy:
\begin{equation}  \label{Wil.8} 
\Wil(\gamma_t)=\int_{\sphere^1} 1+|\kappa_{\gamma_t}^{\sphere^2}|^2 \, d\mu_{\gamma_t} 
< \Wil(\gamma_0)\leq 8,  \quad \forall \, t \in [0,T_{\textnormal{max}}(\gamma_0)),
\end{equation}  
we can combine estimates \eqref{Wil.up.down}, \eqref{Length.up.down} and \eqref{Wil.8} and arrive at the estimates: 
\begin{eqnarray}  \label{Wil.down.8} 
\int_{\sphere^1} |\kappa^{\rel^2}_{\Stereo(\gamma_{t_{j_l}})}|^2 \, d\mu_{\Stereo(\gamma_{t_{j_l}})}              
\leq 8 \,C^2(\varrho) \,\tilde C(\varrho)    \nonumber \\ 
\textnormal{and} \qquad 
\textnormal{length}(\Stereo(\gamma_{t_{j_l}})) 
\leq 8\, \tilde C(\varrho), 
\end{eqnarray}     
for every $l \in \nat$. Since we also know that 
the traces of the projected curves $\Stereo(\gamma_{t_{j_l}})$ 
are contained in some compact subset $K(\varrho)$ of $\rel^2$, 
we can apply Theorem 3.1 in \cite{Bellettini.1993}
and infer that at least some subsequence of the projected paths $\Stereo(\gamma_{t_{j_l}})$ - up to smooth reparametrization - converges weakly in 
$W^{2,2}(\sphere^1,\rel^2)$ and thus - by compact 
Morrey-embedding 
$W^{2,2}(\sphere^1)\hookrightarrow C^{1}(\sphere^1)$ - strongly in $C^1(\sphere^1,\rel^2)$ to some closed 
limit curve $\gamma^*:\sphere^1\longrightarrow \rel^2$,
whose trace has to be contained again in the compact 
subset $K(\varrho)$ of $\rel^2$.
In combination with convergence \eqref{convergence.supports} we conclude that the closed path $\Stereo^{-1}(\gamma^*)$ has to parametrize $\textnormal{spt}(\nu)$ in $\sphere^2$, 
which proves that $\textnormal{spt}(\nu)$ and thus 
$\pi(\textnormal{spt}(\mu))$ is not only 
a compact and path-connected, countably $1$-rectifiable 
subset of $\sphere^2$, 
but actually the trace of a closed $C^1$-curve in $\sphere^2$.
Additionally we recall here that the embedded curves 
$\gamma_{t_{j_l}}$ bound simply connected domains $\Omega_{t_{j_l}}$ in $\sphere^2$ whose areas are bounded from below by the positive number $2(\pi-2)$ and from above 
by the number $2(\pi+2)<4 \pi$, for every $l \in \nat$. 
Hence, combining this again with \eqref{convergence.supports} we infer that the closed limit $C^1$-curve $\textnormal{trace}(\Stereo^{-1}(\gamma^*))= \pi(\textnormal{spt}(\mu))$ has to be the topological 
boundary of some compact and connected subset $B$ 
of $\sphere^2$ with measure $0<\Hn(B)<4 \pi$, ruling 
out that $\pi(\textnormal{spt}(\mu))$ might  
only be a point in $\sphere^2$ or homeomorphic 
to some compact interval.
Still we do not know, whether the closed $C^1$-path $\Stereo^{-1}(\gamma^*)$ parametrizing 
$\pi(\textnormal{spt}(\mu))$ is regular, and 
we cannot rule out neither at this point, whether $\pi(\textnormal{spt}(\mu))$ might have ``cusps'', 
which prevents us from reparametrizing $\Stereo^{-1}(\gamma^*)$ into a regular closed curve.  
We therefore do not know at this point, whether the preimage 
$\pi^{-1}(\pi(\textnormal{spt}(\mu)))$ is a  
$C^1$-Hopf-torus in $\sphere^3$, in the sense of 
Definition \ref{Hopf.torus.immersion}. 
But still we know that $\pi^{-1}(\pi(\textnormal{spt}(\mu)))$ is a compact and path-connected topological space - 
being endowed with the relative topology of $\sphere^3$.   \\  
On the other hand, $\textnormal{spt}(\mu)$ is already known 
to be either an embedded $2$-sphere or an embedded 
compact torus in $\sphere^3$. We shall prove now that 
the first case is topologically impossible. 
We assume by contradiction that 
$\textnormal{spt}(\mu)$ would be an embedded 
$2$-sphere in $\sphere^3$, and we recall here that 
the Hopf-fibration is actually the Serre-fibration 
$\sphere^1 \hookrightarrow \sphere^3 \stackrel{\pi}\longrightarrow \sphere^2$, 
whose restriction to the subset 
$\pi^{-1}(\pi(\textnormal{spt}(\mu)))$ of $\sphere^3$ 
is a Serre-fibration as well, i.e. we have as well: 
\begin{equation} \label{Serre.fibration}
\sphere^1 \hookrightarrow \pi^{-1}(\pi(\textnormal{spt}(\mu))) 
\stackrel{\pi}\longrightarrow \pi(\textnormal{spt}(\mu)).
\end{equation}
Now, by assumption the $2$-sphere $\textnormal{spt}(\mu)$ would have to be contained in the preimage $\pi^{-1}(\pi(\textnormal{spt}(\mu)))$, 
with the additional property that the fiber 
of $\pi$ over any chosen $x \in \pi(\textnormal{spt}(\mu))$
must have non-empty intersection with $\textnormal{spt}(\mu)$. 
Moreover, if for some arbitrarily chosen 
$x \in \pi(\textnormal{spt}(\mu))$ the intersection 
$\textnormal{spt}(\mu)\cap \pi^{-1}(\{x\})$ would not 
be open in $\pi^{-1}(\{x\})$, then it had a boundary 
point $Z$ on the great circle $\pi^{-1}(\{x\})$, 
i.e. there was some 
$Z \in \pi^{-1}(\{x\}) \cap \textnormal{spt}(\mu)$ 
such that every ball $B^4_r(Z)$ still intersects 
$\pi^{-1}(\{x\}) \setminus \textnormal{spt}(\mu)$.  
Since we assume that $\textnormal{spt}(\mu)$ is 
entirely contained in $\pi^{-1}(\pi(\textnormal{spt}(\mu)))$, 
there would have to be another fiber $\pi^{-1}(\{y\})$, 
with $y \in \pi(\textnormal{spt}(\mu))$ and $y\not = x$,
which has non-empty intersection with $\textnormal{spt}(\mu)$ 
and intersects the fiber $\pi^{-1}(\{x\})$ at the 
boundary point $Z$ of $\pi^{-1}(\{x\}) \cap \textnormal{spt}(\mu)$ in $\pi^{-1}(\{x\})$. 
Otherwise, $Z$ would either have to be a point 
on the manifold-boundary $\partial (\textnormal{spt}(\mu))$ 
of the surface $\textnormal{spt}(\mu)$ itself, which is obviously empty here, or $\textnormal{spt}(\mu)$ would not be 
covered by the union of fibers over $\textnormal{spt}(\mu)$, 
recalling here that \eqref{Serre.fibration} especially 
implies that  
\begin{equation} \label{fibers.fibration}
\pi^{-1}(\pi(\textnormal{spt}(\mu)))
=\bigcup_{y\in \pi(\textnormal{spt}(\mu))} \pi^{-1}(\{y\}) 
\end{equation}   
is a disjoint union of great circles - in particular  
of closed embedded curves - in $\sphere^3$.  
Hence, this particular point $Z\in \textnormal{spt}(\mu)$ would satisfy: $Z \in \pi^{-1}(\{x\}) \cap \pi^{-1}(\{y\})$, 
for $y \not = x$, although each pair of distinct fibers 
in \eqref{fibers.fibration} has empty intersection.   
Hence, for every $x \in \pi(\textnormal{spt}(\mu))$ the intersection $\pi^{-1}(\{x\}) \cap \textnormal{spt}(\mu)$ would have to be a non-empty, open and - clearly - also closed subset of the great circle $\pi^{-1}(\{x\})$. Since
each $\pi^{-1}(\{x\})$ is a connected space, we would arrive at the equality 
$$ 
\pi^{-1}(\{x\}) \cap \textnormal{spt}(\mu) = \pi^{-1}(\{x\})
\quad \forall \, x \in \pi(\textnormal{spt}(\mu)),  
$$
or equivalently 
$\pi^{-1}(\{x\}) \subseteq \textnormal{spt}(\mu)$ 
for all $x \in \pi(\textnormal{spt}(\mu))$.  
On account of \eqref{fibers.fibration} this result implies that $\textnormal{spt}(\mu)$ would contain 
$\pi^{-1}(\pi(\textnormal{spt}(\mu)))$, 
and therefore finally:
\begin{equation} \label{equality.of.surafces} 
\textnormal{spt}(\mu) 
=\pi^{-1}(\pi(\textnormal{spt}(\mu))),
\end{equation} 
which is wrong in the considered first case, 
because on the one hand $\textnormal{spt}(\mu)$ 
is simply connected, whereas on the other hand $\pi^{-1}(\pi(\textnormal{spt}(\mu)))$ cannot  
have trivial fundamental group. The latter assertion 
is intuitively totally clear, and technically it follows rather quickly from the long homotopy sequence of 
the Serre-fibration in \eqref{Serre.fibration}: 
\begin{eqnarray*} 
\stackrel{\partial_*}\longrightarrow \pi_2(\sphere^1) \stackrel{i_*}\longrightarrow  \pi_2(\pi^{-1}(\pi(\textnormal{spt}(\mu))))
\stackrel{\pi_*}\longrightarrow \pi_2(\pi(\textnormal{spt}(\mu)))  \stackrel{\partial_*}\longrightarrow   \\
\stackrel{\partial_*}\longrightarrow 
\pi_1(\sphere^1)                      \stackrel{i_*}\longrightarrow             
\pi_1(\pi^{-1}(\pi(\textnormal{spt}(\mu))))   \stackrel{\pi_*}\longrightarrow \pi_1(\pi(\textnormal{spt}(\mu)))
\stackrel{\partial_*}\longrightarrow \pi_0(\sphere^1),
\end{eqnarray*}  
which actually simplifies to the short exact sequence 
\begin{equation} \label{short.exact.sequence}  
0 = \pi_2(\pi(\textnormal{spt}(\mu)))   
\stackrel{\partial_*}\longrightarrow 
\ganz \stackrel{i_*}\longrightarrow  
\pi_1(\pi^{-1}(\pi(\textnormal{spt}(\mu))))   \stackrel{\pi_*}\longrightarrow \pi_1(\pi(\textnormal{spt}(\mu)))
\stackrel{\partial_*}\longrightarrow \pi_0(\sphere^1)=0.  
\end{equation} 
In this last step we used the important fact that  $\pi(\textnormal{spt}(\mu))$ can be parametrized  
by a closed $C^1$-curve, implying that 
$\pi(\textnormal{spt}(\mu))$ is not only   
path-connected and $1$-rectifiable, but that
there is a universal covering  
$c:\rel \twoheadrightarrow \pi(\textnormal{spt}(\mu))$, 
namely the composition $c:= \Stereo^{-1}(\gamma^*) \circ p$, 
where $p:\rel \twoheadrightarrow \sphere^1$ denotes 
the canonical universal covering of $\sphere^1$.  
Therefore, any continuous map 
$g:\sphere^2 \longrightarrow \pi(\textnormal{spt}(\mu))$ 
can be lifted against the above covering map $c$ 
to a continuous map 
$\tilde g:\sphere^2 \longrightarrow \rel$, 
satisfying $c \circ \tilde g = g$. Obviously, since  
$\tilde g$ is homotopic to any prescribed constant map 
from $\sphere^2$ into $\rel$, the composition 
$c \circ \tilde g$ is homotopic to some constant map,  
such that any chosen base points can be preserved 
throughout the chosen homotopy. 
This proves already that $\pi_2(\pi(\textnormal{spt}(\mu)))$   
is trivial, as intuitively expected. 
Now, the above short exact sequence 
\eqref{short.exact.sequence} tells us in particular 
that the induced homomorphism 
$\ganz \stackrel{i_*}\longrightarrow  
\pi_1(\pi^{-1}(\pi(\textnormal{spt}(\mu))))$ is injective, 
proving that $\pi_1(\pi^{-1}(\pi(\textnormal{spt}(\mu))))$  
indeed cannot be trivial. \\
Hence, indeed the second case must hold here in which 
$\textnormal{spt}(\mu)$ is some \emph{embedded compact} torus 
in $\sphere^3$. Since such a surface obviously has no  
manifold-boundary neither, we can argue exactly 
as we did in the first case, in order to
infer the equality in \eqref{equality.of.surafces}
from the obvious inclusion
$\textnormal{spt}(\mu) \subseteq \pi^{-1}(\pi(\textnormal{spt}(\mu)))$, as well in this second 
case. However, this result can only be true, if $\pi(\textnormal{spt}(\mu))$ is an embedded closed  
curve in $\sphere^2$ - exactly what we wanted to know above.  
Hence, the compact surface $\textnormal{spt}(\mu)$ 
does not only turn out to be some embedded torus 
in $\sphere^3$, it actually turns out to be 
an embedded $C^1$-Hopf-torus in the sense of 
Definition \ref{Hopf.torus.immersion}, whose profile curve 
is an embedded closed $C^1$-curve in $\sphere^2$
which can actually be parametrized 
by the closed $C^{1}$-path $\Stereo^{-1}(\gamma^*)$, 
as we had figured out above. \\    
In particular, the sequence of embeddings $\{F_{t_{j_l}}\}$ satisfies all requirements of the second part of Theorem \ref{limit.MIWF}, implying that all statements of the second part of Theorem \ref{limit.MIWF} have to hold for 
the considered sequence of embeddings $\{F_{t_{j_l}}\}$
and for their non-trivial limit varifold $\mu$ from 
\eqref{weak.convergence.F.j}.
We can therefore infer from statement 
\eqref{parametrization.mu} that the embedded Hopf-torus $\textnormal{spt}(\mu)$ possesses a uniformly 
conformal bi-Lipschitz parametrization
$f:(\Sigma,g_{\textnormal{poin}}) \stackrel{\cong}\longrightarrow 
\textnormal{spt}(\mu)$, for some zero scalar curvature and unit volume metric $g_{\textnormal{poin}}$ on $\Sigma$, with 
conformal factor $u\in L^{\infty}(\Sigma)$ bounded  
by some suitable constant $\Lambda$ depending on the sequence 
$\{F_{t_{j_l}}\}$ from \eqref{weak.convergence.F.j} 
and on $\mu$ - as already explained in the proof of 
the second part of Theorem \ref{limit.MIWF} - 
and with $\Will(f)=\Will(\mu) < 8\pi$ 
on account of formulae \eqref{Willmore.energy} and 
\eqref{Willmore.mu}. 
\item[2)] We consider here again an arbitrary sequence  
$t_j \nearrow T_{\textnormal{max}}(F_0)$. 
We recall from the proof of the first part 
of this theorem that $\Will(F_t)<8 \pi$ must hold 
for every $t\in (0,T_{\textnormal{max}}(F_0))$ 
and that there has to be some suitable subsequence 
$\{t_{j_l}\}$ of $\{t_j\}$ for which the varifolds 
$\Hn\lfloor_{F_{t_{j_l}}(\Sigma)}$ converge weakly 
as in \eqref{weak.convergence.F.j} to some 
non-trivial integral $2$-varifold $\mu$,  
whose support is an embedded Hopf-torus in $\sphere^3$. 
In particular, all statements of the second part of Theorem \ref{limit.MIWF} are valid for the appropriately chosen subsequence $\{F_{t_{j_l}}\}$, and we can thus especially infer that the immersions $F_{t_{j_l}}$ can be reparametrized by smooth diffeomorphisms 
$\Phi_{j_l}:\Sigma \stackrel{\cong}\longrightarrow \Sigma$ in such a way that the reparametrizations 
$\tilde F_{t_{j_l}}:=F_{t_{j_l}} \circ \Phi_{j_l}$  
are uniformly conformal with respect to certain metrics 
$g_{\textnormal{poin},j_l}$ of vanishing scalar curvature and 
with smooth conformal factors $u_{j_l}$ which are uniformly bounded in $L^{\infty}(\Sigma,g_{\textnormal{poin}})$ and in $W^{1,2}(\Sigma,g_{\textnormal{poin}})$, i.e. with:
\begin{eqnarray} \label{bounded.factors.1} 
\parallel u_{j_l} \parallel_{L^{\infty}(\Sigma,g_{\textnormal{poin}})} 
\leq \Lambda,  \quad \textnormal{for every} \,\,\, l\in \nat,   \\
\textnormal{and} \,\,  
\parallel u_{j_l} \parallel_{W^{1,2}(\Sigma,g_{\textnormal{poin}})} 
\leq C(\Lambda), \quad \textnormal{for every} \,\,\, l\in \nat. 
\label{bounded.factors.2}  
\end{eqnarray} 
Here, again the constant $\Lambda$ does not only depend 
on the limit varifold $\mu$ but also on local geometric properties of the embedded surfaces $F_{t_{j_l}}(\Sigma)\subset \sphere^3$ appearing in \eqref{weak.convergence.F.j} 
- as already pointed out below formula \eqref{bounded.factor}
in the proof of the second part of Theorem \ref{limit.MIWF} - and $g_{\textnormal{poin}}$ is a certain zero scalar curvature and unit volume metric which satisfies 
on account of \eqref{smooth.converg.metrics.mod} 
- again up to extraction of a subsequence: 
\begin{equation} \label{smooth.converg.metrics.2} 
g_{\textnormal{poin},j_l}
\longrightarrow g_{\textnormal{poin}} \quad 
\textnormal{smoothly as} \,\,l\to \infty,
\end{equation}  
as explained below formula \eqref{smooth.converg.metrics}.
Hence, we also infer the important estimates 
\begin{equation}   \label{bounded.f.tjs} 
\parallel \tilde F_{t_{j_l}} 
\parallel_{W^{1,\infty}(\Sigma,g_{\textnormal{poin}})}
+\parallel \tilde F_{t_{j_l}} 
\parallel_{W^{2,2}(\Sigma,g_{\textnormal{poin}})}
\leq \textnormal{Const}(\Lambda),
\,\,\, \textnormal{for every} \,\,\, l\in \nat,        
\end{equation}
from the proof of the second part of Theorem \ref{limit.MIWF}, implying the existence of a particular subsequence 
$\{\tilde F_{t_{j_k}}\}$ of $\{\tilde F_{t_{j_l}}\}$ which 
converges in the three senses \eqref{weak.convergence.2},  
\eqref{W.1.infty.convergence} and \eqref{C0.convergence}
to the uniformly conformal bi-Lipschitz homeomorphism 
$f:(\Sigma,g_{\textnormal{poin}}) \stackrel{\cong}\longrightarrow 
\textnormal{spt}(\mu)$ from the first part of this theorem. 
Since we aim to prove the desired estimate 
\eqref{L.infty.W.4.2.estimate} for this particular 
subsequence $\{\tilde F_{t_{j_k}}\}$ of 
$\{\tilde F_{t_{j_l}}\}$, 
we should relabel here this subsequence 
$\{\tilde F_{t_{j_k}}\}$ again into $\{\tilde F_{t_{j}}\}$, 
just for ease of notation.  
Now combining estimates \eqref{bounded.factors.1}, 
\eqref{bounded.factors.2} and \eqref{bounded.f.tjs} with convergence \eqref{smooth.converg.metrics.2},
we can proceed exactly as in the proof of Theorem 4.1 and 
of Proposition 5.2 in \cite{Palmurella.2022}, in order to prove estimate \eqref{L.infty.W.4.2.estimate.local} below 
and then finally also the desired estimate \eqref{L.infty.W.4.2.estimate}, for every \,$j\in \nat$. 
First of all, comparing the assumptions of Theorem 4.1 in 
\cite{Palmurella.2022} with our knowledge in formulae 
\eqref{bounded.factors.1}--\eqref{bounded.f.tjs} 
it is worth mentioning that we can exchange our estimate \eqref{bounded.factors.2} by the slightly weaker estimate \eqref{Lorentz.estimate} below, 
which also follows from the equations
\begin{equation} \label{Gauss.curvature.2}
-\triangle_{g_{\textnormal{poin},j}}(u_{j})    
= e^{2 u_{j}} \, K_{\tilde F_{t_{j}}^*(g_{\textnormal{euc}})} 
\quad \textnormal{on} \,\,\, \Sigma
\end{equation}  
in \eqref{Gauss.curvature}, together with the elementary estimates
$$ 
\int_{\Sigma} |K_{\tilde F_{t_{j}}^*(g_{\textnormal{euc}})}| 
\,d\mu_{{\tilde F_{t_{j}}^*(g_{\textnormal{euc}})}} 
\leq \frac{1}{2} \,\int_{\Sigma} 
|A_{\tilde F_{t_{j}}^*(g_{\textnormal{euc}})}|^2  \,
d\mu_{{\tilde F_{t_{j}}^*(g_{\textnormal{euc}})}} 
= 2 \, \Will(F_{t_{j}}) <16\,\pi,   \,\,\,  
\forall \,j\in \nat,
$$ 
and with the fifth part of Proposition \ref{L2.estimates.1} or also directly with convergence  \eqref{smooth.converg.metrics.2}, as exactly pointed out in the proof of Theorem 5.4 in \cite{Riviere.Park.City.2013}
\footnote{It should be stressed here that estimate \eqref{Lorentz.estimate} can actually be proved without 
	the strong $L^{\infty}(\Sigma)$-estimate of the conformal 
	factors in \eqref{bounded.factors.1}, but still estimates \eqref{bounded.factors.1} respectively \eqref{uniformly.conformal.R3} are indispensable for this entire proof because of the importance of the closely related estimates \eqref{bounded.f.tjs} respectively \eqref{bounded.f.tjs.new} regarding the basic estimate \eqref{second.fundam.locally.bounded} below, and additionally because of the unavoidable application of 
	estimate \eqref{uniformly.conformal.R3} in 
	the end of this proof, adopting the proofs of 
	Theorem 4.1 and Proposition 5.2 of \cite{Palmurella.2022}.}: 
\begin{equation} \label{Lorentz.estimate}
	\parallel \nabla^{g_{\textnormal{poin}}}(u_{j}) \parallel_{L^{2,\infty}(\Sigma)}    
	\leq C(g_{\textnormal{poin}}) \, \Will(F_{t_{j}}) 
	\leq C(g_{\textnormal{poin}}) \, \Will(F_{0}), \,\,\,  \forall \,j\in \nat.
\end{equation}  
Here, ``$L^{2,\infty}(\Sigma)$'' is classical terminology 
denoting some particular \emph{Lorentz space} 
which especially satisfies
$L^{2}(\Sigma,g_{\textnormal{poin}}) \hookrightarrow L^{2,\infty}(\Sigma,g_{\textnormal{poin}})$;
see for example \cite{Tartar.1998} or Section 3.3 in \cite{Helein.2004}. 
Now we use our hypothesis on the mean curvature vectors 
$\vec H_{F_{t_j},\sphere^3}$ - 
which we shall abbreviate here by $\vec H_{F_{t_j}}$ - 
of the embeddings $F_{t_j}:\Sigma \longrightarrow \sphere^3$ 
to remain uniformly bounded for all $j \in \nat$, i.e. that there is some large number $K$ such that 
$\parallel \vec H_{F_{t_j}} \parallel_{L^{\infty}(\Sigma)}\leq K$ holds for every $j \in \nat$. 
Since in our considered situation every immersion $F_{t}$ parametrizes some Hopf-torus, 
for $t\in [0,T_{\textnormal{max}}(F_0))$, we can combine our 
uniform bound on $\parallel \vec H_{F_{t_j}} \parallel_{L^{\infty}(\Sigma)}$ 
with formulae \eqref{second.fundam.} and \eqref{H} 
in Proposition \ref{Hopf.Willmore.prop} and infer that the entire second fundamental forms $A_{F_{t_j},\sphere^3}$ - which we shall abbreviate here by $A_{F_{t_j}}$ - 
of the embeddings $F_{t_j}:\Sigma \longrightarrow \sphere^3$ can be uniformly bounded:
\begin{eqnarray}  \label{second.fundam.bounded.0} 
\parallel |A_{F_{t_j}}|^2 
\parallel_{L^{\infty}(\Sigma)} 
\equiv \parallel g^{ih}_{F_{t_j}} g^{kl}_{F_{t_j}}\, 
\langle (A_{F_{t_j}})_{ik},(A_{F_{t_j}})_{hl} \rangle_{\rel^4}  
\parallel_{L^{\infty}(\Sigma)} \leq K^2+2 
\end{eqnarray}
for every $j \in \nat$. Since these scalars are invariant with respect to smooth reparametrizations of the embeddings $F_{t_j}$, statement \eqref{second.fundam.bounded.0} 
implies that also 
\begin{eqnarray}  \label{second.fundam.bounded} 
	\parallel |A_{\tilde F_{t_{j}}}|^2 
	\parallel_{L^{\infty}(\Sigma)} 
	\equiv \parallel g^{ih}_{\tilde F_{t_{j}}} 
	g^{kl}_{\tilde F_{t_{j}}}\, 
	\langle (A_{\tilde F_{t_{j}}})_{ik},(A_{\tilde F_{t_{j}}})_{hl} \rangle_{\rel^4}  
	\parallel_{L^{\infty}(\Sigma)} \leq K^2+2, 
\end{eqnarray}
for every $j\in \nat$.  
Now, on account of the uniform convergence \eqref{C0.convergence}
and on account of the conformal invariance of the flow  
\eqref{Moebius.flow}, we can assume - as in \eqref{containoe} - that the images of the sequence 
$\{\tilde F_{t_{j}}\}$ are contained   
in $\sphere^3\setminus B^4_{\delta}((0,0,0,1))$ for some 
$\delta>0$. We may therefore apply the stereographic 
projection $\Stereo$ to the entire sequence $\{\tilde F_{t_{j}}\}$ and obtain new embeddings 
$f_j:=\Stereo \circ \tilde F_{t_{j}}:\Sigma \longrightarrow \rel^3$ which are again uniformly conformal, i.e. satisfy 
as in \eqref{uniform.conformal.R3}: 
\begin{eqnarray}  \label{uniformly.conformal.R3}
f_j^*(g_{\textnormal{euc}}) = e^{2\tilde u_j}\, g_{\textnormal{poin},j}
\,\,\, \textnormal{with} \,\,\, 
\parallel \tilde u_j \parallel_{L^{\infty}(\Sigma)} 
\leq \tilde \Lambda \quad \textnormal{and} \\  \label{uniformly.conformal.R3.2}
\parallel \nabla^{g_{\textnormal{poin}}}(\tilde u_{j}) \parallel_{L^{2,\infty}(\Sigma,g_{\textnormal{poin}})} 
\leq \tilde \Lambda \quad 
\textnormal{for every} \,\,\, j\in \nat 
\end{eqnarray} 	
and for some large constant $\tilde \Lambda$, which depends on the originally chosen subsequence $\{F_{t_{j_l}}\}$ satisfying \eqref{weak.convergence.F.j} and on the limit varifold $\mu$ - as the constant $\Lambda$ in \eqref{bounded.factors.1} did - and additionally on $g_{\textnormal{poin}}$ and $\delta$. 
Moreover, we also obtain
the estimates in \eqref{bounded.f.tjs} for the new embeddings 
$f_j:\Sigma \longrightarrow \rel^3$:  
\begin{equation}   \label{bounded.f.tjs.new} 
\parallel f_j \parallel_{W^{1,\infty}(\Sigma,g_{\textnormal{poin}})}
+ \parallel f_{j} \parallel_{W^{2,2}(\Sigma,g_{\textnormal{poin}})}
\leq \textnormal{Const}(\Lambda,\delta),
\quad  \forall \,\, j\in \nat,        
\end{equation}
and also estimate \eqref{second.fundam.bounded} 
for the new embeddings $f_j:\Sigma \longrightarrow \rel^3$: 
\begin{eqnarray}  \label{second.fundam.bounded.2} 
	\parallel |A_{f_j}|^2 
	\parallel_{L^{\infty}(\Sigma)} 
	\equiv \parallel g^{ih}_{f_j} g^{kl}_{f_j}\, 
	\langle (A_{f_j})_{ik},(A_{f_j})_{hl} \rangle_{\rel^3}  
	\parallel_{L^{\infty}(\Sigma)} 
	\leq C(K,\delta) 
\end{eqnarray}
for every $j \in \nat$. Combining now 
\eqref{bounded.f.tjs.new} and \eqref{second.fundam.bounded.2} 
we finally infer that for every small $\varepsilon_0>0$ 
there is some small $R_0>0$, depending on 
$\varepsilon_0$, $K$ and $\delta$, such that 
\begin{eqnarray}  \label{second.fundam.locally.bounded} 
\int_{B^{g_{\textnormal{poin}}}_{R_0}(x_0)} |A_{f_j}|^2 \, d\mu_{{f_{j}^*(g_{\textnormal{euc}})}}
\leq \varepsilon_0, \quad  \forall \,j\in \nat,   
\end{eqnarray}    
{\bf independently of $x_0 \in \Sigma$}, where ``$B^{g_{\textnormal{poin}}}_r(x_0)$'' denotes the open geodesic disc of radius $r$ about the center point $x_0$ 
in $\Sigma$ with respect to the fixed zero scalar curvature metric $g_{\textnormal{poin}}$. 
Now we fix some $x_0\in \Sigma$ arbitrarily,
we consider some small $\varepsilon_0>0$ and some small 
$R_0>0$ as in \eqref{second.fundam.locally.bounded}, and we 
introduce isothermal charts 
$\psi_j:B^{2}_1(0) \stackrel{\cong}\longrightarrow U_j(x_0)
\subseteq B^{g_{\textnormal{poin}}}_{R_0}(x_0)$ with respect to $g_{\textnormal{poin},j}$ on 
open neighborhoods $U_j(x_0)$ of $x_0$ satisfying  
$\psi_j(0)=x_0$, for each $j\in \nat$, 
as already above in \eqref{phi.k}. 
Hence, we consider here again harmonic and 
bounded functions $v_j$ on $B^{2}_1(0)$, such that         
\begin{eqnarray} \label{estimate.v.j}
	\psi_j^*g_{\textnormal{poin},j}= 
	e^{2v_j} \,g_{\textnormal{euc}}  
	\,\,\,\textnormal{on} \,\,\,B^2_1(0), \,\, 
	\textnormal{with} \,\,\,
	\parallel v_j \parallel_{L^{\infty}(B^2_1(0))} \leq 
	C(g_{\textnormal{poin}},R_0) \,\,\, \forall \, j\in \nat, 
	\nonumber \\
	\textnormal{and with} \,\,\, \parallel 
	\nabla^s(v_j) \parallel_{L^{\infty}(B^2_{7/8}(0))} 
    \leq C(g_{\textnormal{poin}},R_0,s) \,\,\,
    \forall \, j\in \nat,\,\,
\end{eqnarray}
and for each fixed $s\in \nat$, 
where we have again used convergence \eqref{smooth.converg.metrics.2}. In particular, 
the compositions $f_j\circ \psi_j:B^2_{1}(0) 
\longrightarrow \rel^3$ are uniformly conformal 
with respect to the Euclidean metric on $B_1^2(0)$: 
\begin{equation}  \label{conformal.composition} 
(f_j\circ \psi_j)^*(g_{\textnormal{euc}}) 
= \psi_j^*(e^{2\tilde u_j}\, g_{\textnormal{poin},j})   
= e^{2\tilde u_j\circ \psi_j} \, \psi_j^*g_{\textnormal{poin},j}
= e^{2\tilde u_j\circ \psi_j+2v_j} \, g_{\textnormal{euc}}
\,\,\, \textnormal{on} \,\, B_1^2(0). 
\end{equation}     
Now, using the isothermal charts $\psi_j:B^{2}_1(0) \stackrel{\cong}\longrightarrow U_j(x_0)$ statement \eqref{second.fundam.locally.bounded} implies: 
\begin{eqnarray}  \label{second.fundam.locally.bounded.2} 
\int_{B^2_{1}(0)} |A_{f_j\circ \psi_j}|^2 \, 
d\mu_{(f_j\circ \psi_j)^*(g_{\textnormal{euc}})}
\leq \varepsilon_0,  \,\,\,  
\textnormal{for every}\,\,\, j\in \nat.     
\end{eqnarray} 
Since we use several different references on gauge theory
in this proof, we mention here Section 5.1 in \cite{Helein.2004}, 
explaining that condition \eqref{second.fundam.locally.bounded.2} 
can also be formulated in terms of the usual 
Euclidean metric and the $\Lno$-measure on 
$B^2_1(0)\subset \rel^2$, taking equation \eqref{conformal.composition} and the conformal 
invariance of the Dirichlet functional into account:
\begin{eqnarray}  \label{second.fundam.locally.bounded.3} 
\int_{B^2_{1}(0)} |A_{f_j\circ \psi_j}|_{g_{\textnormal{euc}}}^2 \, d\Lno 
\leq \varepsilon_0, \,\,\,  \textnormal{for every}\,\,\, 
j\in \nat. 
\end{eqnarray} 
Now we aim at estimating 
$\parallel \nabla(\tilde u_j\circ \psi_j+v_j) \parallel_{L^{2}(B_{r}^2(0))}$, for sufficiently small 
radii $r \in (0,R_0)$, in terms of the controllable quantity $\varepsilon_0$ and the given upper bounds $K$ and 
$\Lambda$, uniformly for $j\in \nat$; 
see \eqref{Cauchy.4} below. 
To this end, we shall firstly follow the proof of Theorem 5.5 in \cite{Riviere.Park.City.2013}, and then we will
combine this approach with estimates \eqref{uniformly.conformal.R3.2}, \eqref{estimate.v.j} and \eqref{second.fundam.locally.bounded.3} and with 
Theorem 2.1 in \cite{Giaquinta.1984}, p. 78, on sufficiently small discs about $0$ in $\rel^2$. First of all, on account of estimate \eqref{second.fundam.locally.bounded.3} we may apply 
\emph{H\'elein's lifting theorem}, Theorem 4.2 in \cite{Riviere.Park.City.2013}, in order to obtain pairs of functions $e_j^1,e_j^2 \in
W^{1,2}(B^2_1(0),\sphere^2)$ which satisfy both: 
\begin{eqnarray}  \label{e.j}
N_j = e_j^1 \times e_j^2 \quad \textnormal{on} \,\,\,
B^2_1(0) \,\,\, \textnormal{and also}     \nonumber  \\
\int_{B^2_1(0)} |\nabla e_j^1|^2 + 
|\nabla e_j^2|^2 \, d\Lno 
\leq C \, \int_{B^2_1(0)} |A_{f_j\circ \psi_j}|^2 \, d\Lno
\leq C \varepsilon_0, \,\,\,  \textnormal{for every}\,\,\, j\in \nat,     \quad  
\end{eqnarray} 
having used here already estimate \eqref{second.fundam.locally.bounded.3}, where 
$N_j:B^2_1(0) \longrightarrow \sphere^2$ denote 
unit normals along the conformal embeddings 
$f_j\circ \psi_j$, and where $C$ is an absolute constant. 
Now we use Theorem 3.8 in \cite{Riviere.Park.City.2013} 
and estimate \eqref{e.j} and infer, 
that the unique weak solutions $\mu_j$
of the Dirichlet boundary value problems:
\begin{eqnarray}  \label{boundarie}
- \triangle_{\textnormal{euc}}(\mu_j) 
= \sum_{k=1}^3 \det(\nabla(e_j^1)_k,\nabla(e_j^2)_k)   	
\equiv \star N_j^*\textnormal{vol}_{\sphere^2}	   
= e^{2\tilde u_j\circ \psi_j+2v_j} \, 
K_{(f_j\circ \psi_j)^*g_{\textnormal{euc}}}  
\,\, \textnormal{on} \,\, B_1^2(0)     \quad            \\ 
\textnormal{and} \quad \mu_j=0 \quad \textnormal{on} \,\,\, \partial B_1^2(0), \quad \nonumber
\end{eqnarray}	
can be estimated in $W^{1,2}(B_1^2(0))\cap L^{\infty}(B_1^2(0))$: 
\begin{equation}    \label{estimate.Riviere}
\parallel \mu_j \parallel_{L^{\infty}(B_1^2(0))}
+ \parallel \mu_j \parallel_{W^{1,2}(B_1^2(0))}
\leq  C \varepsilon_0, \,\,\, \textnormal{for every} 
\,\,\, j\in \nat, 
\end{equation}   
similarly to estimates \eqref{bounded.factors.1} and \eqref{bounded.factors.2}, but estimates \eqref{estimate.Riviere} are local and therefore more precise. 
See here also Theorem 6.1 in \cite{Kuwert.Schaetzle.2012}
and its even more suitable variant in Proposition 5.1 of \cite{Schaetzle.Conf.factor.2013}.  
Now, on account of formulae (2.48), (2.51), (4.10) and (5.21) in \cite{Riviere.Park.City.2013} - compare here also to equations \eqref{Liouville}, \eqref{Gauss.curvature} and \eqref{Gauss.curvature.2} above - the conformal factors 
$\lambda_j =\tilde u_j\circ \psi_j+v_j$ of the 
conformal embeddings $f_j \circ \psi_j$ satisfy 
exactly equation \eqref{boundarie} on $B_1^2(0)$, 
implying that the differences $\lambda_j-\mu_j$ are
real-valued, harmonic functions on $B_1^2(0)$, for every 
$j\in \nat$. Hence, as in the proof of Theorem 5.5 in \cite{Riviere.Park.City.2013} we can 
combine Cauchy-estimates and the mean-value theorem, 
in order to estimate on every disc $B^2_r(0)$, 
for $r\in \big{(}0,\frac{1}{8}\big{)}$: 
\begin{eqnarray}   \label{Cauchy}
\parallel \nabla(\lambda_j-\mu_j) \parallel_{L^{\infty}(B^2_r(0))}
= \parallel \nabla(\lambda_j-\mu_j-
\overline{(\lambda_j-\mu_j)}_{B^2_{2r}(0)})
\parallel_{L^{\infty}(B^2_r(0))}      \nonumber                \\
\leq \frac{C}{r^3} \, \int_{B^2_{2r}(0)} 
\big{|}\lambda_j-\mu_j-
\overline{(\lambda_j-\mu_j)}_{B^2_{2r}(0)}\big{|} \, d\Lno, 
\end{eqnarray}  
for every $j\in \nat$,
where $\overline{(\lambda_j-\mu_j)}_{B^2_{2r}(0)}$ denotes the mean value of $\lambda_j-\mu_j$ over $B^2_{2r}(0)$ with respect to the Lebesgue measure $\Lno$. Now, we fix some 
$p\in (1,2)$ and combine \eqref{Cauchy} with H\"older's inequality and with Poincar\'e's inequality, in order to obtain:
\begin{eqnarray}   \label{Cauchy.2} 
\parallel \nabla(\lambda_j-\mu_j) \parallel_{L^{\infty}(B^2_r(0))}    \quad  \nonumber    \\
\leq \frac{C}{r^3} \, \pi^{1-\frac{1}{p}} \,
(2r)^{2-\frac{2}{p}}  \,
\Big{(}\int_{B^2_{2r}(0)} \,
\big{|}\lambda_j-\mu_j-\overline{(\lambda_j-\mu_j)}
_{B^2_{2r}(0)}
\big{|}^p \, d\Lno \Big{)}^{1/p}           \nonumber \\
\leq C_p \, r^{-3+2-\frac{2}{p}+1}\,
\Big{(}\int_{B^2_{2r}(0)}\,
|\nabla(\lambda_j-\mu_j)|^p \, d\Lno \Big{)}^{1/p}   \quad  \\
= C_p \, r^{-\frac{2}{p}}\,
\Big{(}\int_{B^2_{2r}(0)}\,
|\nabla(\lambda_j-\mu_j)|^p \, d\Lno \Big{)}^{1/p}, \nonumber  
\end{eqnarray} 
for every $r\in \big{(}0,\frac{1}{8}\big{)}$, taking 
the exact scaling behaviour of Poincar\'e's inequality 
in line \eqref{Cauchy.2} into account. Using again the 
harmonicity of $\lambda_j-\mu_j$ respectively of its 
gradient $\nabla (\lambda_j-\mu_j)$, we easily obtain from estimate \eqref{Cauchy.2} and from Theorem 2.1 and 
Remark 2.2 in \cite{Giaquinta.1984}, p. 78 - but exchanging 
here the interior $L^2$-estimates in Remark 2.2 by interior $L^p$-estimates: 
\begin{eqnarray}   \label{Cauchy.3} 
	\parallel \nabla(\lambda_j-\mu_j) \parallel_{L^{\infty}(B^2_r(0))}
	\leq \tilde C_p \,r^{-\frac{2}{p}} \,
	r^{\frac{2}{p}} \,\Big{(}\int_{B^2_{\frac{1}{2}}(0)}
	|\nabla(\lambda_j-\mu_j)|^p \, d\Lno \Big{)}^{1/p}  \nonumber \\
	= \tilde C_p \,\Big{(}\int_{B^2_{\frac{1}{2}}(0)}
	|\nabla(\lambda_j-\mu_j)|^p \, d\Lno \Big{)}^{1/p},  
\end{eqnarray} 
for every \,$j\in \nat$, for every 
$r \in \big{(}0,\frac{1}{8}\big{)}$ and 
for the fixed $p\in (1,2)$. Hence, we infer from estimates \eqref{uniformly.conformal.R3.2}, \eqref{estimate.v.j}, \eqref{estimate.Riviere} and \eqref{Cauchy.3}:
\begin{eqnarray}   
\Big{(}\int_{B^2_{r}(0)}
|\nabla(\lambda_j-\mu_j)|^2 \, d\Lno \Big{)}^{1/2}
\leq  \sqrt \pi \, r \,
\parallel \nabla(\lambda_j-\mu_j) \parallel_{L^{\infty}(B^2_r(0))}  \quad  \nonumber  \\ 
\leq \sqrt \pi \,\tilde C_p \,r \,
\Big{(}\int_{B^2_{\frac{1}{2}}(0)}
|\nabla(\lambda_j-\mu_j)|^p \, d\Lno \Big{)}^{1/p}
\leq \textnormal{Const}(\tilde \Lambda,\varepsilon_0,K,\delta,p) \,\,r             \quad             \label{Cauchy.4}  
\end{eqnarray} 
for every $j\in \nat$ and for every 
$r \in \big{(}0,\frac{1}{8}\big{)}$, where we have combined  
estimates \eqref{uniformly.conformal.R3.2}, 
\eqref{estimate.v.j} and \eqref{estimate.Riviere} with the 
continuity of the embedding 
$L^{2,\infty}(B^2_{\frac{1}{2}}(0)) \hookrightarrow  
L^{p}(B^2_{\frac{1}{2}}(0))$, for the fixed $p \in (1,2)$, 
in order to obtain the last inequality in \eqref{Cauchy.4}.
See here Section 3.2 in \cite{Riviere.Park.City.2013} and 
the literature mentioned there.   
Hence, on account of \eqref{Cauchy.4} we can determine 
some small radius $r_0\in \big{(}0,\frac{1}{8}\big{)}$, 
depending only on $\Lambda,\varepsilon_0,K,p$, 
$g_{\textnormal{poin}}$ and on $\delta$, 
such that the integral $\int_{B_{r_0}^2(0)}
|\nabla(\lambda_j-\mu_j)|^2 \, d\Lno$ is smaller 
than $\varepsilon_0^2$. 
Combining this with estimate \eqref{estimate.Riviere} 
we finally infer, that the conformal 
factors $\lambda_j=\tilde u_j\circ \psi_j+v_j$ of 
the conformal embeddings $f_j \circ \psi_j$ satisfy:
\begin{eqnarray}   \label{final.estim} 
\Big{(}\int_{B_{r_0}^2(0)}
|\nabla(\tilde u_j\circ \psi_j+v_j)|^2 \, d\Lno \Big{)}^{1/2} 
\leq (C+1)\,\varepsilon_0,
\,\,\, \textnormal{for every} \,\,\, j\in \nat,  
\end{eqnarray} 
where $C>1$ is the same absolute constant as in  \eqref{estimate.Riviere}, where $r_0$ depends only on $\Lambda,K,\varepsilon_0,p$, $g_{\textnormal{poin}}$ 
and on $\delta$, and where $\varepsilon_0$ had to be  
chosen sufficiently small in \eqref{second.fundam.locally.bounded}. On account of 
convergence \eqref{smooth.converg.metrics.2} 
and estimate \eqref{estimate.v.j}, and on account of the 
conformal invariance of the Dirichlet-integral
estimate \eqref{final.estim} implies immediately:  
\begin{eqnarray}   \label{final.estim.2} 
\Big{(}\int_{B^{g_{\textnormal{poin}}}_{\varrho_0}(x_0)}
|\nabla(\tilde u_j)|^2 \, 
d\mu_{g_{\textnormal{poin}}} \Big{)}^{1/2} 
\leq (C+2)\,\varepsilon_0,
\,\,\, \textnormal{for every} \,\,\, j\in \nat,  
\end{eqnarray}  
for the same absolute constant $C$ as in 
\eqref{estimate.Riviere} and \eqref{final.estim}, 
where $B^{g_{\textnormal{poin}}}_{\varrho}(x_0)$ denotes an open geodesic disc of radius $\varrho$ about the fixed center point $\psi_j(0)=x_0\in \Sigma$ with respect to the metric 
$g_{\textnormal{poin}}$ from \eqref{smooth.converg.metrics.2}, and where $\varrho_0$ is a sufficiently small positive 
number which only depends on $\Lambda,K,\varepsilon_0,p$, 
$g_{\textnormal{poin}}$ and on $\delta$ - just as $r_0$ did - but not on the choice of $x_0$
\footnote{Estimates \eqref{final.estim} and \eqref{final.estim.2} 
are also asserted in Theorem 2.2 in \cite{Palmurella.2022}, 
where the reader is advised to check our basic reference \cite{Riviere.Park.City.2013} as well.}. 
Gathering all estimates in 
\eqref{uniformly.conformal.R3}--\eqref{second.fundam.locally.bounded.3} and in \eqref{final.estim}--\eqref{final.estim.2} we can actually apply the entire reasoning of the proof of Theorem 4.1 in \cite{Palmurella.2022}, in particular estimates (4.10) and (4.13)--(4.16) in Proposition 4.7 and Lemma 4.9 of \cite{Palmurella.2022}, in order to obtain here the estimate 
\begin{eqnarray}  \label{L.infty.W.4.2.estimate.local}  
\parallel \nabla^{g_{\textnormal{poin}}}(f_j) 
\parallel_{W^{3,2}(B^{g_{\textnormal{poin}}}_{\varrho_0}(x_0))}^2 \leq \textnormal{Const} \cdot
\Big{(} \int_{B^{g_{\textnormal{poin}}}_{R_0}(x_0)} 
|\nabla_{L^2} \Will(f_j)|^2 \, 
d\mu_{f_j^*(g_{\textnormal{euc}})} \,+ 1 \,\Big{)}  \quad             
\end{eqnarray}
for every $j\in \nat$, provided $\varepsilon_0$ had been 
chosen sufficiently small in \eqref{second.fundam.locally.bounded}, where the small radius $\varrho_0$ had been determined in \eqref{final.estim.2}
and where the constant in \eqref{L.infty.W.4.2.estimate.local} depends on $g_{\textnormal{poin}},\Will(F_0),K,\Lambda,\delta$ 
and on the choice of $\varepsilon_0$, but not 
on the choice of $x_0$. 
Now, the final step of the proof works as in the end of the proof of Proposition 5.2 in \cite{Palmurella.2022}.
Since the center $x_0\in \Sigma$ of the open
geodesic disc $B^{g_{\textnormal{poin}}}_{\varrho_0}(x_0)$ 
in $(\Sigma,g_{\textnormal{poin}})$ had been chosen arbitrarily in $\Sigma$ and since $\Sigma$ is compact, we infer from estimate \eqref{L.infty.W.4.2.estimate.local}
- covering $\Sigma$ with finitely many appropriate coordinate patches $\psi^i_j: B^2_{r_0}(0) 
\stackrel{\cong}\longrightarrow \psi^i_j(B^2_{r_0}(0)) 
\subset \Sigma$ as in \eqref{estimate.v.j}, for $i=1,\ldots,N=N(r_0,g_{\textnormal{poin}},\Sigma)$ 
with $N$ being independent of $j\in \nat$ on account of  convergence \eqref{smooth.converg.metrics.2} - that estimate \eqref{L.infty.W.4.2.estimate} actually holds globally 
on $\Sigma$ for the sequence of embeddings 
$f_j = \Stereo \circ \tilde F_{t_{j}}:\Sigma \longrightarrow \rel^3$, with a large constant depending only on  $g_{\textnormal{poin}},\Will(F_0),K,\Lambda,\Sigma$ and on $\delta$, similarly to the constant in \eqref{L.infty.W.4.2.estimate.local}.
In order to obtain estimate \eqref{L.infty.W.4.2.estimate} 
for the original embeddings 
$\tilde F_{t_{j}}:\Sigma \longrightarrow \sphere^3$   
we only have to apply now the inverse stereographic projection 
$\Stereo^{-1}:\rel^3 \longrightarrow \sphere^3\setminus \{(0,0,0,1)\}$ to $f_j$, explicitly given by: 
$$ 
(x,y,z)\mapsto \frac{1}{x^2+y^2+z^2+1} \, \big{(}2x,2y,2z,x^2+y^2+z^2-1\big{)}. 
$$
It easily follows from Section 3 of \cite{Cairns.Sharpe.Webb.1994} 
combined with formula (5) in \cite{Jakob_Moebius_2016} 
that the non-linear map 
\begin{equation} \label{conformal.invariant.oper} 
\textnormal{Imm}_{\textnormal{uf}}(\Sigma,\rel^n) 
\ni f \mapsto \int_{\Sigma} \frac{1}{|A^0_f|^4} \, 
|\nabla_{L^2} \Will(f)|^2 \, d\mu_{f^*(g_{\textnormal{euc}})} \in \rel
\end{equation} 
is a conformally invariant operator, for any fixed $n\geq 3$. 
Hence, estimate \eqref{L.infty.W.4.2.estimate.local}
immediately implies estimate \eqref{L.infty.W.4.2.estimate}, 
up to checking the facts that $|A^0_{\tilde F_{t_{j}}}|^2 
\equiv |A^0_{\tilde F_{t_{j}},\sphere^3}|^2$ 
and also $|A^0_{f_{j}}|^2$ remain uniformly bounded on $\Sigma$, for all $j \in \nat$.  
Indeed, the traces $|A^0_{\tilde F_{t_{j}},\sphere^3}|^2$ 
remain bounded from above in terms of the uniform upper bound in \eqref{second.fundam.bounded} 
and they also remain bounded from below by the number $2$ 
on account of formula \eqref{second.fundam.trace.free} - recalling here that the considered flow line $\{F_t\}$ 
consists of smooth parametrizations of Hopf-tori in $\sphere^3$ on account of Proposition \ref{correspond.flows}. Obviously, these two estimates yield an upper and a lower bound for the trace $|A^0_{f_{j}}|^2$ on $\Sigma$ 
for each $f_j = \Stereo \circ \tilde F_{t_{j}}:
\Sigma \longrightarrow \rel^3$ in terms 
of $K$ and $\delta$, similarly to the elementary argument 
leading to estimate \eqref{second.fundam.bounded.2}. 
\item[3)] Now we additionally assume that the speed $|\frac{d}{dt}\Will(F_{t})|$
of ``energy decrease'' remains uniformly bounded 
by some large constant $W$ at every time $t=t_{j}$, and 
we consider a weakly/weakly* convergent subsequence 
$\{\tilde F_{t_{j_k}}\}$ in $W^{2,2}((\Sigma,g_{\textnormal{poin}}),\rel^4)$ and in $W^{1,\infty}((\Sigma,g_{\textnormal{poin}}),\rel^4)$ 
as in \eqref{weak.convergence.2} and 
\eqref{W.1.infty.convergence},
which we had obtained from the original sequence 
$\{F_{t_{j_l}}\}$ in part (2) of Theorem \ref{limit.MIWF}
by extraction of an appropriate subsequence 
$\{F_{t_{j_k}}\}$ and by appropriate reparametrization 
of each embedding $F_{t_{j_k}}$. 
Hence, we have here additionally: 
\begin{equation} \label{bounded.energy.decrease}     
\int_{\Sigma} \frac{1}{|A^0_{F_{t_{j_k}}}|^4}\, 
|\nabla_{L^2} \Will(F_{t_{j_k}})|^2 \, 
d\mu_{F_{t_{j_k}}^*(g_{\textnormal{euc}})} 
=2\,\Big{|}\frac{d}{dt}\Will(F_t)\Big{|}
\lfloor_{t=t_{j_k}}\leq 2W, \,\,\, \forall \, k \in \nat.   
\end{equation}
On account of the uniform upper bound \eqref{second.fundam.bounded.0}, 
we obtain from \eqref{bounded.energy.decrease} that 
$$  
\int_{\Sigma} |\nabla_{L^2} \Will(F_{t_{j_k}})|^2 \, 
d\mu_{F_{t_{j_k}}^*(g_{\textnormal{euc}})} 
\leq (K^2+2)^2\,2W, \quad \forall \, k \in \nat,   
$$
and thus also 
\begin{equation} \label{bounded.energy.decrease.2}     
\int_{\Sigma} |\nabla_{L^2} \Will(\tilde F_{t_{j_k}})|^2 \, d\mu_{\tilde F_{t_{j_k}}^*(g_{\textnormal{euc}})} 
\leq (K^2+2)^2\,2W, \quad \forall \, k \in \nat, 
\end{equation}
on account of the invariance of the differential  
operator $F\mapsto \nabla_{L^2} \Will(F)$ with respect to smooth reparametrization and on account of the definition
of $\tilde F_{t_{j}}$ below formula \eqref{smooth.converg.metrics}.  
Combining now estimates \eqref{L.infty.W.4.2.estimate} 
and \eqref{bounded.energy.decrease.2} both estimating  
the sequence $\{\tilde F_{t_{j_k}}\}$
and recalling that the embeddings $\tilde F_{t_{j_k}}$ 
converge in the senses \eqref{weak.convergence.2}, \eqref{W.1.infty.convergence} and \eqref{C0.convergence}, 
then the ``principle of subsequences'' yields that 
$\{\tilde F_{t_{j_k}}\}$ also converges weakly in $W^{4,2}((\Sigma,g_{\textnormal{poin}}),\rel^4)$ 
to the conformal bi-Lipschitz parametrization
$f:(\Sigma,g_{\textnormal{poin}}) \stackrel{\cong}\longrightarrow 
\textnormal{spt}(\mu)$ of the corresponding limit Hopf-torus 
$\textnormal{spt}(\mu)$ from the first part of this
theorem, i.e. to the parametrization of the support of the limit varifold $\mu$ of the weakly convergent subsequence $\{F_{t_{j_l}}\}$ from line \eqref{weak.convergence.F.j}.
From Rellich's embedding theorem, A 8.4 in \cite{Alt.2015},
we immediately infer also strong convergence of the  
sequence $\{\tilde F_{t_{j_k}}\}$ to $f$ in  $W^{3,2}((\Sigma,g_{\textnormal{poin}}),\rel^4)$, 
as $k\to \infty$. Since this implies together with formula 
\eqref{Willmore.energy} and with convergence \eqref{W.1.infty.convergence} in particular \, 
$\Will(\mu)= \Will(f) = 
\lim_{k\to \infty} \Will(\tilde F_{t_{j_k}})
= \lim_{k\to \infty} \Will(F_{t_{j_k}})$,   
thus all conditions of the third part of Theorem \ref{limit.MIWF} are satisfied by the sequence $\{F_{t_{j_k}}\}$, and therefore statement \eqref{no.A.concentration} has to hold here, 
just as asserted, for the reparametrized 
sequence $\{\tilde F_{t_{j_k}}\}$ or equivalently 
for the original sequence $\{F_{t_{j_k}}\}$ itself, 
taking the invariance  
of the functional in \eqref{no.A.concentration} 
with respect to smooth reparametrization into account.  
Finally, we infer from estimates \eqref{L.infty.W.4.2.estimate} and \eqref{bounded.energy.decrease.2} together 
with the compactness of the embedding  
$W^{4,2}((\Sigma,g_{\textnormal{poin}}),\rel^4)
\hookrightarrow C^{2,\alpha}((\Sigma,g_{\textnormal{poin}}),\rel^4)$, 
for any $\alpha \in (0,1)$, 
that $\{\tilde F_{t_{j_k}}\}$ converges to $f$ in $C^{2,\alpha}((\Sigma,g_{\textnormal{poin}}),\rel^4)$, 
as $k\to \infty$. 
\end{itemize}
\noindent 
For the proof of Theorem \ref{limit.at.infinity} we firstly 
recall here Theorem 1 of \cite{Ruben.MIWF.IV}.
\begin{proposition}  \label{Center.manifold} 
	Let $\Sigma$ be a smooth compact torus, and let 
	$F^*:\Sigma \stackrel{\cong}\longrightarrow 
	M\big{(}\frac{1}{\sqrt 2}(\sphere^1 \times \sphere^1)\big{)}$
	be a smooth diffeomorphic parametrization of a compact torus in $\sphere^3$, which is conformally equivalent to the standard Clifford torus $\frac{1}{\sqrt 2}(\sphere^1 \times \sphere^1)$ via some conformal transformation 
	$M\in \textnormal{M\"ob}(\sphere^3)$, 
	and let some $\beta \in (0,1)$ and $m \in \nat$ be fixed. Then, there is some small neighborhood $W=W(\Sigma,F^*,m)$ about $F^*$ in $h^{2+\beta}(\Sigma,\rel^4)$, such that for every $C^{\infty}$-smooth initial immersion 
	$F_1:\Sigma \longrightarrow \sphere^3$, which is contained in $W$, the unique flow line $\{\PP(t,0,F_1)\}_{t \geq 0}$ of the MIWF exists globally and converges - up to smooth reparametrization - fully to a smooth and diffeomorphic parametrization of a torus in $\sphere^3$, which is again conformally equivalent to the standard Clifford torus in $\sphere^3$. This full convergence takes place with respect to the $C^m(\Sigma,\rel^4)$-norm and at an exponential rate, 
	as $t\nearrow \infty$. 
\end{proposition}
\noindent  
\underline{{\bf Proof of Theorem \ref{limit.at.infinity}}}:\\\\
\noindent 
We consider in Theorem \ref{limit.at.infinity} 
some arbitrarily chosen flow line $\{F_t\}$ of the MIWF which starts moving with Willmore energy strictly below $4\pi^2$. Hence recalling the monotonicity of the energy 
$t\mapsto \Will(F_t)$ from formula \eqref{Willmore.monotonicity} we can infer from the second part of Remark \ref{Simple.maps} that each immersion $F_t$ is a \emph{simple parametrization} of its image in $\sphere^3$, in the sense of Definition \ref{Simple.map}.   
Next we employ the strongest assumption of Theorem \ref{limit.at.infinity}, namely that the considered flow line $\{F_t\}$ is global, i.e. that it satisfies $T_{\textnormal{max}}(F_0)=\infty$.  
Again because of the monotonicity
of the energy $t\mapsto \Will(F_t)$ along the considered 
flow line $\{F_t\}$, we know that 
$v:=\lim_{t\to \infty} \Will(F_t)$ has to exist. 
Now there are two possibilities for the limit $v$: 
(i) $v \geq 8\pi$ or (ii) $v < 8\pi$. 
Let's start discussing the first case. 
Applying here the same argument which  
we have already used at the beginning of the proof of Theorem \ref{limit.MIWF} between formulae \eqref{Willmore.monotonicity} and \eqref{Energy.smaller.8pi} 
- combining inequality \eqref{Willmore.monotonicity} with the real analyticity of the function 
$t \mapsto \Will(F_t)$ due to Theorem 3 in \cite{Ruben.MIWF.III} - we either have a 
\emph{stationary flow line}, i.e. precisely (a) 
$8\pi\leq v \leq \Will(F_t)=\Will(F_0)< 
4\pi^2$ for all $t\geq 0$, 
or (b) a strictly monotonically decreasing flow line of 
the MIWF, i.e. satisfying: 
\begin{equation} \label{strict.monotonicity} 
8\pi \leq v< \Will(F_{t_2})<\Will(F_{t_1}) \quad 
\textnormal{for every pair} \quad t_2 >t_1\geq 0. 
\end{equation} 
Now, because of Propositions \ref{Hopf.Willmore.prop} and \ref{correspond.flows} the first case (a) would immediately 
imply the existence of a stationary flow line 
$\{\gamma_t\}$ of flow \eqref{elastic.energy.flow} 
with $8\leq \Wil(\gamma_t)=\Wil(\gamma_0)< 4\pi$ 
for all $t\geq 0$. This would contradict the second part of Proposition 6 in \cite{Ruben.MIWF.II}, 
stating that there is no elastic curve on 
$\sphere^2$ with elastic energy in $(2\pi, 4\pi)$. 
Hence, we must have here statement \eqref{strict.monotonicity}. Now we again use our main
assumption on the flow line $\{F_t\}$ to be global, 
namely we integrate inequality \eqref{Willmore.monotonicity} 
from $0$ to $\infty$, and we thus conclude: 
\begin{eqnarray*}  
\int_0^{\infty} \int_{\Sigma} 
\frac{1}{|A^0_{F_{t}}|^4}\, 
|\nabla_{L^2} \Will(F_{t})|^2 \,  
d\mu_{F_{t}^*(g_{\textnormal{euc}})} \, dt    
=- 2 \, \lim_{T\to \infty} \int_0^{T}
\frac{d}{dt} \Will(F_{t}) \, dt        =    \\
=2\, \lim_{T\to \infty}(\Will(F_{0})-\Will(F_{T})) 
< 2\,\Will(F_{0}).  
\end{eqnarray*} 
Hence, there has to be some sequence $t_j\nearrow \infty$
such that 
\begin{equation}  \label{speed.to.zero}  
\int_{\Sigma} \frac{1}{|A^0_{F_{t_j}}|^4}\, 
|\nabla_{L^2} \Will(F_{t_j})|^2 \, d\mu_{F_{t_j}^*(g_{\textnormal{euc}})} 
= 2\, \Big{|}\frac{d}{dt}\Will(F_{t})\Big{|}\lfloor_{t=t_j}
\longrightarrow 0,
\end{equation} 
as $j\to \infty$. 
Now, if some smooth reparametrization $\{\tilde F_t\}$ 
of the flow line $\{F_t\}$ would fully converge in $C^4(\Sigma,\rel^4)$ to some $C^4$-immersion $F^*$, then 
first of all $\sup_{\Sigma} |A^0_{F_{t_j}}|^2$ would remain uniformly bounded on $\Sigma$ for every $j$, and secondly 
convergence \eqref{speed.to.zero} would imply: 
$$
0 \longleftarrow \int_{\Sigma} |\nabla_{L^2} \Will(F_{t_j})|^2 \, d\mu_{g_{F_{t_j}}}
= \int_{\Sigma} |\nabla_{L^2} \Will(\tilde F_{t_j})|^2 \, d\mu_{g_{\tilde F_{t_j}}} 
\longrightarrow \int_{\Sigma} |\nabla_{L^2} \Will(F^*)|^2 \, d\mu_{g_{F^*}},
$$ 
as $j\to \infty$, showing that $\nabla_{L^2} \Will(F^*)=0$ 
on $\Sigma$, i.e. that $F^*$ would actually have to be a smooth Willmore immersion with Willmore energy $\Will(F^*)=v$, 
simply because    
\begin{equation} \label{Will.in.the.end}
\Will(F^*) = \lim_{j \to \infty} \Will(\tilde F_{t_j}) 
= \lim_{j \to \infty} \Will(F_{t_j}) 
=v \in \Big{[}8\pi, 4\pi^2\Big{)}.
\end{equation}
On the other hand, recalling Definition \ref{Hopf.torus.immersion} and Proposition \ref{closed.lifts}
we can easily construct profile curves 
$\gamma_{t_j}:\sphere^1 \longrightarrow \sphere^2$ 
of the smooth Hopf-torus-immersions $F_{t_j}$ 
respectively $\tilde F_{t_j}$ which  
converge in $C^4(\sphere^1,\rel^3)$ 
to some regular path $\gamma^*\in C^4(\sphere^1,\rel^3)$. 
Combining now convergence \eqref{speed.to.zero} with 
\eqref{First.Var.Will.Wil} we would obtain for the sequence
$\{\gamma_{t_j}\}$ and its limit $\gamma^*$: 
$$ 
\int_{\sphere} \frac{1}{(\kappa_{\gamma^*}^2+1)^2} \, 
|\nabla_{L^2} \Wil(\gamma^*)|^2 \, d\mu_{\gamma^*} 
\longleftarrow 
\int_{\sphere} \frac{1}{(\kappa_{\gamma_{t_j}}^2+1)^2} \, 
|\nabla_{L^2} \Wil(\gamma_{t_j})|^2 \, d\mu_{\gamma_{t_j}} 
\longrightarrow 0  \quad \textnormal{as} \, j\to \infty,   
$$ 
implying that the regular path $\gamma^*$ would have to be 
a smooth elastic curve, necessarily having 
elastic energy 
\begin{equation} \label{Wil.in.the.end}
\Wil(\gamma^*) = \lim_{j \to \infty} \Wil(\gamma_{t_j}) 
=\lim_{j \to \infty} \frac{1}{\pi} \Will(F_{t_j}) 
= \frac{v}{\pi} \in  \Big{[}8,4\pi\Big{)}, 
\end{equation} 
where we applied in the latter two equations 
firstly identity \eqref{Will.Wil} - again using
the fact that each $\gamma_{t_j}$ is the profile 
curve of the Hopf-torus-immersion $F_{t_j}$ - and 
then the third equation in \eqref{Will.in.the.end}.  
But the conclusion in \eqref{Wil.in.the.end} again 
contradicts the second part of Proposition 6 
in \cite{Ruben.MIWF.II}, similarly to our argument above, 
dealing with subcase (a) of case (i).
Hence, indeed in the first case ``$v\geq 8\pi$'' 
no smooth reparametrization of the considered global 
flow line $\{F_t\}$ can fully converge in 
$C^4(\Sigma,\rel^4)$ to some  
$C^4$-immersion $F^*$, as $t \to \infty$, thus 
showing the asserted ``divergent behaviour'' of 
$\{F_t\}$ as $t \to \infty$.  \\
Now, in the second alternative 
``$\lim_{t\to \infty} \Will(F_t) = v< 8 \pi$'' we 
could prove the strict monotonicity of 
$t\mapsto \Will(F_t)$ along the considered flow line 
$\{F_t\}$ just as we did above in the first case 
``$v\geq 8 \pi$'', and therefore we can pick here some sufficiently large time $T_0>0$ such that 
$\Will(F_t) < \Will(F_{T_0})<8 \pi$ holds for every 
$t>T_0$. Therefore, in the case ``$v< 8 \pi$'' 
we can assume without loss of generality 
that the considered flow line $\{F_t\}$ of the MIWF 
would start moving in a smooth 
parametrization $F_0:\Sigma \longrightarrow \sphere^3$ 
of a Hopf-torus in $\sphere^3$ with $\Will(F_0)< 8\pi$,
implying that every Hopf-torus-immersion $F_t$ 
which belongs to the considered flow line of the MIWF 
has to be an embedding for all times $t\geq 0$.  \\ 
Now we recall the second a-priori assumption of this 
theorem, namely that there is some constant $K>0$ 
such that $\parallel \vec H_{F_{t},\sphere^3} \parallel_{L^{\infty}(\Sigma)}\leq K$ for every  
$t\in [0,\infty)$, implying here automatically estimates \eqref{second.fundam.bounded} and \eqref{second.fundam.bounded.2} on account of 
formula \eqref{second.fundam.} in Proposition 
\ref{Hopf.Willmore.prop}.
Now, as in our discussion of the alternative  
``$v \geq 8 \pi$'' we can deduce convergence  \eqref{speed.to.zero} from the assumption on $\{F_t\}$ 
to be global. However, on account of statement \eqref{speed.to.zero} the sequence $\{F_{t_j}\}$ actually
satisfies all requirements of the third part of Theorem \ref{singular.time.MIWF.Hopf.tori}. 
Hence, we can insert exactly the divergent sequence of times 
$t_j \nearrow T_{\textnormal{max}}(F_0)=\infty$ 
satisfying \eqref{speed.to.zero} into the third part of Theorem \ref{singular.time.MIWF.Hopf.tori}, and any reparametrized subsequence $\{\tilde F_{t_{j_k}}\}$ converging weakly/weakly* 
as in \eqref{weak.convergence.2} and \eqref{W.1.infty.convergence} - which we had considered 
in the second and third part of Theorem \ref{singular.time.MIWF.Hopf.tori} - 
converges even weakly in $W^{4,2}((\Sigma,g_{\textnormal{poin}}),\rel^4)$
and strongly in $W^{3,2}((\Sigma,g_{\textnormal{poin}}),\rel^4)$
to a uniformly conformal bi-Lipschitz homeomorphism 
$f$ between $(\Sigma,g_{\textnormal{poin}})$ and a Hopf-torus  
$\textnormal{spt}(\mu)$ in $\sphere^3$, where 
$g_{\textnormal{poin}}$ is some appropriate smooth 
metric on $\Sigma$ of vanishing scalar curvature. Now we argue as in \eqref{containoe} and apply the stereographic 
projection $\Stereo:\sphere^3 \setminus \{(0,0,0,1)\} 
\stackrel{\cong} \longrightarrow \rel^3$.  
Hence, we can conclude that also the smooth embeddings  
$\tilde f_k:= \Stereo \circ \tilde F_{t_{j_k}}$ converge 
weakly in $W^{4,2}((\Sigma,g_{\textnormal{poin}}),\rel^3)$
and strongly in $W^{3,2}((\Sigma,g_{\textnormal{poin}}),\rel^3)$ 
to the uniformly conformal bi-Lipschitz homeomorphism 
$\tilde f:=\Stereo \circ f$ between $(\Sigma,g_{\textnormal{poin}})$ and  
$\Stereo(\textnormal{spt}(\mu)) \subset \rel^3$. 
Moreover, we can use here again the conformal invariance of  
the operator 
$f\mapsto \int_{\Sigma} \frac{1}{|A^0_f|^4} \, 
|\nabla_{L^2} \Will(f)|^2 \, d\mu_{f^*(g_{\textnormal{euc}})}
$
from formula \eqref{conformal.invariant.oper}, and therefore \eqref{speed.to.zero} also implies: 
$$
\int_{\Sigma} \frac{1}{|A^0_{\tilde f_k}|^4}\, 
|\nabla_{L^2} \Will(\tilde f_{k})|^2 \, 
d\mu_{\tilde f_{k}^*(g_{\textnormal{euc}})} 
= \int_{\Sigma} \frac{1}{|A^0_{\tilde F_{t_{j_k}}}|^4}\, 
|\nabla_{L^2} \Will(\tilde F_{t_{j_k}})|^2 \, 
d\mu_{\tilde F_{t_{j_k}}^*(g_{\textnormal{euc}})} \longrightarrow 0 
$$ 
as $k\to \infty$. Combining this convergence with estimate  
\eqref{second.fundam.bounded.2} we obtain: 
\begin{equation}  \label{first.variation.to.zero}  
\int_{\Sigma} |\nabla_{L^2} \Will(\tilde f_{k})|^2 \, 
d\mu_{\tilde f_{k}^*(g_{\textnormal{euc}})}   
\longrightarrow 0, 	\quad \textnormal{as} \,\,\,k\to \infty.
\end{equation} 
Now, similarly to the argument in line \eqref{distributional.Willmore.limit}, we can 
compute here by means of formulae \eqref{distributional.Willmore}, \eqref{Riviere.Bernard.Willmore} and \eqref{first.variation.to.zero} and by means of the 
above mentioned strong convergence of $\{\tilde f_k\}$ in $W^{3,2}((\Sigma,g_{\textnormal{poin}}),\rel^3)$ to 
$\tilde f$:   
\begin{eqnarray}  \label{distributional.Willmore.limit.2} 
0\longleftarrow \langle \nabla_{L^2}\Will(\tilde f_k),\varphi 
\rangle_{L^2(\Sigma,\tilde f_k^*(g_{\textnormal{euc}}))} 
= \langle \nabla_{L^2}\Will(\tilde f_k),\varphi \rangle_{\Dom'(\Sigma)} 	
\longrightarrow 
\langle \nabla_{L^2}\Will(\tilde f),\varphi \rangle_{\Dom'(\Sigma)},	 
\end{eqnarray}  
as $k \to \infty$, for every fixed 
$\varphi \in C^{\infty}(\Sigma,\rel^3)$. Hence, we infer 
from \eqref{distributional.Willmore.limit.2} that 
$\nabla_{L^2}\Will(\tilde f)\equiv 0$ 
in the distributional sense of \eqref{distributional.Willmore}, i.e. that the uniformly conformal bi-Lipschitz homeomorphism $\tilde f$ is 
\emph{weakly Willmore} on $\Sigma$ in the sense of 
Corollary 7.3 in \cite{Riviere.Park.City.2013} respectively 
of Definition VII.3 in \cite{Riviere.2011}. 
We can therefore immediately infer from Theorem 7.11 in 
\cite{Riviere.Park.City.2013} respectively from 
Corollary VII.6 in \cite{Riviere.2011} that 
$\tilde f$ is actually a smooth diffeomorphism between 
$(\Sigma,g_{\textnormal{poin}})$ and an embedded classical Willmore surface in $\rel^3$, where we have strongly relied 
on the fact that $\tilde f$ has already been known to be 
a uniformly conformal bi-Lipschitz homeomorphism onto $\Stereo(\textnormal{spt}(\mu))$. 
Hence applying now inverse stereographic projection from $\rel^3$ to $\sphere^3\setminus \{(0,0,0,1)\}$ and recalling the statement of the third part of Theorem \ref{singular.time.MIWF.Hopf.tori}, 
the original limit embedding  
$f:\Sigma \stackrel{\cong} \longrightarrow \textnormal{spt}(\mu)$ turns out to parametrize a smooth Willmore-Hopf-torus in $\sphere^3$. Since we assume here 
that $\Will(F_t)<8\pi$ for every $t>0$, we must have 
- combining statement \eqref{Willmore.monotonicity} 
with the proven Willmore conjecture, 
Theorem A in \cite{Marques.Neves.2014} - 
in the limit as $t_{j_k}\nearrow \infty$: 
$\Will(f) \in [2\pi^2,8\pi)$, and thus by formula 
\eqref{Will.Wil}: $\Wil(\gamma) \in [2\pi,8)$ for the elastic 
energy of any smooth profile curve $\gamma$ of the Willmore-Hopf-torus $\textnormal{spt}(\mu)$. 
Moreover, we infer here from formula \eqref{First.Var.Will.Wil} that any smooth profile curve $\gamma$ of the Willmore-Hopf-torus $\textnormal{spt}(\mu)$ is an elastic curve in $\sphere^2$, i.e. solves:
$2\,\Big{(}\nabla^{\perp}_{\frac{\gamma'}{|\gamma'|}}
\Big{)}^2(\vec{\kappa}_{\gamma})
+ |\vec{\kappa}_{\gamma}|^2 \vec{\kappa}_{\gamma}
+ \vec{\kappa}_{\gamma} \equiv 0$ on $\sphere^1$. 
Hence, again applying the second part of Proposition 6 
in \cite{Ruben.MIWF.II} we can infer here from 
$\Wil(\gamma) \in [2\pi,8)$ that actually 
$\Wil(\gamma) = 2\pi$ has to hold, 
implying that any chosen profile curve $\gamma$ 
of the Willmore-Hopf-torus $\textnormal{spt}(\mu)$
must parametrize some great circle in $\sphere^2$. Hence, 
$\textnormal{spt}(\mu)=\pi^{-1}(\textnormal{trace}(\gamma))$
has to be the Clifford torus in $\sphere^3$ - at least up to 
some appropriate conformal transformation of $\sphere^3$.
We can therefore conclude that  
$f:(\Sigma,g_{\textnormal{poin}}) 
\stackrel{\cong} \longrightarrow 
M\big{(}\frac{1}{\sqrt 2}(\sphere^1 \times \sphere^1)\big{)}$
is smooth and diffeomorphic, 
where $M$ is an appropriate M\"obius-transformation 
of $\sphere^3$. Now we can again conclude from 
Theorem \ref{singular.time.MIWF.Hopf.tori}
that the reparametrized embeddings   
$\tilde F_{t_{j_k}}$ converge also in $C^{2,\alpha}((\Sigma,g_{\textnormal{poin}}),\rel^4)$ 
to the limit embedding $f$, as $k\to \infty$, for any 
fixed $\alpha \in (0,1)$.  
Hence, we can apply here the above Proposition \ref{Center.manifold} to $F^*:=f$ and $F_1:=\tilde F_{t_{j_{k^*}}}$ for some sufficiently large index $k^*=k^*(m)\in \nat$, with $\beta\in (0,\alpha)$ 
arbitrarily chosen, such that 
$\tilde F_{t_{j_{k^*}}} \in W(\Sigma,f,m)$ 
in the terminology of Proposition \ref{Center.manifold}, 
where we have used the fact that 
$C^{2,\alpha}((\Sigma,g_{\textnormal{poin}}),\rel^4)
\hookrightarrow h^{2+\beta}((\Sigma,g_{\textnormal{poin}}),\rel^4)$ 
is a continuous embedding for $0<\beta<\alpha$. 
Hence, we obtain from Proposition \ref{Center.manifold} 
the existence of some smooth family of smooth diffeomorphisms 
$\theta^{j_{k^*}}_s:\Sigma \stackrel{\cong}\longrightarrow \Sigma$, for $s\geq t_{j_{k^*}}$, depending on  
$F^*=f$ and $F_1:=\tilde F_{t_{j_{k^*}}}$, such that the reparametrized flow line 
$\{\PP(s,t_{j_{k^*}},\tilde F_{t_{j_{k^*}}}) 
\circ \theta^{j_{k^*}}_s\}_{s\geq t_{j_{k^*}}}$ of the MIWF 
- being only a smooth solution to the \emph{relaxed MIWF-equation} \eqref{generalized.MIWF} below - satisfies: 
\begin{equation}   \label{full.convergence.m}  
\PP(s,t_{j_{k^*}},\tilde F_{t_{j_{k^*}}}) 
\circ \theta^{j_{k^*}}_s \longrightarrow G  \quad \textnormal{in} \,\,\,C^m(\Sigma,\rel^4) 
\end{equation}	
fully as $s\to \infty$ and at an exponential rate,
where $G:\Sigma \stackrel{\cong}\longrightarrow
M^*\big{(}\frac{1}{\sqrt 2}(\sphere^1 \times \sphere^1)\big{)}$ 
is a smooth diffeomorphism and where $M^*$ is another appropriate M\"obius transformation of $\sphere^3$.  
Now, without loss of generality we may assume here 
as in the proof of Proposition \ref{Center.manifold}, 
i.e. of Theorem 1 in \cite{Ruben.MIWF.IV}, that the reference-immersion $F^*=f$ parametrizes the standard 
Clifford torus $\CC$ itself. 
Moreover, a closer inspection of the preparation 
for the proof of Proposition \ref{Center.manifold}, 
especially the section between 
formulae (7) and (25) in \cite{Ruben.MIWF.IV}, 
\footnote{This technique is actually an adaption  
of the standard method of \emph{normal graph representations} of immersions over some fixed 
smooth immersion into $\rel^n$. See here p. 31 in 
\cite{Ruben.MIWF.IV}, Section 5 of 
\cite{Skorzinski.2015} and Section 4 of \cite{Pruess.Simonett.2013} for more precise information.} 
and of Part (iv) of the proof of Proposition \ref{Center.manifold} itself in \cite{Ruben.MIWF.IV}  
one can quickly figure out that the reparametrized flow line 
$\{\PP(s,t_{j_{k^*}},\tilde F_{t_{j_{k^*}}}) 
\circ \theta^{j_{k^*}}_s\}_{s\geq t_{j_{k^*}}}$ actually satisfies \eqref{full.convergence.m}, because it can be written as a normal graph over the embedding 
$f:\Sigma \stackrel{\cong}\longrightarrow \frac{1}{\sqrt 2}(\sphere^1 \times \sphere^1)$ of the Clifford torus $\CC$ in $\sphere^3$ for each $s\geq t_{j_{k^*}}$ by means of the exponential map $\exp^{\sphere^3}$ restricted to the normal bundle $N\CC$ along $\CC$ within $T\sphere^3$, i.e. because  we have here:   
\begin{equation} \label{graph.representation.1}
\PP(s,t_{j_{k^*}},\tilde F_{t_{j_{k^*}}}) 
\circ \theta^{j_{k^*}}_s(x)  
= \exp^{\sphere^3}_{f(x)}\big{(}\rho_s(f(x))\, \nu_{\CC}(f(x))\big{)}  
\quad \forall \,x \in \Sigma, 
\end{equation} 
for every $s\geq t_{j_{k^*}}$, where $\{\rho_s\}$ is some smooth family of smooth real-valued functions on $\CC$, 
and $\nu_{\CC}$ is a fixed section of Euclidean length $1$ 
into the normal bundle $N\CC$ of $\CC$ within $T\sphere^3$. 
Now, recalling here again that  
$\tilde F_{t_{j_{k^*}}}=F_{t_{j_{k^*}}}\circ \Phi_{j_{k^*}}$
we can use the invariance of the MIWF with respect to 
time-independent smooth reparametrizations and formulate statement \eqref{graph.representation.1}  
more accurately as:
\begin{equation} \label{graph.representation.2}
\PP(s,t_{j_{k^*}},F_{t_{j_{k^*}}}) \circ 
(\Phi_{j_{k^*}} \circ \theta^{j_{k^*}}_s)(x) = \exp^{\sphere^3}_{f(x)}\big{(}\rho_s(f(x))\, 
\nu_{\CC}(f(x))\big{)} \quad \forall \,x \in \Sigma, 
\end{equation} 
for every $s\geq t_{j_{k^*}}$. Now, statement \eqref{graph.representation.2} shows on account of 
$F_{t_{j_{k^*}}}=\PP(t_{j_{k^*}},0,F_{0})$ that 
$\{\rho_s(f(\,\cdot\,))\, \nu_{\CC}(f(\,\cdot\,))\}
_{s\geq t_{j_{k^*}}}$ is the 
\underline{unique family of smooth sections} 
into the pullback bundle $f^*N\CC$ of $N\CC$ 
which represents the considered flow line $\{F_s\}=\{\PP(s,0,F_{0})\}$ 
as a \emph{normal graph along} $f$, at least for every 
$s\geq t_{j_{k^*}}$, and furthermore that 
$\{\Phi_{j_{k^*}} \circ \theta^{j_{k^*}}_s\}$ is the 
\underline{unique family of smooth reparametrizations} of 
the considered flow line 
$\{\PP(s,0,F_{0})\}_{s\geq t_{j_{k^*}}}$ 
such that the reparametrized surface 
$\PP(s,0,F_{0})\circ (\Phi_{j_{k^*}} \circ \theta^{j_{k^*}}_s)$ can indeed be written as a 
\emph{normal graph along} $f$ - 
as in \eqref{graph.representation.2} by means of the exponential map - at least for each $s\geq t_{j_{k^*}}$.
\footnote{See here especially Theorem 5.1 in \cite{Skorzinski.2015} for detailed constructions and 
explanations, at least concerning the standard situation 
of immersions into some $\rel^n$.}
On account of this uniqueness we can conclude 
that the value of the initial time $s=t_{j_{k^*}}$ in \eqref{graph.representation.1} and \eqref{graph.representation.2}, 
i.e. the size of the index $k^*=k^*(m)$ respectively the choice of $m$ in the formulation of the asserted theorem, 
did neither affect the choice of smooth 
sections into the pullback bundle $f^*N\CC$ on the 
right hand side of equation \eqref{graph.representation.2} nor the choice of smooth reparametrizations on the left hand side 
of equation \eqref{graph.representation.2}, whenever they are 
defined, i.e. for any $s\geq t_{j_{k^*}}$. 
Therefore formulae \eqref{full.convergence.m}-\eqref{graph.representation.2}  
give rise to a well-defined and smooth family 
of smooth diffeomorphisms 
$\Theta_t:\Sigma \longrightarrow \Sigma$, for $t\geq 0$, 
which coincides with the above constructed family  
$\{\Phi_{j_{k^*}} \circ \theta^{j_{k^*}}_t\}$ at every 
$t \geq t_{j_{k^*}}$ - no matter how $m$ was chosen in the 
formulation of the asserted theorem respectively how large 
the index $k^*=k^*(m)$ had to be chosen before formula  
\eqref{full.convergence.m} - such that the reparametrization of the considered flow line $\{F_{t}\}_{t\geq 0}$ 
of the MIWF by $\{\Theta_t\}_{t\geq 0}$ 
converges as in \eqref{full.convergence.m}: 
$$  
F_{t} \circ \Theta_t \longrightarrow G  
\quad \textnormal{in} \,\,\,C^m(\Sigma,\rel^4) 
$$
fully as $t\to \infty$ and at an exponential rate, 
simultaneously for any chosen $m \in \nat$, where $G$ 
smoothly parametrizes some conformally transformed 
Clifford torus in $\sphere^3$.  

\section{Existence of singularities} \label{Singular.flow.lines.do.exist} 

As announced in the introduction, we arrive here at our 
counterpart of Lemma 3.8 in \cite{Dall.Acqua.Schaetzle.Mueller.2024} respectively 
of Theorems 4.1 and 5.1 in \cite{Blatt.2009}, where 
``singularities'' - more precisely of divergent flow lines - of the elastic energy flow in the upper half-plane respectively of the classical Willmore flow were 
constructed respectively detected. However, we won't   
be able here to prove an optimal energy threshold 
which would distinguish between convergent and 
possibly divergent flow lines of \eqref{elastic.energy.flow}, as in \cite{Dall.Acqua.Schaetzle.Mueller.2024} and \cite{Blatt.2009}.    
\begin{theorem} \label{Singularities.do.exist} 
There is some smooth and regular path 
\footnote{Here we are going to construct some 
appropriate closed initial curve $\gamma^*$ 
whose elastic energy is $8 \pi -\varepsilon$, for 
some arbitrarily small $\varepsilon >0$.}  
$\gamma^*:\sphere^1 \longrightarrow \sphere^2$ 
such that the corresponding flow line 
$\{\gamma_t\}_{t \in [0,T_{\textnormal{Max}}(\gamma^*))}$ of the degenerate elastic energy flow \eqref{elastic.energy.flow} with $\gamma_0=\gamma^*$ cannot be global and additionally subconverge smoothly - up to smooth reparametrization - 
to some elastic curve $\gamma_{\infty}$ on $\sphere^2$, 
i.e. there is some flow line 
$\{\gamma_t\}_{t \in [0,T_{\textnormal{Max}}(\gamma^*))}$ 
of flow \eqref{elastic.energy.flow} for which there cannot 
hold $T_{\textnormal{Max}}(\gamma^*) = \infty$ and additionally exist some sequence $t_j \nearrow \infty$ and smooth diffeomorphisms $\varphi_j:\sphere^1 \stackrel{\cong}\longrightarrow  \sphere^1$ such that: 
\begin{equation} \label{smooth.subconvergence.1} 
\gamma_{t_j} \circ \varphi_j 
\longrightarrow \gamma_{\infty} \quad  
\textnormal{in}	\,\,\,C^k(\sphere^1,\rel^3) \,\,\,
\textnormal{for every}\,\, k \in \nat,  
\end{equation}  
as $j \to \infty$, where $\gamma_{\infty}:\sphere^1 \longrightarrow \sphere^2$ is a closed elastic curve on $\sphere^2$.    
\end{theorem}   
\emph{Proof:\,}   
We will explicitly construct some appropriate smooth 
initial curve 
$\gamma^*:\sphere^1 \longrightarrow \sphere^2$ such 
that the corresponding unique flow line 
$\{\gamma_t\}_{t \in [0,T_{\textnormal{Max}}(\gamma^*))}$ 
of the degenerate elastic energy flow \eqref{elastic.energy.flow} cannot be global and additionally satisfy \eqref{smooth.subconvergence.1} with $\gamma_{\infty}$ being a closed elastic curve on $\sphere^2$, for any choice of divergent times $t_j \nearrow \infty$ and of smooth diffeomorphisms $\varphi_j:\sphere^1 \stackrel{\cong}\longrightarrow \sphere^1$. In order to evolve the key ideas chronologically one firstly has to remember the fact that the \emph{non-geodesic} closed elastic curves on $\sphere^2$ have been exactly classified - up to rotations and reflections of $\sphere^2$ - in the first part of Proposition 6 of \cite{Ruben.MIWF.II}, following essentially the classification of \emph{free} elastica on $\sphere^2$ in Section 3 of \cite{Langer.Singer.1984}.
On account of this classification every \emph{non-geodesic} elastic curve on $\sphere^2$ can be characterized by its number $n \in \nat$ of consequtive lobes and 
its number $m\in \nat$ of trips along some fixed 
great circle on $\sphere^2$ before it finally closes up, 
where $m$ and $n$ have to be coprime positive integers 
and satisfy $\frac{m}{n}\in (0,2-\sqrt{2})$. 
We therefore adopt the notation ``$\gamma_{(m,n)}$'' from Proposition 6 in \cite{Ruben.MIWF.II}, where $\gamma_{(m,n)}$
represents the isometry class of all elastic curves 
on $\sphere^2$ having $n$ consequtive lobes while  
they perform exactly $m$ trips along some fixed great 
circle on $\sphere^2$. 
Now, in a first attempt to determine an appropriate 
initial curve $\gamma^*$ verifying the assertion of   
Theorem \ref{Singularities.do.exist} it appeared to be 
a natural idea to pick some sufficiently small $m^*>1$, 
e.g. simply $m^*=2$, then to determine the unique 
elastic curve $\gamma_{(m^*,n^*)}$ with minimal elastic 
energy within the countable set of all  
elastic curves $\gamma_{(m^*,n)}$ having trip number $m=m^*$,  and to perturb exactly this curve slightly in such a way that 
the elastic energy $\Wil$ strictly decreases.
A second, more subtle approach would be here to 
focus rather on the \emph{geodesic} but on the 
\emph{non-geodesic} elastic curves $\gamma_{(m,n)}$  
and hence to slightly perturb some particular $f$-fold cover 
of the equator in such a way that the elastic energy $\Wil$ strictly decreases. 
As we will see below, the first approach fails 
but at least leads to Corollary \ref{minimal.energy.fixed.m}, playing an important technical role later on in this proof,   whereas the second approach indeed works out, taking here
exactly $f=4$. Now, in order to make this strategy work, we will need two key-ingredients. 
(a) Some precise, rather technical computations arising in the proof of the second part of Proposition 6 in \cite{Ruben.MIWF.II}, (b) Langer's and Singer's insight \cite{Langer.Singer.1987}
\footnote{The two assertions quoted here are actually 
only special cases of Langer's and Singer's Theorem 3.1 
in \cite{Langer.Singer.1987}, where instability 
of closed elastica on spheres has been 
precisely investigated.}  
that each \emph{non-geodesic} elastic curve on $\sphere^2$
- performing only one loop through its entire trace -  
and additionally the $f$-fold cover of any great circle, 
for each $f\geq 4$, is an \emph{unstable} critical point of the elastic energy $\Wil$. 
Regarding ingredient (b) we should recall more precisely 
that for each \emph{non-geodesic} elastica 
$\gamma:\sphere^1 \longrightarrow \sphere^2$  
a geometrically natural choice of variation 
$\vec F_{\gamma}$ along $\gamma$ which slightly 
reduces the elastic energy in $\gamma$ is of the form 
$\phi \, \vec \kappa_{\gamma}$, for 
any smooth non-vanishing function $\phi:\sphere^1\longrightarrow \rel$. In other words, 
the vector field $\phi \, \vec \kappa_{\gamma}$ along 
$\gamma$ satisfies:
\begin{equation}  \label{Unstable.elastic.curves}
(\delta^2\Wil)_{\gamma}
(\vec F_{\gamma},\vec F_{\gamma}) <0.
\end{equation} 
In order to achieve inequality \eqref{Unstable.elastic.curves} for some $f$-fold cover 
$\gamma:=E \oplus E \oplus \ldots \oplus E$ of some fixed 
great circle $E$, for any fixed $f\geq 4$,
one can choose one of the two unit normal vector fields 
$N_{E}$ along the chosen equator $E$ and consider 
a particular smooth section 
$\vec F_{\gamma}:=
\phi_1\, N_{E} \oplus \phi_2 N_{E} \oplus \ldots 
\oplus \,\phi_f\, N_{E} \in \Gamma(\gamma^*T\sphere^2)$ which  is normal along the $f$-fold cover $\gamma$ of the great circle $E$, for some suitable collection of $f$ smooth 
functions $\phi_j: \sphere^1\setminus \{1\} \longrightarrow \,\rel$, as explained and proved on p. 147 in \cite{Langer.Singer.1987}.   
Regarding ingredient (a) we recall here from Section 4 in 
\cite{Langer.Singer.1987} and from the proof of Proposition 6 in \cite{Ruben.MIWF.II} that the wavelength 
$\Lambda(\gamma_{(m,n)})=\frac{m}{n}2 \pi$ of the 
closed non-geodesic elastica ``$\gamma_{(m,n)}$'' can be 
expressed as the function (102) in \cite{Ruben.MIWF.II} 
of the modulus $p \in \big{(}0,\frac{1}{\sqrt{2}}\big{)}$ of the Jacobi elliptic function appearing in formula (97) of \cite{Ruben.MIWF.II} - concretely parametrizing 
the squared curvature along the corresponding   
elastic curve ``$\gamma_{(m,n)}$'' or ``$\gamma(p(m,n))$'' - 
and this function $p\mapsto \Lambda(\gamma(p))$ is strictly monotonically decreasing by formula (103) in \cite{Ruben.MIWF.II}. This does not only give us a  
one-to-one correspondence between all quotients 
$\frac{m}{n}\in (0,2-\sqrt{2})$ with $\textnormal{gcd}(m,n)=1$ 
and exactly those moduli $p=p(m,n)\in \big{(}0,\frac{1}{\sqrt{2}}\big{)}$ which 
produce the squared curvatures of all 
\emph{non-geodesic closed elastic curves} on $\sphere^2$, 
but it also gives rise to the possibility to compute 
the energy $\Wil$ of any chosen non-geodesic elastic curve $\gamma_{(m,n)}$ in terms of certain elliptic integrals depending on the unique, corresponding modulus $p=p(m,n)$. Indeed, in formula (106) of \cite{Ruben.MIWF.II} the author computed exactly: 	
\begin{equation} \label{energy.directly}    
\Wil(\gamma_{(m,n)}) = \frac{8n}{\sqrt{2-4(p(m,n))^2}} \,\big{(}2E(p(m,n))-K(p(m,n))\big{)}.  
\end{equation} 	              
Moreover, on p. 38 of \cite{Ruben.MIWF.II} it is 
explained why the function  
$f(p):=\frac{1}{\sqrt{1-2p^2}} \,(2E(p)-K(p))$ - 
appearing in formula \eqref{energy.directly} above - 
is strictly monotonically increasing on its domain $\big{(}0,\frac{1}{\sqrt{2}}\big{)}$. 
These facts give rise to the following important 
\begin{corollary} [Vertical energy pattern of the $\{\gamma_{(m,n)}\}$-scheme] \label{minimal.energy.fixed.m} 
Let some $\bar m \in \nat$ be fixed. 
The sequence of numbers 
$\{\Wil(\gamma_{(\bar m,n)})\}_{n\geq 2}$ is monotonically 
increasing, where we only consider those indices  
$n\geq 2$ for which the elastic curves 
$\gamma_{(\bar m,n)}$ actually exist.   
Particularly, for every fixed integer $\bar m \geq 1$
the subset of $\rel_+$ consisting of the elastic energies 
of all non-geodesic closed elastic curves on $\sphere^2$ 
with fixed trip number $\bar m$ contains its infimum, 
and this minimal value is $\Wil(\gamma_{(\bar m,\bar n)})$, where $\bar n$ is the minimal number $n$ of lobes 
which is coprime with $\bar m$ and still satisfies 
$\frac{\bar m}{n}<2-\sqrt{2}$. Expressed in formulae this means: 
	\begin{equation} \label{n.minimal.energy.fixed.m}
	\inf \Big{\{} \Wil(\gamma_{(\bar m,n)}) \,\Big{|}\,
	n\in \nat \,\wedge \, \textnormal{gcd}(\bar m,n)=1 \,\wedge \,\frac{\bar m}{n} <2-\sqrt{2}  \,\Big{\}} = \Wil(\gamma_{(\bar m,\bar n)}), 
	\end{equation} 
where exactly $\bar n=\min\big{\{}n\in \nat \,|\,
\textnormal{gcd}(\bar m,n)=1 \,\wedge \, 
\frac{\bar m}{n} <2-\sqrt{2} \,\big{\}}$.
\end{corollary}  	 
\emph{Proof:\,}  
For fixed $\bar m \in \nat$ the wavelength $\Lambda(\gamma_{(\bar m,n)})=\frac{\bar m}{n} 2 \pi$ 
obviously increases strictly monotonically as 
$n$ drops down to its minimal possible value $\bar n$, 
such that still $\frac{\bar m}{n}<2-\sqrt{2}$ 
and $\textnormal{gcd}(m,n)=1$ hold, as formulated  
in \eqref{n.minimal.energy.fixed.m}. Hence, by formulae (102) and (103) in \cite{Ruben.MIWF.II} the corresponding modulus  
$p=p(\bar m,n)\in \big{(}0,\frac{1}{\sqrt{2}}\big{)}$ 
has to drop strictly monotonically 
$p(\bar m,n) \searrow p(\bar m,\bar n)$ as 
$n \searrow \bar n$. Now, combining this result with the fact that the function $f(p)=\frac{1}{\sqrt{1-2p^2}} \,(2E(p)-K(p))$ actually increases strictly monotonically 
on $\big{(}0,\frac{1}{\sqrt{2}}\big{)}$, we finally infer: 
$$  
f(p(\bar m,n)) \searrow f(p(\bar m,\bar n)) \quad  
\textnormal{as} \,\,\,n \searrow \bar n. 
$$     
Consequently, also the product 
$\frac{8n}{\sqrt{2}}\,f(p(\bar m,n))$ drops strictly monotonically to the value 
$\frac{8\bar n}{\sqrt{2}}\,f(p(\bar m,\bar n))$ 
as \,$n \searrow \bar n$. Expressed in terms of formula 
\eqref{energy.directly} this means: 
$$    
\Wil(\gamma_{(\bar m,n)}) =  
\frac{8n}{\sqrt{2-4(p(\bar m,n))^2}} \,
\big{(}2E(p(\bar m,n))-K(p(\bar m,n))\big{)}
\searrow \Wil(\gamma_{(\bar m,\bar n)}) \quad  
\textnormal{as} \,\,\,n \searrow \bar n.
$$  
Since - up to isometries of $\sphere^2$ - there aren't any more non-geodesic closed elastic curves on $\sphere^2$ than the particular elastica $\gamma_{(m,n)}$ appearing in the classification of Proposition 6 in \cite{Ruben.MIWF.II}, 
the set of elastic energies of all non-geodesic closed elastic 
curves on $\sphere^2$ with fixed trip number $\bar m$
is exactly the countable set 
$\big{\{} \Wil(\gamma_{(\bar m,n)}) \,\Big{|}\,
n\in \nat \,\wedge \, \textnormal{gcd}(\bar m,n)=1 \,\wedge \,\frac{\bar m}{n} <2-\sqrt{2}  \,\big{\}}\subset \rel_+$. Hence, the previous computation already proves both 
assertions of the corollary.  	    \\\\
\noindent 
Now we pick up our key ingredient (b) and start proving Theorem \ref{Singularities.do.exist}. 
We try to perform a topological contradiction argument, 
and in order to keep this argument as simple as possible, we should focus on some $\Wil$-unstable elastica in $\sphere^2$ 
which travels more than only once about $\sphere^2$  
along some fixed equator, but still has the least possible 
elastic energy. We claim that {\bf the $4$-fold equator} is exactly this desired path, i.e. the path $\gamma^*=E\oplus E \oplus E \oplus E$ which traverses some fixed equator $E$ with constant speed exactly $4$ times. First of all,  
this critical point of $\Wil$ is unstable on account of statement \eqref{Unstable.elastic.curves}. Actually 
statement \eqref{Unstable.elastic.curves} tells us 
that each \emph{non-geodesic} elastic curve $\gamma_{(m,n)}$ is $\Wil$-unstable as well, but we can easily prove here that most of these candidates have much higher elastic energies than our chosen $\gamma^*$ - a fact which will be very important below.   
We estimate for each fixed number $m^*>4$ of trips 
about some fixed great circle in $\sphere^2$
the elastic energy $\Wil(\gamma_{(m^*,n)})$ from 
below by means of formula \eqref{energy.directly}, 
combined with the monotonicity of the function $f(p)=\frac{1}{\sqrt{1-2p^2}} \,(2E(p)-K(p))$ 
on $\big{(}0,\frac{1}{\sqrt{2}}\big{)}$ - as mentioned 
above - and with the particular facts that 
$f(p)\geq f(0)= \frac{\pi}{2}$ and that each admissible 
pair of numbers $(m,n)$ satisfies:
$0<\frac{m}{n}<2-\sqrt{2}$. On account of the latter fact 
we have $n>\frac{3}{2}\,m^*$ for any admissible $n$, 
and thus formula \eqref{energy.directly} yields:
\begin{equation} \label{Wil.estimate.below}     
\Wil(\gamma_{(m^*,n)}) 
> \frac{8\cdot \frac{3}{2}\,m^*}{\sqrt{2}} \,
f(p(m^*,n)) \geq \frac{6\,m^*}{\sqrt{2}} \,\pi 
\geq \frac{30}{\sqrt{2}} \,\pi \approx 66.64324, 
\end{equation} 
for every admissible pair $(m^*,n)$ and $m^*>4$,  
which is obviously much larger than 
$\Wil(\gamma^*)=8\pi \approx 25.13274$. Combining this 
insight with Corollary \ref{minimal.energy.fixed.m} for 
$m=1,2,3,4$ and with the concrete minimal values  
\begin{eqnarray} \label{special.values}
\Wil(\gamma_{(1,2)})\approx 19.17, \quad 
\Wil(\gamma_{(2,5)})\approx 55.01, \quad  
\Wil(\gamma_{(3,7)})\approx 74.97, \quad 
\Wil(\gamma_{(4,7)})\approx 62.89
\end{eqnarray}
\footnote{See the author's computations on p. 40 in \cite{Ruben.MIWF.II}.} 
of energies of elastic curves with $m=1,2,3,4$ 
trips about some fixed equator,  
we see that indeed the $4$-fold equator with energy about $25.133$ is the cheapest choice of an elastic curve
which is $\Wil$-unstable and performs more than only 
one single trip about $\sphere^2$ along some fixed great circle. Now we apply the statement below inequality \eqref{Unstable.elastic.curves} to the $4$-fold cover
$\gamma^*:= E \oplus E \oplus E \oplus E$ of the 
chosen great circle $E$. We therefore choose 
one of the two unit normals $N_{E}$ along the  
great circle $E$ within $T\sphere^2$, and we consider 
an unstable variational section 
$\vec F_{\gamma^*}:= 
\phi_1\, N_{E} \oplus \phi_2 N_{E} \oplus  
\phi_3\, N_{E} \oplus \,\phi_4\, N_{E}$ 
of the normal bundle - within $\Gamma((\gamma^*)^*T\sphere^2)$ - of the $4$-fold cover $\gamma^*$ of the great circle $E$, 
for some particular collection of $4$ smooth 
functions $\phi_j: \sphere^1\setminus \{1\} \longrightarrow \,\rel$, as introduced below inequality \eqref{Unstable.elastic.curves}. Since the value of the energy 
$\varepsilon \mapsto 
\Wil(\exp_{\gamma^*}(\varepsilon \vec F_{\gamma^*}))$ has 
to decrease below $\Wil(\gamma^*)$ for 
$\varepsilon \in (0,\varepsilon_0)$, provided 
$\varepsilon_0>0$ is sufficiently small,   
we should focus on the smooth and regular closed curves 
$\{\gamma^{\varepsilon}\}_{\varepsilon \in [0,\varepsilon_0)}
:=\{\exp_{\gamma^*}(\varepsilon \vec F_{\gamma^*})\}
_{\varepsilon \in [0,\varepsilon_0)}$ 
with $\gamma^0 \equiv E\oplus E \oplus E \oplus E=\gamma^*$. Hence, these paths travel $4$-times about $\sphere^2$ along 
the chosen great circle $E$ - just as $\gamma^*$ itself does - and their energies satisfy: 
\begin{equation} \label{distorted.curve} 
\Wil(\gamma^{\varepsilon}) <\Wil(\gamma^*) = 8 \pi
\quad \forall \, \varepsilon \in (0,\varepsilon_0),
\end{equation} 
for some sufficiently small $\varepsilon_0>0$. 
Here and in the sequel we silently use the fact that the direction $N_{E}$ of the distortion along $E$ appearing in 
\eqref{distorted.curve} is a smooth section of the pullback bundle $(\gamma^*)^*T\sphere^2$. 
On account of inequality \eqref{distorted.curve} 
we shall consider now the unique flow line 
$\{\gamma^{\varepsilon}_t\}$ of flow 
\eqref{elastic.energy.flow} starting at such 
a distorted closed path $\gamma^{\varepsilon}$, for some 
arbitrarily chosen $\varepsilon \in (0,\varepsilon_0)$. \\ 
Now we shall assume by contradiction that $\{\gamma^{\varepsilon}_t\}_{t \geq 0}$ would exist globally   
and subconverge to an elastic curve in $\sphere^2$. 
Precisely we assume that there are $t_j \nearrow \infty$ 
and smooth diffeomorphisms $\varphi_j:\sphere^1 \stackrel{\cong}\longrightarrow \sphere^1$
such that \eqref{smooth.subconvergence.1} 
holds for any $k \in \nat$, where the limit path 
$\gamma^{\varepsilon}_{\infty}:\sphere^1 \longrightarrow \sphere^2$ parametrizes a smooth elastic curve on $\sphere^2$, 
possibly performing several loops through its trace ! 
Since $\gamma^{\varepsilon}_0$ satisfies \eqref{distorted.curve}, we could infer from the monotonicity 
of the elastic energy along the global flow line 
$\{\gamma^{\varepsilon}_t\}$ that  
\begin{equation}  \label{Wil.goes.down}
\Wil(\gamma^{\varepsilon}_{\infty})
\leq \Wil(\gamma^{\varepsilon}) < \Wil(\gamma^*)=8\pi,
\end{equation} 
which shows again on account of Corollary 
\ref{minimal.energy.fixed.m}, estimate \eqref{Wil.estimate.below} and the values in 
\eqref{special.values} that our particular 
elastic limit curve $\gamma^{\varepsilon}_{\infty}$ in \eqref{smooth.subconvergence.1} has to be 
one of the following $4$ curves: 
$$
(A):\,\, \gamma_{(1,2)}\quad  (B): \,\, E \quad 
(C): \,\, E \oplus E \quad (D): \,\, E \oplus E\oplus E,
$$ 
at least up to appropriate isometries of $\sphere^2$.
Now there are the following two approaches. The first 
alternative would be to introduce and use here 
\underline{Arnold's} \cite{Arnold.1993} 
\underline{homotopy invariants} on the set     $\textnormal{Imm}_{\textnormal{gen}}(\sphere^1,\sphere^2)$ 
of generic \footnote{A regular closed path  
$\gamma:\sphere^1 \longrightarrow \rel^2$ is called  
\emph{generic}, iff its self-intersections are 
transversal and of multiplicity at most $2$, 
i.e. transversal ``double points''.
See here p. 1 in \cite{Arnold.1993} or Section 2.1 in \cite{Polyak.1996}.}   
closed paths in $\sphere^2$ - an open subset of the more familiar set $\textnormal{Imm}(\sphere^1,\sphere^2)$ of all 
smooth and regular closed paths in $\sphere^2$. 
The second alternative would be to rely on one's 
optimism that our degenerate elastic energy flow 
\eqref{elastic.energy.flow} - being started in 
$\gamma^{\varepsilon}$ - would not touch every point of the $2$-sphere but only points on the punctured sphere 
$\sphere^2\setminus \{b\}$ for some suitably chosen 
``base point'' $b$ and that the stereographic projection 
$\Psi_b:\sphere^2 \setminus \{b\} \longrightarrow \rel^2$  
of the entire flow line $\{\gamma^{\varepsilon}_{t}\}$
into $\rel^2$ together with \underline{Whitney-Graustein's 
tangent-rotation number} ``$\textnormal{ind}$'' 
on $\textnormal{Imm}(\sphere^1,\rel^2)$ 
will suffice, in order to produce the desired contradiction. 
We will see below that these two alternatives have to be combined with each other !  \\          
Let's start discussing the \emph{first alternative}, 
thus motivating our application of 
\emph{Arnold's invariants} on the set  $\textnormal{Imm}_{\textnormal{gen}}(\sphere^1,\rel^2)$
of generic closed paths in $\rel^2$. 
These invariants decide whether two different generic closed paths in $\rel^2$ can be connected by a continuous family of generic closed paths, and moreover they  
measure the ``distance'' between the corresponding equivalence 
classes algebraically.
Furthermore, as explained in Section 2.4 of \cite{Polyak.1996} 
one can use these three invariants on $\textnormal{Imm}_{\textnormal{gen}}(\sphere^1,\rel^2)$ 
in order to build three canonical   
\emph{generic homotopy invariants} 
$J^{\pm}_{\sphere^2}$, $\textnormal{St}_{\sphere^2}$ 
mapping $\textnormal{Imm}_{\textnormal{gen}}(\sphere^1,\sphere^2)$ 
into $\nat_0$, and these invariants are explicitly expressed  
in Section 2.4 of \cite{Polyak.1996} as follows: 
\begin{equation} \label{invariants}  
J^{\pm}_{\sphere^2}(\gamma) 
= J^{\pm}_{a}(\gamma) + \frac{1}{2} \big{(}\textnormal{ind}_a(\gamma)\big{)}^2 \quad \vee \quad  \textnormal{St}_{\sphere^2}(\gamma) 
= \textnormal{St}_a(\gamma) - 
\frac{1}{4}\cdot \big{(}\textnormal{ind}_a(\gamma)\big{)}^2.  
\end{equation}  
Here $a$ denotes some \underline{arbitrarily fixed} point on 
$\sphere^2 \setminus \textnormal{trace}(\gamma)$ serving as the ``north pole'' on $\sphere^2$ which 
is sent to ``infinity'' by some appropriate stereographic 
projection $\Psi_a:\sphere^2\setminus \{a\} \longrightarrow \rel^2$. Obviously we need some canonical way to
replace the given path $\gamma$ by some diffeomorphic image 
in $\rel^2$ - e.g. by $\Psi_a \circ \gamma$ - such that   
the invariants $J^{\pm}, \textnormal{St}:
\textnormal{Imm}_{\textnormal{gen}}(\sphere^1,\rel^2)
\longrightarrow \ganz$   
together with Arnold's entire theory \cite{Arnold.1993} 
classifying smooth generic closed curves in $\rel^2$
become applicable. 
Hence, in the sequel Arnold's intuitive notation  $\textnormal{ind}_a(\gamma)$, $J^{\pm}_{a}(\gamma)$ and $\textnormal{St}_a(\gamma)$ means simply 
$\textnormal{ind}(\Psi_a\circ \gamma)$, 
$J^{\pm}(\Psi_a\circ \gamma)$ and      
$\textnormal{St}(\Psi_a\circ \gamma)$, using here some suitable stereographic projection 
$\Psi_a:\sphere^2 \setminus \{a\} \longrightarrow \rel^2$, 
which is uniquely determined by the point $a$ up to some 
rotation about the line through $a$ and the origin.\\
On the one hand Arnold's original variants $J^{\pm}$ 
and $\textnormal{St}$ on $\textnormal{Imm}_{\textnormal{gen}}(\sphere^1,\rel^2)$  
were characterized in Paragraph 1 of 
\cite{Arnold.1993} in terms of their exact behaviour 
in three types of \emph{non-generic transgressions} or \emph{perestroikas}; see also Section 2.3 of 
\cite{Polyak.1996}. On the other hand, on account of formula 
\eqref{invariants} and the entire discussion in 
Paragraphs 1, 3, 6, 7 of \cite{Arnold.1993} we can see 
that the invariant $\textnormal{St}_{\sphere^2}(\gamma)$ - typically called ``strangeness'' - can be explicitly computed 
for some chosen path $\gamma \in \textnormal{Imm}_{\textnormal{gen}}(\sphere^1,\sphere^2)$, 
if one knows how to choose some point 
$a \in \sphere^2\setminus \textnormal{trace}(\gamma)$ 
in such a way that the index 
$\textnormal{ind}(\Psi_a \circ \gamma)$ of the
stereographically projected, planar path 
$\Psi_a \circ \gamma$ is well-known - or at least 
computable - and if one can transform the generic planar path 
$\Psi_a \circ \gamma$ regularly homotopically into 
one of the so-called \emph{standard curves} 
$K_{\textnormal{ind}(\Psi_a \circ \gamma)}$ or 
$A_{\textnormal{ind}(\Psi_a \circ \gamma)}$ 
in $\rel^2$ without experiencing any \emph{triple point perestroikas} during such a homotopy; 
see here also Figures 5, 27 and 56 in \cite{Arnold.1993}. \\ 
Now, in our particular situation we can 
indeed calculate the three invariants 
appearing in \eqref{invariants} of any elastic 
curve appearing in (A)--(D) and of the initial curve $\gamma^{\varepsilon}_{0}$. First of all, we can choose the ``base-point'' $a\in \sphere^2$ appearing in formulae \eqref{invariants} as one of the two intersection points of the symmetry axis of the equator $E$ respectively of the elastic curve $\gamma_{(1,2)}$, and thus we can apply 
stereographic projection into $\rel^2$ in a fairly 
natural way. The $f$-fold equator is obviously not 
generic for any $f>1$, but still we can slightly perturb its
stereographically projected image in $\rel^2$, 
in order to obtain a  
``\emph{close generic approximation}'' to this projected, 
planar curve in $\textnormal{Imm}_{\textnormal{gen}}(\sphere^1,\rel^2)$
in a canonical way, 
namely the standard planar curve $A_{f}$ of index $f$ 
and maximal strangeness-value: 
$\textnormal{St}(A_f) = \frac{(f-1)\,f}{2}$; 
see here Figure 27 in Paragraph 7 of \cite{Arnold.1993}. 
This means precisely that in every $C^2(\sphere^1,\rel^2)$-neighbourhood 
of the projected $f$-fold equator - say parametrized with 
constant speed - there is some smooth generic path
which is regularly homotopic to some smooth regular parametrization of the standard curve $A_{f}$ 
- having rotation-index $\pm f$ and maximal strangeness value 
$\frac{(f-1)\,f}{2}$ - without producing any type of \emph{perestroika} during a suitably chosen homotopy. 
\footnote{See here Section 2.3 in \cite{Polyak.1996} and 
the first two paragraphs in \cite{Arnold.1993}.}    
On the other hand for any fixed $f\geq 1$ 
there is some sufficently small open $C^2(\sphere^1,\rel^2)$-neighbourhood about the 
stereographically projected, arc-length parametrized 
$f$-fold equator, such that in this neighbourhood 
there is no generic path which is regularly homotopic 
to some smooth regular parametrization of a standard curve 
- for example $A_{p}$ - whose rotation-index is $\pm p$ 
for some $p \not= f$. 
We could also propose to choose another type of 
standard curve in $\textnormal{Imm}_{\textnormal{gen}}(\sphere^1,\rel^2)$ 
being $C^2$-close to the stereographic projection of the $f$-fold equator, but any such ``generic representative'' 
would have rotation index $f$ ! Hence, on account of  
Whitney-Graustein's famous index-theorem - proving that the 
\emph{tangent-rotation number} ``$\textnormal{ind}$'' 
is the unique invariant on entire
$\textnormal{Imm}(\sphere^1,\rel^2)$ 
\footnote{See here for example page 5 in \cite{Arnold.1993}.}-
we infer that there is some regular homotopy in $\rel^2$ between any two choices of 
``\emph{close generic approximations}'' 
to the stereographic projection of the $f$-fold equator, 
and we therefore deduce from Theorem 1 in \cite{Arnold.1993}
that the strangeness values $\textnormal{St}$ of two  
such choices may only differ by a whole number, 
counting the number of - either positive or negative - triple-point-perestroikas during any such 
regular homotopy. We may therefore sloppily interpret 
the $f$-fold equator $E \oplus E \oplus \ldots \oplus E$ 
in $\sphere^2$ for any chosen $f\geq 1$ as an element of  $\textnormal{Imm}_{\textnormal{gen}}(\sphere^1,\sphere^2)$ 
having a well-defined strangeness value 
\footnote{See here also Paragraph 7, particularly pp. 29--31, in \cite{Arnold.1993}.} modulo $\ganz$:
\begin{equation}  \label{St.sphere.f.fold.Equat}
\textnormal{St}_{\sphere^2}(E \oplus E \oplus \ldots \oplus E)
\equiv \Big{[}\textnormal{St}(A_f) 
- \frac{1}{4} \cdot \big{(}\textnormal{ind}(A_f)\big{)}^2 \Big{]} \,\, \textnormal{mod} \, \ganz 
= \frac{f^2-\,2f}{4} \,\, \textnormal{mod} \, \ganz,
\end{equation} 
on account of the second formula in \eqref{invariants} 
and Theorem 1 in \cite{Arnold.1993}. 
The possibility to attach a unique strangeness value 
to the $f$-fold equator in $\sphere^2$ at least modulo $\ganz$ - namely $\frac{f^2-\,2f}{4} \, \textnormal{mod} \, \ganz$ 
by formula \eqref{St.sphere.f.fold.Equat} -  
will be useful in the sequel in order to rule out the 
two possibilities (A) and (B) for the limit curve  
$\gamma^{\varepsilon}_{\infty}$ of the considered flow 
line $\{\gamma^{\varepsilon}_t\}$, whereas possibilities 
(C) and (D) will have to be taken care of individually.    \\
Now, following exactly the same considerations we can slightly perturb the initial curve 
$\gamma^{\varepsilon} = 
\exp_{\gamma^*}(\varepsilon \vec F_{\gamma^*})$ of our 
candidate flow line $\{\gamma^{\varepsilon}_{t}\}$, 
in order to obtain a smooth \emph{generic} path, 
which still satisfies inequality \eqref{distorted.curve}
and whose stereographic projection into $\rel^2$ is 
regularly homotopic to some closed generic curve - 
as e.g. $A_4$ or $K_4$ - without producing any type of \emph{perestroikas} during a suitably chosen homotopy. 
Hence, we may assume that 
\begin{equation}  \label{index.gamma.0.a}  
\textnormal{ind}_a(\gamma^{\varepsilon}_{0}) = 
\textnormal{ind}_a(E \oplus E \oplus E \oplus E) = \pm 4  
\end{equation}
holds for a generic point $a \in \sphere^2 \setminus \textnormal{trace}(\gamma^{\varepsilon}_{0})$,
and that there holds additionally by \eqref{invariants} and Theorem 1 in \cite{Arnold.1993}: 
\begin{equation} \label{invariants.gamma.epsilon.0}  
\textnormal{St}_{\sphere^2}(\gamma^{\varepsilon}_{0}) 
\equiv \Big{[}\textnormal{St}(A_4) 
- \frac{1}{4} \cdot \big{(}\textnormal{ind}(A_4)\big{)}^2 \Big{]} \, \textnormal{mod}\,\, \ganz 
\equiv  \frac{4^2 -\,2\cdot 4}{4} \, \textnormal{mod}\,\, \ganz \equiv 0 \,\, \textnormal{mod}\,\, \ganz, 
\end{equation} 
similarly to our computation in formula \eqref{St.sphere.f.fold.Equat}.
Just as in formula \eqref{St.sphere.f.fold.Equat} 
the mod $\ganz$-value obtained in \eqref{invariants.gamma.epsilon.0}
only depends on $\gamma^{\varepsilon}_{0}$, but not 
on any further, arbitrary choices. 
Applying the above ``perturbative argument'', which
resulted in the general formula \eqref{St.sphere.f.fold.Equat}, especially to the $4$ remaining possibilities (A)--(D) for the limit elastic curve $\gamma^{\varepsilon}_{\infty}$ we easily obtain:
\[ \textnormal{St}_{\sphere^2}(\gamma^{\varepsilon}_{\infty}) 
\equiv \Big{[}\textnormal{St}(A_{\textnormal{ind}_a
(\gamma^{\varepsilon}_{\infty})}) - 
\frac{1}{4} \cdot \big{(}\textnormal{ind}_a(\gamma^{\varepsilon}_{\infty})
\big{)}^2 \Big{]} \, \textnormal{mod} \,\ganz
= \left\{ \begin{array}{r@{\quad:\quad}l} 
 \frac{3}{4} \, \textnormal{mod} \,\ganz   &  
 \textnormal{Case} \,\, (A), (B), (D)   \\ 
0 \, \textnormal{mod} \,\ganz    &     \gamma^{\varepsilon}_{\infty} = E \oplus E       
\end{array} \right.    \]           
Now, since the considered flow line $\{\gamma^{\varepsilon}_t\}$ is assumed 
to be global and since each single curve $\gamma^{\varepsilon}_t$ is immersed 
- by definition of a flow line - the restriction 
$\{\gamma^{\varepsilon}_t\}_{t \in [0,T]}$ of the 
entire flow line can obviously be interpreted as a regular homopoty $\gamma^{\varepsilon}_t:\sphere^1 \times [0,T]  
\longrightarrow \sphere^2$, for every positive $T$.
Now suppose, we are given some arbitrary regular homotopy 
$H:[0,1]\times \sphere^1 \longrightarrow \sphere^2$
between two generic closed curves $H(0,\,\cdot\,)$ 
and $H(1,\,\cdot\,)$ on $\sphere^2$. 
Recalling the definition of $\textnormal{St}_{\sphere^2}$  
in \eqref{invariants}, one sees that the number of triple-point-perestroikas during the given homotopy $H$ would only be correctly computed by the difference between  
$\textnormal{St}_{\sphere^2}(H(0,\,\cdot\,))$ and 
$\textnormal{St}_{\sphere^2}(H(1,\,\cdot\,))$, if one 
is able to remove at least one suitable point $b$ from $\sphere^2$ which is not contained in the image of $H$, 
i.e. if $H$ does not cover the entire $2$-sphere.
Otherwise one would not be able to first map the 
entire homotopy $H$ diffeomorphically into $\rel^2$
and then apply all properties of 
Arnold's invariants $J^{\pm}$ and $\textnormal{St}$ on $\textnormal{Imm}_{\textnormal{gen}}(\sphere^1,\rel^2)$ 
- the building blocks of each invariant in \eqref{invariants} - as they are exactly established in Paragraphs 1--4 of \cite{Arnold.1993}.
Let's show here why this \underline{technical obstruction} turns out to be unproblematic in our special setting. Without loss of generality we may assume that the Hopf-torus $\pi^{-1}(\textnormal{trace}(\gamma^{\varepsilon}_0))$ 
is compactly contained in $\sphere^3 \setminus \{(0,0,0,1)\}$,
and we can therefore map
$\pi^{-1}(\textnormal{trace}(\gamma^{\varepsilon}_0))$  
by means of standard stereographic projection $\Stereo$
from $\sphere^3\setminus \{(0,0,0,1)\}$ into $\rel^3$. 
Now we choose some smooth parametrization 
$F_0:\Sigma \longrightarrow  
\Stereo(\pi^{-1}(\textnormal{trace}(\gamma^{\varepsilon}_0)))$,
and we obtain the unique flow line 
$\{F_t\}_{t\in [0,T_{\textnormal{max}}(F_0))}$ of the MIWF in $\rel^3$ starting in $F_0$. Now $T_{\textnormal{max}}(F_0)$ cannot be finite here, because otherwise the corresponding flow line $\{\Stereo^{-1}\circ F_t\}_{t\in [0,T_{\textnormal{max}}(F_0))}$ of the MIWF in $\sphere^3$
would break down at the finite time $T_{\textnormal{max}}(F_0)$ as well, which would violate our assumption on the original flow line $\{\gamma^{\varepsilon}_t\}$ of flow \eqref{elastic.energy.flow} to be global; compare here also with Proposition \ref{correspond.flows}. 
Now, for every $T>0$ the product $[0,T]\times \Sigma$ is compact, and therefore the image of the restriction 
$F:[0,T]\times \Sigma \longrightarrow \rel^3$ 
of the considered flow line $\{F_t\}$ of the MIWF in 
$\rel^3$ has to be compact. 
Hence, the restriction
$\{\Stereo^{-1}\circ F_t\}_{t\in [0,T]}$ of 
the corresponding flow line of the MIWF in $\sphere^3$ 
cannot exhaust all of $\sphere^3$, 
and therefore the restriction 
$\{\gamma^{\varepsilon}_t\}_{t\in [0,T]}$ 
of the original global flow line 
$\{\gamma^{\varepsilon}_t\}_{t\in [0,\infty)}$ 
of flow \eqref{elastic.energy.flow} to any compact 
time interval $[0,T]$ cannot exhaust all of $\sphere^2$
- again because of Proposition \ref{correspond.flows} - 
just as required.     
\footnote{One should remark here that this geometrically motivated argument relies strongly on the conformal invariance of the MIWF and would not work out for the classical elastic energy flow \eqref{classical.elastic.energy.flow}, 
simply because its flow lines cannot be related in a  one-to-one fashion to flow lines of another geometric 
flow whose ambient space is non-compact.} \\   
Now we can choose a suitable point 
$b_T \in \sphere^2\setminus \bigcup_{0\leq t\leq T} 
\textnormal{trace}(\gamma^{\varepsilon}_t)$, and 
we infer from formula \eqref{invariants} respectively 
Section 2.4 in \cite{Polyak.1996}, from equation \eqref{invariants.gamma.epsilon.0}, from 
Whitney's famous index-theorem and from Theorem 1 in \cite{Arnold.1993} respectively 
Property 2.4 in Section 2.3 of \cite{Polyak.1996}: 
\begin{eqnarray} \label{invariants.gamma.epsilon.T}  
\textnormal{St}_{\sphere^2}(\gamma^{\varepsilon}_{T})
=\textnormal{St}_{b_T}(\gamma^{\varepsilon}_{T}) 
- \frac{1}{4} \cdot \big{(}\textnormal{ind}_{b_T}
(\gamma^{\varepsilon}_{T})\big{)}^2=    \nonumber  \\  
=\textnormal{St}_{b_T}(\gamma^{\varepsilon}_{0}) + k_T
- \frac{1}{4} \cdot \big{(}\textnormal{ind}_{b_T}(\gamma^{\varepsilon}_{0})
\big{)}^2  
=\textnormal{St}_{\sphere^2}(\gamma^{\varepsilon}_{0}) + k_T 
\equiv 0 \,\,\, \textnormal{mod}\,\, \ganz
\end{eqnarray}
for every fixed positive time $T$.  
Here, the integer $k_T$ counts the difference between positive 
and negative triple point perestroikas which the  
planar curves $\textnormal{trace}(\Psi_{b_T} \circ \gamma^{\varepsilon}_{t})\subset \rel^2$ encounter, 
as $t$ runs from $0$ to $T$; see here especially pp. 4--5 
in the first paragraph of \cite{Arnold.1993}.   \\ 
Moreover, in the two cases (A) and (B) we know that 
the limit path $\gamma^{\varepsilon}_{\infty}$ is simply closed, i.e. not only a regular path but a smooth 
embedding, and therefore \emph{generic} in particular. 
Hence, in cases (A) and (B) we can guarantee  
that there is some sufficiently large index $J$ - depending 
on the rate of convergence in \eqref{smooth.subconvergence.1} - such that for each $j\geq J$ the reparametrized 
closed path $\gamma^{\varepsilon}_{t_j} \circ \varphi_j$ 
is (a) a smooth embedding and (b) isotopic to the limit   embedding $\gamma^{\varepsilon}_{\infty}$ 
\footnote{In cases (C) and (D) these strong properties  
might be relaxed appropriately, but still the genericity of the paths $\gamma^{\varepsilon}_{t_j}$ would be unclear.}. 
Assertion (a) follows immediately from  
the assumed $C^k$-convergence \eqref{smooth.subconvergence.1} 
of our global flow line $\{\gamma^{\varepsilon}_t\}$ 
and from the fact that \emph{embeddedness} of closed regular 
paths in $\sphere^2$ is an open property with respect to the $C^2(\sphere^1,\rel^3)$-norm, and Assertion (b) can be 
quickly derived from convergence \eqref{smooth.subconvergence.1} by means of the linear homotopy 
$\tilde L_j(x,\tau):=
(1-\tau) \, \gamma^{\varepsilon}_{t_j} \circ \varphi_j(x) 
+ \tau \, \gamma^{\varepsilon}_{\infty}(x)$ 
connecting $\gamma^{\varepsilon}_{t_j} \circ \varphi_j$ 
with $\gamma^{\varepsilon}_{\infty}$ in every 
$C^k(\sphere^1,\rel^3)$, which can be centrally and 
bijectively projected from $\rel^3$ into $\sphere^2$ 
for each $j\geq J$: 
$$ 
L_j(x,\tau):=
\frac{(1-\tau) \, \gamma^{\varepsilon}_{t_j} \circ \varphi_j(x) + \tau \, \gamma^{\varepsilon}_{\infty}(x)}
{|(1-\tau) \,\gamma^{\varepsilon}_{t_j}\circ \varphi_j(x) 
+ \tau \, \gamma^{\varepsilon}_{\infty}(x)|},
\qquad \textnormal{for} \,x\in \sphere^1, \,\, 
\tau \in [0,1].  
$$  
Since we know by \eqref{smooth.subconvergence.1} that 
$\parallel \gamma^{\varepsilon}_{t_j} \circ \varphi_j
-\gamma^{\varepsilon}_{\infty}
\parallel_{C^2(\sphere^1,\rel^3)}$ 
tends to zero as $j\to \infty$, whence the path  
$L_j(\,\cdot \,,\tau):\sphere^1 \longrightarrow \sphere^2$ 
is indeed a closed smooth embedding into $\sphere^2$ 
for \underline{every} $\tau \in [0,1]$, provided $j \geq J$ with $J$ chosen sufficiently large, just as asserted.     
Hence, combining this particular isotopy 
with Theorem 1 in \cite{Arnold.1993}, Whitney's index-theorem and with equation \eqref{invariants.gamma.epsilon.T}, we finally conclude in cases (A) and (B): 
$$
\textnormal{St}_{\sphere^2}(\gamma^{\varepsilon}_{\infty})
= \textnormal{St}_{\sphere^2}(\gamma^{\varepsilon}_{t_j})  
\equiv 0 \, \textnormal{mod} \, \ganz,  
$$ 
for every $j\geq J$, contradicting our first result for $\textnormal{St}_{\sphere^2}(\gamma^{\varepsilon}_{\infty})
\, \textnormal{mod} \, \ganz$ below equation \eqref{invariants.gamma.epsilon.0}.
Hence, ``$\gamma^{\varepsilon}_{\infty}=E \oplus E$'' or 
``$\gamma^{\varepsilon}_{\infty}=E \oplus E \oplus E$''
- up to some isometry of $\sphere^2$ - are the only 
remaining possibilities. 
In order to derive a contradiction in these two cases 
as well, we will exploit the advantage of 
our second alternative to rely on the invariance of the \emph{tangent-rotation number} 
$\textnormal{ind}:\textnormal{Imm}(\sphere^1,\rel^2)
\longrightarrow \ganz$ with respect to regular homotopy only, thus avoiding \emph{further technical difficulties} arising from \emph{genericity} of closed curves.\\ 
Let's start with case (C):  
$\gamma^{\varepsilon}_{\infty}=E \oplus E$.    
Hence, again we choose some arbitrarily large index 
$J \in \nat$, set $T:=t_{j}$ for any $j \geq J$ and exploit 
the important fact that there is some 
$b_T \in \sphere^2\setminus \bigcup_{0\leq t\leq T} \textnormal{trace}(\gamma^{\varepsilon}_t)$
such that the projection  
$\Psi_{b_T}\circ \gamma^{\varepsilon}_T$ is an element 
of $\textnormal{Imm}(\sphere^1,\rel^2)$ in particular. 
Now, as explained in Lemma 3.9 of \cite{Dall.Acqua.Schaetzle.Mueller.2024}
the tangent-rotation number 
$\textnormal{ind}:\textnormal{Imm}(\sphere^1,\rel^2)
\longrightarrow \ganz$ can also be expressed analytically 
- as a particular path-integral - and is 
therefore not only the unique homotopy invariant on  
$\textnormal{Imm}(\sphere^1,\rel^2)$ but also  
a contiuous function from it into $\ganz$, 
if we consider here $\textnormal{Imm}(\sphere^1,\rel^2)$ 
as a topological subspace of $C^{2}(\sphere^1,\rel^2)$.  
Now, either we have the convenient situation in which 
$b_T \not \in \textnormal{trace}(E \oplus E)$ or the 
inconvenient situation in which the point $b_T$ 
sits exactly on the equator $\textnormal{trace}(E)$, 
for the fixed large time $T=t_j$. 
In the first subcase we will obviously have 
$\textnormal{ind}(\Psi_{b_T}\circ \gamma^{\varepsilon}_T)
=\textnormal{ind}(A_2)= \pm 2$, 
simply combining the assumed $C^k$-convergence \eqref{smooth.subconvergence.1} with the continuity of 
$\textnormal{ind}:(\textnormal{Imm}(\sphere^1,\rel^2), 
\parallel\, \cdot \,\parallel_{C^2(\sphere^1,\rel^2)}) 
\longrightarrow \ganz$.
But in the second subcase it is not that obvious, which 
index the projected curve 
$\Psi_{b_T}\circ \gamma^{\varepsilon}_T$ might have, 
independently of the size of the chosen $J$ or $T$.  
Fortunately, $E \oplus E$ is still a topologically   
fairly primitive path, and thus one can easily infer 
from convergence \eqref{smooth.subconvergence.1} 
that there are only two different cases, up to 
regular homotopy between elements of 
$\textnormal{Imm}(\sphere^1,\rel^2)$: 
(a) the path $\Psi_{b_T}\circ \gamma^{\varepsilon}_T$
is again regularly homotopic to the double loop 
$A_2$ respectively $K_2$ and has therefore rotation-index 
$\pm 2$, or (b) the path $\Psi_{b_T}\circ \gamma^{\varepsilon}_T$ is regularly homotopic to $K_0$, 
a ``large figure eight''
\footnote{See here again the explanations on p. 6 in 
\cite{Arnold.1993} or on p. 993 in \cite{Polyak.1996} 
and Figure 56 in \cite{Arnold.1993}.}, 
and has therefore  
$\textnormal{ind}(\Psi_{b_T}\circ \gamma^{\varepsilon}_T)
=0$, provided $T=t_j$ was chosen sufficiently large. 
Case (D) in which $\gamma^{\varepsilon}_{\infty}=
E \oplus E \oplus E$ can be handled similarly. 
Again, we either have the convenient situation in which 
$b_T \not \in \textnormal{trace}(E \oplus E \oplus E)$ 
or the inconvenient situation in which the point $b_T$ 
sits exactly on the equator $\textnormal{trace}(E)$, 
for the chosen time $T=t_j$. 
In the first subcase we will obviously have 
$\textnormal{ind}(\Psi_{b_T}\circ \gamma^{\varepsilon}_T)
=\textnormal{ind}(A_3)= \pm 3$, 
similarly to the above reasoning in the easy subcase  
of case (C). In the second subcase of case (D) 
it is again more difficult to guess which index the 
projected curve $\Psi_{b_T}\circ \gamma^{\varepsilon}_T$ 
might have, for any large $T=t_j$.  
Fortunately, in this concrete situation one can still 
infer by means of Figure 56 in \cite{Arnold.1993} 
that convergence \eqref{smooth.subconvergence.1}  
can again only result in two qualitatively different cases, 
at least up to regular homotopy: 
(a) the path $\Psi_{b_T}\circ \gamma^{\varepsilon}_T$
is again regularly homotopic to the triple loop 
$A_3$ and has therefore rotation-index $\pm 3$, 
or (b) the path $\Psi_{b_T}\circ \gamma^{\varepsilon}_T$ is regularly homotopic to another homotopy-type of 
closed generic paths in the plane, 
namely to a closed generic path having two double points 
- just as $A_3$ - but rotation-index $\pm 1$ instead of 
$\pm 3$, provided $J$ was chosen sufficiently large; 
see here especially the first row of Figure 56 and 
Figure 28 in \cite{Arnold.1993}. Now, the most important consequence of the possibility to choose some  
$b_T \in \sphere^2\setminus \bigcup_{0\leq t\leq T} \textnormal{trace}(\gamma^{\varepsilon}_t)$
is that 
$\{\Psi_{b_T}\circ \gamma^{\varepsilon}_t\}_{t\in [0,T]}$ 
is a well-defined regular homotopy between 
the regular paths $\Psi_{b_T}\circ \gamma^{\varepsilon}_0$ and 
$\Psi_{b_T}\circ \gamma^{\varepsilon}_T$ in $\rel^2$.  
We can therefore use the invariance of the tangent-rotation number with respect to \emph{regular homotopy on entire} $\textnormal{Imm}(\sphere^1,\rel^2)$ together 
with the obvious fact that 
$$
\textnormal{ind}\big{(}\Psi_{b_T} \circ 
\gamma^{\varepsilon}_0 \big{)} = 
\textnormal{ind}\big{(}\Psi_{b_T} \circ 
(\exp_{\gamma^*}(\varepsilon \vec F_{\gamma^*}))\big{)}  
= \pm 4,
$$
as already pointed out above in equation \eqref{index.gamma.0.a}, in order to arrive at 
$$
\pm 4 =\textnormal{ind}(\Psi_{b_T} \circ \gamma^{\varepsilon}_0)
=\textnormal{ind}(\Psi_{b_T} \circ \gamma^{\varepsilon}_T) 
\in \{0,\pm 2\} \cup \{\pm 1,\pm 3\}=\{0,\pm 1,\pm 2,\pm 3\},
$$ 
which is again wrong, and Theorem 
\ref{Singularities.do.exist} is proved.  
\begin{remark}  \label{no.classical.approach}
We mentioned in the introduction that 
on the one hand Hamilton's and Huisken's \cite{Huisken.1984} early papers about the Ricci- and mean-curvature flow    
started to stress the general picture that 
one should either aim at optimal breakdown-criteria and lifespan-estimates for some challenging geometric flow, 
as e.g. the Inverse Willmore flow \cite{Mayer}, 
or prove global existence and subconvergence to some 
smooth limit immersion of every flow line, as e.g. in  \cite{Dall.Acqua.Pozzi.2014}, \cite{Dall.Acqua.Pozzi.2018} 
dealing with the classical elastic energy flow \eqref{classical.elastic.energy.flow}, but that on the other hand {\bf this entire methodology apparently fails} when applied to the MIWF. Here we shall finally prove this assertion by means of the restriction of the MIWF to smooth Hopf-tori in $\sphere^3$ which reduces to our degenerate variant \eqref{elastic.energy.flow} of the classical elastic energy flow \eqref{classical.elastic.energy.flow} in $\sphere^2$ via the Hopf-fibration; compare here with Proposition \ref{correspond.flows}. More concretely, 
in this remark we will 
{\bf prove the fundamental assertion} that our final results, 
Theorem \ref{limit.at.infinity} and Corollary \ref{limit.at.infinity.1}, {\bf cannot be achieved} 
by means of any adaption of Hamilton's and Huisken's  classical methods, being applied to the degenerate elastic energy flow \eqref{elastic.energy.flow} or to the MIWF directly - certainly a surprising key-insight about the MIWF and its simplified, subcritical variant \eqref{elastic.energy.flow}. 
Here, this approach would be to combine formulae (2.14), (2.18), (2.21) and (2.23) in Lemma 2.3 of \cite{Dall.Acqua.Pozzi.2018} with our flow equation \eqref{elastic.energy.flow}, in order to 
obtain - as a first step - explicit evolution equations for the curvature vector $\vec \kappa_{\gamma_t}$ and 
for the arclength $d\mu_{\gamma_t}$ along a general flow line $\{\gamma_t\}_{t\in [0,T]}$ of the degenerate elastic energy flow \eqref{elastic.energy.flow} and to combine them with the following trick, the analog of Proposition 4.2 in \cite{Mayer}. 
	\begin{proposition} \label{evolution.nabla.Phi}
	Let $\{\gamma_t\}_{t\in [0,T]}$ be a smooth flow line of evolution equation \eqref{elastic.energy.flow}, and let $\{\Phi_t\}_{t\in [0,T]}$ 
	be a family of smooth normal vector fields along 
	$\{\gamma_t\}_{t\in [0,T]}$, satisfying:\\ 
	$\nabla_t^{\perp}(\Phi_t) + 
	\frac{2}{(\kappa_{\gamma_t}^2+1)^2} \, 
	(\nabla^{\perp}_s)^4(\Phi_t) = Y_t$
	for some smooth normal vector field $Y_t$ along 
	$\{\gamma_t\}_{t\in [0,T]}$. Then its covariant derivative $\Psi_t:=\nabla^{\perp}_s\Phi_t$ satisfies the equation
	\begin{eqnarray*}  
			\nabla_t^{\perp}(\Psi_t) + 
			\frac{2}{(\kappa_{\gamma_t}^2 + 1)^2} \, \big{(}\nabla^{\perp}_s\big{)}^4(\Psi_t) = 
			\nabla^{\perp}_s(Y_t) 
			+ \frac{8}{(\kappa_{\gamma_t}^2 + 1)^3} \, 
			\langle \vec \kappa_{\gamma_t}, 
			\nabla_s^{\perp}(\vec \kappa_{\gamma_t}) \rangle \, \big{(}\nabla^{\perp}_s\big{)}^4(\Phi_t)    \\
			-\Big{[} \frac{2}{(\kappa_{\gamma_t}^2+1)^2} \, 
			\langle \big{(}\nabla_s^{\perp}\big{)}^2(\vec \kappa_{\gamma_t}),
			\vec \kappa_{\gamma_t} \rangle 
			+\frac{\kappa_{\gamma_t}^2}{1+\kappa_{\gamma_t}^2} \Big{]} \cdot 
			\nabla_s^{\perp}(\Phi_t)                  \\  
			+\Big{[} \frac{8}{(\kappa_{\gamma_t}^2+1)^3} \, 
			\langle \vec \kappa_{\gamma_t}, 
			\nabla_s^{\perp}(\vec \kappa_{\gamma_t}) \rangle\, 
			\big{(}\nabla_s^{\perp}\big{)}^2(\vec \kappa_{\gamma_t}) 
			- \frac{2}{(\kappa_{\gamma_t}^2+1)^2} \, 
			\big{(}\nabla_s^{\perp}\big{)}^3(\vec \kappa_{\gamma_t})
			- \frac{1}{1+\kappa_{\gamma_t}^2}\,
			\nabla_s^{\perp}(\vec \kappa_{\gamma_t}) \Big{]}\,
			\langle \vec \kappa_{\gamma_t}, \Phi_t \rangle   \\                        
			+\frac{2\,\vec \kappa_{\gamma_t}}{(1+\kappa_{\gamma_t}^2)^2} \,  
			\big{\langle} \big{(} \nabla^{\perp}_{s} \big{)}^3 
			(\vec \kappa_{\gamma_t}),\Phi_t \big{\rangle}               
			- \frac{8 \,\vec \kappa_{\gamma_t}}{(1+\kappa_{\gamma_t}^2)^3} \,\langle \vec \kappa_{\gamma_t}, 
			\nabla_s^{\perp}(\vec \kappa_{\gamma_t}) \rangle\,  
			\langle \big{(} \nabla^{\perp}_{s} \big{)}^2 
			(\vec \kappa_{\gamma_t}),\Phi_t \rangle  \\
			+ \frac{\vec    \kappa_{\gamma_t}}{1+\kappa_{\gamma_t}^2} \,
			\langle \nabla^{\perp}_{s}(\vec \kappa_{\gamma_t}),\Phi_t 
			\big{\rangle}
		\end{eqnarray*} 
		for every $t\in [0,T]$, where we abbreviated above  
		$\nabla^{\perp}_s:=
		\nabla^{\perp}_{\frac{\gamma'}{|\gamma'|}}$ 
		for ease of notation. 
	\end{proposition}  
	\noindent 
	Now, following Sections 3 and 4 in \cite{Dall.Acqua.Pozzi.2014} and Sections 2.2 and 4.1 in \cite{Dall.Acqua.Pozzi.2018} and also in view of the uniform bounds \eqref{Length.total} and \eqref{bound.below} we introduce some technically useful notations.   
	\begin{definition}   \label{P.a.b.c}	
	For a fixed smooth, closed and regular curve $\gamma:\sphere^1 \longrightarrow \sphere^2$, and integers $b \geq 2$ and $a \geq 0$, $c\geq 0$,  
	we call ``$P^{a,c}_b(\vec \kappa_{\gamma})$'' any finite 
	linear combination of products 
	$$ 
	(\nabla^{\perp}_s)^{i_1}(\vec \kappa_{\gamma}) \ast 
	\ldots \ast (\nabla^{\perp}_s)^{i_b}(\vec \kappa_{\gamma})	
	$$
	with $i_1 + \ldots + i_b = a$ and $\max i_j \leq c$, 
	where we abbreviated again $\nabla^{\perp}_s
	:=\nabla^{\perp}_{\frac{\gamma'}{|\gamma'|}}$.  
	\end{definition} 
	\noindent 
	Combining now the mentioned evolution equations for  
	$\vec \kappa_{\gamma_t}$ and $d\mu_{\gamma_t}$ with Proposition \ref{evolution.nabla.Phi}, we obtain the 
	following general evolution equation by induction, which 
	corresponds to Proposition 4.3 in \cite{Mayer} or also 
	to Lemma 3.1 in \cite{Dall.Acqua.Pozzi.2014}. 
	\begin{proposition}  \label{evolution.nabla.k}
	Let $\{\gamma_t\}_{t\in [0,T]}$ be a flow line of evolution equation \eqref{elastic.energy.flow} and 
	$k \in \nat_0$. Then the vector field 
	$(\nabla^{\perp}_s)^k(\vec \kappa_{\gamma_t})$ 
	satisfies the equation  
	\begin{eqnarray}   \label{evolution.nabla.m}
	\nabla_t^{\perp}((\nabla^{\perp}_s)^k(\vec \kappa_{\gamma_t})) + \frac{2}{(\kappa_{\gamma_t}^2+1)^2} \, (\nabla^{\perp}_s)^4((\nabla^{\perp}_s)^k
	   (\vec \kappa_{\gamma_t}))  \nonumber \\
		= \frac{1}{(\kappa_{\gamma_t}^2+1)^3} 
		P^{k+4,k+3}_3(\vec \kappa_{\gamma_t})
		+\sum_{(a,b,d) \in I(k)} 
		\frac{1}{(\kappa_{\gamma_t}^2+1)^d}\,
		P^{a,k+2}_b(\vec \kappa_{\gamma_t}),		
		\end{eqnarray}
		for any $t\in [0,T]$, where we used the notation of Definition \ref{P.a.b.c} and where the set  
		$I(k)$ consists of those triples $(a,b,d) \in \nat_0^3$, such that $a \in \{k+4,k+2,k\}$, $b$ is odd with $b\leq 2k+5$, and $1\leq d \leq k+4$.
	\end{proposition}  
	\noindent 
	Interestingly, trying to infer  
    $L^{\infty}$-$L^{\infty}$-estimates for   
	covariant derivatives $(\nabla^{\perp}_s)^k
	(\vec \kappa_{\gamma_t})$ of only low orders 
	$k=0,1,2,3,\ldots$ from a combination of Proposition \ref{evolution.nabla.k} with Lemma 4.3 in \cite{Dall.Acqua.Pozzi.2014} - a tricky interpolation inequality - one runs into 
	\emph{severe computational difficulties}. 
	More precisely, the usual trick ``to multiply equation \eqref{evolution.nabla.m}
	with $(\nabla^{\perp}_s)^k(\vec \kappa_{\gamma_t})$ and 
	then to integrate by parts twice'' does not elegantly work 
	out here on account of the factor ``$\frac{2}{(\kappa_{\gamma_t}^2+1)^2}$'' on 
	the left hand side of equation \eqref{evolution.nabla.m}, 
	and moreover because of - sloppily speaking - too many terms of relatively high order on the right hand side of equation \eqref{evolution.nabla.m}, 
	preventing us from a successful interpolation 
	via Lemma 4.3 in \cite{Dall.Acqua.Pozzi.2014}. \\ 
    {\bf And indeed}, this direct PDE-argument necessarily 
    has to fail here ! 
	In order to prove this, we recall that a successful application of the above mentioned standard technique would yield a-priori $L^{\infty}$-$L^{\infty}$-estimates for the curvature vector $\vec \kappa_{\gamma_t}$ along any given flow line $\{\gamma_t\}_{t\in [0,T_{\textnormal{max}})}$ of \eqref{elastic.energy.flow} together with all its covariant derivatives 
	$(\nabla^{\perp}_s)^k(\vec \kappa_{\gamma_t})$ and 
	$(\nabla_s)^k(\vec \kappa_{\gamma_t})$ up to the 
	maximal time of existence ``$T_{\textnormal{max}}$'', 
	which is assumed to be finite first of all: 
	\begin{equation}  \label{curvature.bound.infty}
	\parallel \big{(} \nabla_{\frac{\gamma_t'}{|\gamma_t'|}}\Big{)}^k
	(\vec{\kappa}_{\gamma_t}) \parallel_{L^{\infty}(\sphere^1)} \leq C, \quad \textnormal{for} \,\,\, t \in [0,T_{\textnormal{max}}),
	\end{equation}
	for some sufficiently large constant $C=C(\gamma_0,k,T_{\textnormal{max}})$, 
	for each $k\in \nat_0$, exactly as explained
	in steps 1-6 of the proof of Theorem 1.1 in \cite{Dall.Acqua.Pozzi.2014} respectively on  
	page 12 in \cite{Dall.Acqua.Pozzi.2018}. Now, adapting 
	the arguments of step no. 7 within the proof of Theorem 1.1 in \cite{Dall.Acqua.Pozzi.2014} slightly to our 
	degenerate flow equation \eqref{elastic.energy.flow} 
	we could improve a-priori estimates 
	\eqref{curvature.bound.infty} - essentially following 
	also here Dall'Acqua's and Pozzi's induction argument  combining formulae (5.17) and (5.18) in \cite{Dall.Acqua.Pozzi.2014} with estimates \eqref{curvature.bound.infty} and Gronwall's Lemma: 
	\begin{equation}  \label{curvature.bound.infty.2}
	\parallel \Big{(} \nabla_{\gamma_t'}\Big{)}^k
	(\vec{\kappa}_{\gamma_t}) \parallel_{L^{\infty}(\sphere^1)} \leq \tilde C \quad \textnormal{for} \,\,\, t \in [0,T_{\textnormal{max}}),
	\end{equation}
	for another, larger constant 
	$\tilde C=\tilde C(\gamma_0,k,T_{\textnormal{max}})$. 
	Now, as in step no. 8 within the proof of Theorem 1.1 in \cite{Dall.Acqua.Pozzi.2014} we could use estimates   \eqref{curvature.bound.infty.2} in order to extend the 
	considered general flow line  
	$\{\gamma_t\}_{t\in [0,T_{\textnormal{max}})}$   	   
	from $[0,T_{\textnormal{max}})$ to some longer 
	time interval, thus obtaining a contradiction if 
	$T_{\textnormal{max}}$ was actually finite. 
	Hence, we could infer here that any given flow line of \eqref{elastic.energy.flow} 
	with smooth initial path $\gamma_0$ was certainly a global flow line $\{\gamma_t\}_{t\in [0,\infty)}$.  
	Hence, as in step no. 9 within the proof of Theorem 
	1.1 in \cite{Dall.Acqua.Pozzi.2014} we would automatically arrive here at the uniform $L^{\infty}$-estimates 
	(5.22), i.e. at a-priori estimates \eqref{curvature.bound.infty} 
	on the entire ray $[0,\infty)$ - just as in formula (63) of \cite{Ruben.MIWF.II} for flow lines of the classical 
	elastic energy flow \eqref{classical.elastic.energy.flow} 
	on $\sphere^2$ - and exactly as in \cite{Dall.Acqua.Pozzi.2014} or in \cite{Ruben.MIWF.II} we could go on from here and conclude by means of estimates
	\eqref{Length.total} and \eqref{bound.below} -  
	respectively (60) and (62) in \cite{Ruben.MIWF.II} -    
	on the length of every path $\gamma_t$ along the considered flow line that for some appropriate sequence $\{t_j\}$ with $t_j \to \infty$ at least the reparametrizations $\tilde \gamma_{t_j}$ of the paths $\gamma_{t_j}$ ``to their arc-length'' would 
	{\bf converge smoothly}: 
	\begin{equation} \label{smooth.subconvergence} 
	\tilde \gamma_{t_j} := \gamma_{t_j} \circ \varphi_j 
	\longrightarrow \gamma_{\infty} \quad  \textnormal{in}	\,\,\,C^k(\sphere^1,\rel^3),\,\,\, \forall \, k \in \nat, 
	\end{equation}  
	as $j \to \infty$, for some 
	non-constant, smooth closed path $\gamma_{\infty}$ with $|\gamma_{\infty}'|\equiv \textnormal{const}>0$. 
	Moreover, up to minor modifications we could infer now as in \cite{Dall.Acqua.Pozzi.2014} or as in \cite{Ruben.MIWF.II}, pp. 23--24, from convergence \eqref{smooth.subconvergence} and estimates \eqref{curvature.bound.infty} on the entire 
	ray $[0,\infty)$ that the smooth limit curve $\gamma_{\infty}$ in \eqref{smooth.subconvergence}  
	would have to be a {\bf closed elastic curve} in $\sphere^2$.
    However, comparing this result with the statement of Theorem \ref{Singularities.do.exist} above we instantly obtain a contradiction, because we considered here some arbitrary flow line of flow \eqref{elastic.energy.flow} emanating from some smooth, regular closed initial path $\gamma_0$ in $\sphere^2$, without any further conditions. 
\end{remark}

\section{Appendix A}  \label{remarks}
\begin{remark}  \label{no.better} 
We should remark here first of all
that the $W^{4,2}$-regularity which we achieved 
for the conformal parametrization 
$f:(\Sigma,g_{\textnormal{poin}}) 
\stackrel{\cong}\longrightarrow \textnormal{spt}(\mu)$
of a limit torus $\textnormal{spt}(\mu)$ in the fourth part of Theorem \ref{limit.MIWF} and in the third part of Theorem \ref{singular.time.MIWF.Hopf.tori},
is even stronger than the minimally required regularity of  
some umbilic-free initial immersion of the torus $\Sigma$ into 
any $\rel^n$, $n\geq 3$, in which the ``relaxed MIWF-equation'' \eqref{generalized.MIWF} can be uniquely started; see Definition \ref{relaxed.MIWF} below and the short-time existence result in Theorem 2 (1) of \cite{Ruben.MIWF.III}. 
But unfortunately, the conformal parametrization $f$ of 
the limit torus $\textnormal{spt}(\mu)$ in the fourth 
part of Theorem \ref{limit.MIWF} might fail to be umbilic-free, and one obviously cannot 
restart the MIWF in a limit-immersion $f$ possessing 
some umbilic points - at least not within our classical 
framework for the MIWF.   
Yet no single technique could be developed yielding a general criterion which rules out that for some given sequence 
$t_j \nearrow T_{\textnormal{max}}$ there might be 
points $x_j\in \Sigma$ such that    $|A^0_{F_{t_j}}(x_j)|^2\longrightarrow 0$ as $j \to \infty$. 
Only in the special situation of the third part of Theorem \ref{singular.time.MIWF.Hopf.tori} 
we trivially know that the limit Hopf-torus 
$\textnormal{spt}(\mu)$ of some arbitrarily chosen 
singular flow line $\{F_t\}_{t\in [0,T_{\textnormal{max}}(F_0))}$ is umbilic-free, but still it is unclear how to use estimate \eqref{L.infty.W.4.2.estimate} 
and the flow equations \eqref{Moebius.flow} and 
\eqref{elastic.energy.flow}, in order to    
prove that a suitable reparametrization of the 
considered singular flow line 
$\{F_t\}_{t\in [0,T_{\textnormal{max}}(F_0))}$ 
would have to extend to a function of class  
$W^{1,p}([0,T_{\textnormal{max}}(F_0)];L^p(\Sigma,\rel^4)) \cap L^p([0,T_{\textnormal{max}}(F_0)];W^{4,p}(\Sigma,\rel^4))$- say with $p=2$ - such that $\{F_t\}$ would automatically have a well-defined trace $F_T$ in $W^{2,2}(\Sigma,\rel^4)$ at $t=T$. Since the above parabolic $L^p$-space embeds into 
$C^0([0,T_{\textnormal{max}}(F_0)];
W^{4-\frac{4}{p},p}(\Sigma,\rel^4))$, for any $p>1$, 
$F_T$ would have to parametrize the same manifold 
as the uniformly conformal $W^{4,2}$-parametrization $f$ of the limit Hopf-torus $\textnormal{spt}(\mu)$ from the third part of Theorem \ref{singular.time.MIWF.Hopf.tori}.  
Hence, in this situation it would make sense to 
``restart'' the relaxed MIWF \eqref{generalized.MIWF} 
in the umbilic-free $W^{4,2}$-parametrization 
$f:\Sigma \stackrel{\cong}\longrightarrow \textnormal{spt}(\mu)$ by means of Theorem 2 (1) in \cite{Ruben.MIWF.III}; compare here also 
with Remark \ref{no.better.3} below.   
\end{remark}  
\begin{definition} \label{relaxed.MIWF}
	\footnote{The choice of parabolic $L^p$-spaces 
	$X_{T,p}:=W^{1,p}([0,T],L^{p}(\Sigma,\rel^n)) 
	\cap L^{p}([0,T],W^{4,p}(\Sigma,\rel^n))$ in Definition 
	\ref{relaxed.MIWF} with $p>3$ is motivated 
	by the techniques and results of the author's paper  \cite{Ruben.MIWF.III}. In that paper the author 
	investigated the DeTurck-modification of the MIWF - 
	a quasilinear evolution equation of 4th order -
	proving well-posedness of this modified flow, 
	maximal regularity of its strong solutions  
	in $X_{T,p}$ and several properties of its 
	evolution operator mapping 
	$\textnormal{trace}(X_{T,p})=
	W^{4-\frac{4}{p},p}(\Sigma,\rel^n)$ into $X_{T,p}$. 
	In view of modern optimal regularity and semigroup theory, 
	the quasilinear character of the differential operator 
	on the right hand side of the modified MIWF-equation 
	can be employed effectively, if and only if 
	$\textnormal{trace}(X_{T,p})=
	W^{4-\frac{4}{p},p}(\Sigma,\rel^n)$ embeds into  
	$C^{2,\alpha}(\Sigma,\rel^n)$ for some 
	positive $\alpha \leq 2-\frac{6}{p}$, which holds 
	only for $p>3$; see here especially Lemma 1 and  
	Theorems 1, 2 and 4 in \cite{Ruben.MIWF.III}.   
	This rather orthodox approach to the 
	MIWF does not claim to yield any type of ``weak solutions'' of the MIWF under minimal regularity assumptions on their initial immersions, and it should not be confused with Palmurella's and Riviere's ambitious attempt in \cite{Palmurella.2022} and \cite{Palmurella.2024} to construct solutions of some particular weak formulation 
	of the classical Willmore flow equation which are 
	one-parameter families of uniformally conformal 
	$W^{1,\infty}(\Sigma)\cap W^{2,2}(\Sigma)$-immersions.}   
	We call a family $\{F_t\}_{t\in [0,T]}$ of immersions of $\Sigma$ into $\rel^n$ a \emph{relaxed flow line} of the MIWF in $\rel^n$, if the resulting function 
	$F:\Sigma \times [0,T]\longrightarrow \rel^n$ 
	is of class $W^{1,p}([0,T],L^{p}(\Sigma,\rel^n)) 
	\cap L^{p}([0,T],W^{4,p}(\Sigma,\rel^n))$, 
	for some $T\in (0,\infty)$ and $p\in (3,\infty)$,  
	and satisfies the ``relaxed MIWF-equation'': 
	\begin{equation}   \label{generalized.MIWF}
		\big{(}\partial_t F_t(x)\big{)}^{\perp_{F_t}} 
		= -\frac{1}{2} \frac{1}{|A^0_{F_t}(x)|^4} \,
		\Big{(} \triangle_{F_t}^{\perp} \vec H_{F_t}(x) + Q(A^{0}_{F_t})(\vec H_{F_t})(x) \Big{)}
	\end{equation} 
	in a.e. $(x,t)\in \Sigma \times [0,T]$.
	As in equation \eqref{second.fundam.form}, the symbol 
	``$^{\perp_{F_t}}$'' abbreviates the orthogonal projection of the velocity vector $\partial_t F_t(x)$ into the normal 
	space of the immersion $F_t$ within $\rel^n$, at any 
	fixed $x\in \Sigma$.
\end{definition}  
\begin{remark}  \label{no.better.2}   
A second technical issue is here the possibility 
that the limit Hopf-torus $\textnormal{spt}(\mu)$, 
which we had obtained in the third 
part of Theorem \ref{singular.time.MIWF.Hopf.tori},  
might depend on the choice of the particular subsequence 
in (\ref{weak.convergence.mu})--(\ref{Hausdorff.converg})
of the originally considered sequence of embeddings 
$\{F_{t_{j}}\}$ - for any fixed sequence of times 
$t_j \nearrow T_{\textnormal{max}}(F_0)$. 
This is another reason why neither the fourth part of Theorem \ref{limit.MIWF} nor the third part of Theorem \ref{singular.time.MIWF.Hopf.tori} of this article 
can be combined with Theorems 2 (1) and 3 (1) 
in \cite{Ruben.MIWF.III}, in order to yield some necessary condition for a flow line to be singular, i.e. to break down in finite time $T_{\textnormal{max}}(F_0)<\infty$. 
We must therefore confess that 
both the fourth part of Theorem \ref{limit.MIWF} and 
the third part of Theorem \ref{singular.time.MIWF.Hopf.tori}
slightly miss their original aim to rule out 
the development of singularities of the MIWF 
at some finite maximal existence time $T_{\textnormal{max}}(F_0)$.  
More precisely, we admit here that we can neither 
conclude from the fourth part of Theorem \ref{limit.MIWF} 
nor from the third part of Theorem \ref{singular.time.MIWF.Hopf.tori} that  
either the supremum of the mean curvature, i.e. 
$\parallel \vec H_{F_{t}} \parallel_{L^{\infty}(\Sigma)}$, 
or the supremum of $|A^0_{F_{t}}|^2$ over $\Sigma$,   
or the speed of Willmore-energy-decrease, i.e.  $|\frac{d}{dt}\Will(F_{t})|$, have to ``blow up'' along 
a general singular flow line $\{F_{t}\}$ of the MIWF, 
as $t$ approaches the singular time $T_{\textnormal{max}}(F_0)<\infty$ from the past. \\ 
Moreover, for the same reason the third part of Theorem \ref{limit.MIWF} and also the third part of Theorem \ref{singular.time.MIWF.Hopf.tori} 
show us on account of statements \eqref{no.A.concentration}
and \eqref{no.concentration.Mill},
that the phenomenon of {\bf concentration of curvature} in the ambient space $\rel^4$ of embeddings $\{F_{t}\}$ 
moving along the MIWF in $\sphere^3$ is 
probably {\bf not a criterion for the respective 
flow line $\{F_{t}\}$ to develop a singularity} in 
finite time, in contrast to the famous statement of Theorem 1.2 in \cite{Kuwert_Schaetzle_2002} 
about the classical Willmore flow in $\rel^n$, $n\geq 3$, 
and even the initial energy threshold $\Will(F_0)\leq 8\pi$ 
does not improve this picture.  
Hence, unfortunately Theorems \ref{limit.MIWF} and 
\ref{singular.time.MIWF.Hopf.tori} do not support 
the expectation that we might figure out 
optimal geometric criteria for a flow line of the MIWF in $\sphere^3$ to develop a singularity in finite time - not even in a relatively low energy regime - indicating a stark contrast to the behaviour of the classical Willmore flow in $\rel^n$; compare here with Theorem 1.2 in \cite{Kuwert_Schaetzle_2002} and Theorem 5.2 in \cite{Kuwert.Schaetzle.Annals}.   
\end{remark}  
\begin{remark}  \label{no.better.3} 
As already mentioned in Remark \ref{no.better.2} above, 
we only consider sequences 
$\{F_{t_j}\}$ with $t_j \nearrow T_{\textnormal{max}}(F_0)$
in Theorems \ref{limit.MIWF} and \ref{singular.time.MIWF.Hopf.tori},
which we reparametrize in such a way that 
each immersion $\tilde F_{t_j}:=F_{t_j} \circ \Phi_j$ becomes 
uniformly conformal with respect to some smooth metric of zero scalar curvature $g_{\textnormal{poin},j}$ and such
that these metrics $g_{\textnormal{poin},j}$
- up to extraction of a subsequence -  
converge smoothly to some fixed smooth metric $g_{\textnormal{poin}}$ of zero scalar curvature on 
the abstract torus $\Sigma$. Instead we could  
pragmatically choose $\Sigma$ to be the 
\underline{standard Clifford torus} in $\sphere^3$, 
and we could try to reparametrize every single 
embedding $F_{t}$ of the flow line $\{F_{t}\}$, i.e. 
for every single $t\in [0,T_{\textnormal{max}}(F_0))$, 
in such a way that the new, reparametrized family of embeddings $\{\tilde F_{t}\}_{t \in [0,T_{\textnormal{max}}(F_0))}$ 
still solves the ``relaxed MIWF-equation'' \eqref{generalized.MIWF}, which is a quasi-linear 
parabolic differential equation. 
As in Section 5 of \cite{Palmurella.2022} the motivation 
for this idea is obvious: 
We try to use a systematical and continuous gauge of 
all metrics $\{F_t^*(g_{\textnormal{euc}})\}_{t\in [0,T_{\textnormal{max}}(F_0))}$ into uniformly conformal metrics  
$$
\{\tilde F_t^*(g_{\textnormal{euc}})\}_{t\in [0,T_{\textnormal{max}}(F_0))}
=\{e^{2u_t}\, g_{\textnormal{poin}}\}_{t\in [0,T_{\textnormal{max}}(F_0))},
$$
where $g_{\textnormal{poin}}$ should be some fixed smooth metric of zero scalar curvature on the Clifford torus, 
in order to prove - under appropriate, mild conditions
on the considered flow line $\{F_t\}_{t\in [0,T_{\textnormal{max}}(F_0))}$ - by means of a parabolic bootstrap argument $C^{\infty}$-smoothness of 
the ``relaxed flow line'' $\{\tilde F_t\}$ up to $t=T_{\textnormal{max}}(F_0)$ and actually including $t=T_{\textnormal{max}}(F_0)$, provided
$T_{\textnormal{max}}(F_0)$ is finite. 
Following many classical examples in the field of geometric flows we could use here Theorem 2 (1) and Theorem 3 (1) in \cite{Ruben.MIWF.III} in order to obtain a contradiction to the finiteness of $T_{\textnormal{max}}(F_0)$, and there would actually have to hold $T_{\textnormal{max}}(F_0)=\infty$. Hence, this technique might rule out singularities along $\{F_t\}$, at least under appropriate additional conditions on the original flow line $\{F_t\}$ of the MIWF. 
However, comparing the ``distributional Willmore flow'' for 
immersions of $\sphere^2$ into $\rel^3$ in 
\cite{Palmurella.2022} with the MIWF \eqref{Moebius.flow}
we discover at least two big problems when trying
to adapt such a program to the MIWF: 
	\begin{itemize} 
	\item[1)] First of all, following Definition 1.8 in 
	\cite{Palmurella.2022} and thus trying to write 
	the MIWF-equation \eqref{Moebius.flow} ``in a weak form'', namely in terms of the 2nd order distribution 
	\eqref{distributional.Willmore}  
	instead of the well-known 4th order differential operator  
	\eqref{Euler.Lagrange.operator}, 
    does not work out, simply because the leading factor 
    $\frac{1}{|A^0_{F_t}|^4}$ on the right hand side of 
    \eqref{Moebius.flow} prevents us from successfully integrating by parts. 
	Hence, this entire approach to the MIWF seems to 
	break down from the very beginning. 	  		    
    But even if we could somehow overcome this first obstruction, then we would neither have formula (2.2) in Theorem 2.4 of \cite{Palmurella.2022} nor the estimates of Theorem 2.8 of \cite{Palmurella.2022} at our disposal, because the main results of M\"uller's and De Lellis' papers \cite{De.Lellis.Mueller.2005} and \cite{De.Lellis.Mueller.2006} and of Kuwert's and 
	Scheuer's article \cite{Kuwert.Scheuer.2020}   
	only hold for the classical Willmore flow moving  
	immersions of the $2$-sphere into $\rel^3$ respectively $\rel^n$, whereas we would have to consider immersions of a fixed torus - for example of the Clifford torus - into $\sphere^3$. We therefore do not have any tool which might help us to control either the deviation 
	$\parallel e^{u_t}-1 \parallel_{L^{\infty}(\Sigma)}$ 
	of the conformal factors $e^{u_t}$ of          
	$\tilde F_t^*(g_{\textnormal{euc}})$ on $\Sigma$ 
	or the deviation of areas and barycenters of
	the immersions $\tilde F_t$ in terms of the initial Willmore energy $\Will(F_0)$ and initial area 
	$\Area(F_0)$, nor in terms of any other controllable geometric quantity, for all $t\in [0,T_{\textnormal{max}}(F_0))$
	\footnote{See here the proofs of Theorem 1.9 and Proposition 5.3 in \cite{Palmurella.2022} for 
	precise information.}. 
	\item[2)] 
	Another fundamental problem lies in the fact that first of all the embeddings $F_t$ are uniformly conformal with respect to varying metrics $g_{\textnormal{poin}}(t)$ of zero scalar curvature, and that the moduli space $\Mod_1$ is isomorphic to $\quat/\textnormal{PSL}_2(\ganz)$
	by Theorem 2.7.2 in \cite{Jost.2006}, 
	whereas there is - up to conformal automorphisms - 
	only one conformal class on $\sphere^2$ 
	by Corollary 5.4.1 in \cite{Jost.2006}. 
	Therefore the reasoning yielding Lemma 1.1 in \cite{Palmurella.2022} breaks down, i.e. we cannot 
	easily achieve an explicit and useful  
	formula - just as formula (1.8) in \cite{Palmurella.2022} -
	for the tangential vector field $\{U(F_t)\}$ along $\{F_t\}$ which yields the concretely modified flow equation 
	$$ 
	\partial_t\tilde F_t \stackrel{!}= 
	- \frac{1}{|A^0_{\tilde F_t}|^4} \,\nabla_{L^2} \Will(\tilde F_t) + U(\tilde F_t) \quad \textnormal{on} \,\,\, \Sigma \times [0,T_{\textnormal{max}}(F_0))
	$$       
	whose unique smooth solution would be the 
	correctly gauged flow line $\{\tilde F_t\}$, with 
	$\tilde F_0=F_0$.
	This problem was actually addressed in the recent contribution \cite{Palmurella.2024}, however, 
	as we have observed above in the first part of this 
	remark, the basic idea to replace the original 
	Euler-Lagrange-operator \eqref{Euler.Lagrange.operator} by its distributional counterpart \eqref{distributional.Willmore} on the right hand 
	side of \eqref{Moebius.flow} does not work out here, 
	and we therefore cannot even start to investigate 
	the question whether the short-time existence result in Theorem 1 of \cite{Palmurella.2024} - with or without its lifespan-estimate - might hold for the MIWF, as well. 
	\end{itemize}
\end{remark}

\section{Appendix B}

In this appendix we recall and quote Lemma 3 in \cite{Ruben.MIWF.II}, i.e. the existence of horizontal  
$C^{\infty}$-lifts of some arbitrary closed path $\gamma:\sphere^1 \to \sphere^2$ of class $C^{\infty}$ 
with respect to fibrations of the type $\pi \circ F$, for 
smooth parametrizations $F:\Sigma \longrightarrow
\pi^{-1}(\textnormal{trace}(\gamma))$ of the Hopf-torus 
$\pi^{-1}(\textnormal{trace}(\gamma)) \subset \sphere^3$. 
The following lemma is an important tool in the proof of Propositions \ref{Hopf.Willmore.prop} and \ref{correspond.flows} of this article. 
\begin{lemma} \label{closed.lifts} [Lemma 3 in \cite{Ruben.MIWF.II}]
Let $\gamma: \sphere^1 \longrightarrow \sphere^2$ be a regular path of class $C^{\infty}$, and let 
$F:\Sigma \longrightarrow \sphere^3$ be a
smooth immersion, parametrizing the Hopf-torus
$\pi^{-1}(\textnormal{trace}(\gamma)) \subset \sphere^3$, 
see here Definition \ref{Hopf.torus.immersion}.
\begin{itemize}
\item[1)] For every fixed $s^* \in \sphere^1$
and $q^* \in \pi^{-1}(\gamma(s^*))\subset \sphere^3$
there is a unique horizontal smooth lift $\eta^{(s^*,q^*)}:
\textnormal{dom}(\eta^{(s^*,q^*)}) \longrightarrow \pi^{-1}(\textnormal{trace}(\gamma))
\subset \sphere^3$, defined on a non-empty, 
open and connected subset
$\textnormal{dom}(\eta^{(s^*,q^*)}) \subsetneq \sphere^1$, of $\gamma:\sphere^1 \longrightarrow \sphere^2$ with respect to the Hopf-fibration $\pi$, such that $\textnormal{dom}(\eta^{(s^*,q^*)})$ contains the point $s^*$ and such that $\eta^{(s^*,q^*)}$ attains the value
$q^*$ in $s^*$; i.e. $\eta^{(s^*,q^*)}$ is a smooth path in the torus $\pi^{-1}(\textnormal{trace}(\gamma))$, which intersects the fibers of $\pi$ perpendicularly and satisfies:
	\begin{displaymath} 
	(\pi \circ \eta^{(s^*,q^*)})(s) = \gamma(s)  \qquad
	\forall \,s \in \textnormal{dom}(\eta^{(s^*,q^*)})  \quad
	\textnormal{and} \quad \eta^{(s^*,q^*)}(s^*) = q^*,
	\end{displaymath}
and there is exactly one such function $\eta^{(s^*,q^*)}$
mapping the open and connected subset 
$\textnormal{dom}(\eta^{(s^*,q^*)}) \subsetneq \sphere^1$ into $\pi^{-1}(\textnormal{trace}(\gamma)) \subset \sphere^3$.
\item[2)] There is some $\epsilon=\epsilon(F,\gamma)>0$ such that for every fixed $s^* \in \sphere^1$ and every $x^*\in (\pi \circ F)^{-1}(\gamma(s^*))$
there is a horizontal smooth lift $\eta_F^{(s^*,x^*)}$ of $\gamma\lfloor_{\sphere^1 \cap B_{\epsilon}(s^*)}$ with respect to the fibration $\pi \circ F :\Sigma \longrightarrow \textnormal{trace}(\gamma) \subset \sphere^2$,
attaining the value $x^*$ in $s^*$, i.e. $\eta_F^{(s^*,x^*)}$ is a smooth path in the torus $\Sigma$ which intersects the fibers of $\pi \circ F$ perpendicularly and satisfies:
\begin{displaymath} 
(\pi \circ F \circ \eta_F^{(s^*,x^*)})(s) = \gamma(s)  \qquad \forall \,s \in \sphere^1 \cap B_{\epsilon}(s^*) \quad \textnormal{and} \quad \eta_F^{(s^*,x^*)}(s^*)=x^*.
\end{displaymath}	
This implies in particular that for the above $\epsilon=\epsilon(F,\gamma)>0$ the function 
$\eta_F \mapsto F \circ \eta_F$ maps
the set $\LL(\gamma\lfloor_{\sphere^1 \cap B_{\epsilon}(s^*)},\pi \circ F)$
of horizontal smooth lifts of $\gamma\lfloor_{\sphere^1 \cap B_{\epsilon}(s^*)}$
with respect to $\pi \circ F$ surjectively onto the set
$\LL(\gamma\lfloor_{\sphere^1 \cap B_{\epsilon}(s^*)},\pi)$ of horizontal smooth lifts of $\gamma\lfloor_{\sphere^1 \cap B_{\epsilon}(s^*)}$ with respect to $\pi$.
\end{itemize}
\end{lemma}
\emph{Proof:\,} See the proof of Lemma 3 in \cite{Ruben.MIWF.II}.

\section{Appendix C}
In the proof of the first part of Theorem \ref{singular.time.MIWF.Hopf.tori} we employed the 
following two GMT-results, which can be quickly derived 
from the general ``monotonicity formula'' for  
$n$-rectifiable varifolds in $\rel^{n+m}$ 
with locally bounded first variations; 
see Paragraph 17 in \cite{Simon.1984}.
\begin{proposition} \label{14.5} 
	Let $\nu_j$ be $n$-rectifiable varifolds defined 
	on an open subset $\Omega \subseteq \rel^{n+m}$,
	$n,m\in \nat$, with locally bounded first variations 
	$\delta \nu_j$ and such that for every open ball 
	$B_{\varrho}:=B^{n+m}_{\varrho}(x_0) \subseteq \Omega$
	there holds:
	$$ 
	\nu_j(\Omega) + \varrho^{1-\alpha n-\beta}\, 
	\nu_j(B_{\varrho})^{\alpha-1} \, 
	\parallel \delta \nu_j \parallel(B_{\varrho}) 
	\leq \Lambda \quad \forall j \in \nat,   
	$$  	
	where $\alpha, \beta$ can be any pair of positive 
	numbers and $\Lambda>>1$ arbitrarily large. 
	If moreover $\nu_j \longrightarrow \nu$ weakly as Radon measures on $\Omega$, then the $n$-dimensional Hausdorff-density $\theta^n(\nu)$ of $\nu$ exists 
	in every point of $\Omega$, and for any convergent 
	sequence $x_j \longrightarrow x_0 \in \Omega$ there holds:  
	$$ 
	\theta^n(\nu,x_0) \geq 
	\limsup_{j\to \infty} \theta^n(\nu_j,x_j).      
	$$
\end{proposition}
\begin{proposition} \label{14.7} 
Let $\nu_j$ be $n$-rectifiable varifolds defined 
on an open subset $\Omega \subseteq \rel^{n+m}$,
$n,m\in \nat$, with locally bounded first variations 
$\delta \nu_j$ and such that for every open ball 
$B_{\varrho}:=B^{n+m}_{\varrho}(x_0) \subseteq \Omega$
there holds:
$$ 
\nu_j(\Omega) + \varrho^{1-\alpha n-\beta}\, 
\nu_j(B_{\varrho})^{\alpha-1} \, 
\parallel \delta \nu_j \parallel(B_{\varrho}) 
\leq \Lambda \quad \forall j \in \nat,   
$$  	
where $\alpha, \beta$ can be any pair of positive 
numbers and $\Lambda>>1$ arbitrarily large. 
If moreover the $n$-dimensional Hausdorff-densities 
$\theta^n(\nu_j)$ of $\nu_j$ satisfy 
$\theta^n(\nu_j) \geq 1$ on 
$\textnormal{spt}(\nu_j)$, for each $j\in \nat$, and if 
$\nu_j \longrightarrow \nu$ weakly as Radon measures on $\Omega$, then there holds: 
$$ 
\textnormal{spt}(\nu_j) \longrightarrow  
\textnormal{spt}(\nu)  \quad 
\textnormal{in Hausdorff distance locally in} \,\, \Omega,    
$$ 
as $j\to \infty$, and even more precisely: 
$$ 
\textnormal{spt}(\nu)=\{\,x\in \Omega \,|\, 
\exists\, x_j \in \textnormal{spt}(\nu_j) \,\, 
\textnormal{for every} \,j\in \nat \,\, 
\textnormal{such that} \,\, x_j \longrightarrow x\,\}.     
$$ 
\end{proposition}
\noindent      \\
{\small \underline{Acknowledgements}\\\\
	The M\"obius-invariant Willmore flow was originally discovered by Professor Ben Andrews. 
	I would like to thank Professor Itai Shafrir and Professor Yehuda Pinchover for their support and hospitality at the Israel Institute of Technology. I was funded by the Ministry of Absorption of the State of Israel throughout the years 2019--2022.}


\begin{thebibliography}{3}
	\bibitem{Alt.2015} Alt, H.W.: Linear Functional Analysis. 
	Universitext, Springer publisher, Berlin-Heidelberg, 2015.
	\bibitem{Andrews.Hopper.2011} Andrews, B., Hopper, C.: 
	The Ricci flow in Riemannian Geometry.
	A complete proof of the differentiable
	1/4-pinching Sphere Theorem. 
	Springer publisher, Berlin-Heidelberg, 2011. 
	\bibitem{Arnold.1993} Arnold, V.: 
	Plane Curves, their Invariants, Perestroikas and 
	Classifications. Lecture notes, FIM-ETH Z\"urich, 
	1993.
	\bibitem{Bellettini.1993} Bellettini, G., Dal Maso, G., 
	Paolini, G.: 
	Semicontinuity and relaxation properties of a curvature
	depending functional in 2D. Annali della Scuola Normale Superiore di Pisa, Classe di Scienze $4^e$ s\'erie {\bf 20}, 247--297 (1993).  
	\bibitem{Bernard.2016} Bernard, Y.: 
	Noether's theorem and the Willmore functional. 
	Adv. Calc. Var. {\bf 9} (3), 217--234 (2016).
	\bibitem{Blatt.2009} Blatt, S.: A singular example for 
	the Willmore flow. Analysis {\bf 29}, 407--430 (2009).
	\bibitem{Brezis.Coron.1984} Brezis, H., 
	Coron, J.-M.: Multiple solutions of $H$-Systems and 
	Rellich's Conjecture. Comm. on Pure and Applied Mathematics {\bf 37}, 149--187 (1984). 
	\bibitem{Cairns.Sharpe.Webb.1994} 
	Cairns, G., Sharpe, R.W., Webb, L.: 
	Conformal Invariants for curves and surfaces in three dimensional space forms. 
	Rocky Mountain J. of Math. {\bf 24}, No. 3, 933--959 (1994).  
	\bibitem{Dall.Acqua.Pozzi.2014} Dall'Acqua, A., Pozzi, P.: A Willmore-Helfrich $L^2$-flow of curves with natural
	boundary conditions. Comm. Anal. Geom. {\bf 22}, no. 4, 617--669.
	\bibitem{Dall.Acqua.Pozzi.2018} Dall'Acqua, A., Laux, T., 
	Chun-Chi Lin, Pozzi, P., et. al.: The elastic flow of curves on the sphere. Geometric flows, Vol. {\bf 3}, 1--13 (2018).
	\bibitem{Dall.Acqua.Schaetzle.Mueller.2024} 
	Dall'Acqua, A., M\"uller, M., Sch\"atzle, R., Spener, A.: 
	The Willmore flow of tori of revolution. Anal. PDE 
	{\bf 17}, 3079--3124 (2024). 
	\bibitem{De.Lellis.Mueller.2005} De Lellis, C., M\"uller, S.: Optimal rigidity estimates for nearly umbilical surfaces. J. Differ. Geom. {\bf 69}, no. 1, 75--110 (2005).
	\bibitem{De.Lellis.Mueller.2006} De Lellis, C., M\"uller, S.: A $C^0$-estimate for nearly umbilical surfaces. 
	Calc. Var. {\bf 26}, no. 3, 283--296 (2006).	
 	\bibitem{Giaquinta.1984} Giaquinta, M.: 
 	Multiple Integrals in the Calculus of Variations and Nonlinear Elliptic Systems. Annals of Mathematics Studies {\bf 105}, Princeton University Press, Princeton, New Jersey, 1983. 
	\bibitem{Helein.2004} Helein, F.: Harmonic maps, conservation laws and moving frames. Second edition. Cambridge Tracts in Mathematics {\bf 150}, Cambridge University Press, 2004.
	\bibitem{Huisken.1984} Huisken, G.: Flow by mean curvature 
	of convex surfaces into spheres. 
	J. Differential Geometry {\bf 20}, 237--266 (1984). 
	\bibitem{Jakob_Moebius_2016} Jakob, R.: Short-time existence of the M\"obius-invariant Willmore flow. J. Geom. Anal. {\bf 28}(2), 1151--1181 (2018).
	\bibitem{Ruben.MIWF.III} Jakob, R.: Functional analytic properties and regularity of the M\"obius-invariant Willmore flow in $\mathbf{R}^n$. Partial Differential Equations and Applications {\bf 3}, Article No. 67 (2022), arXiv:2011.03832.
	\bibitem{Ruben.MIWF.II} Jakob, R.: 
	The Willmore flow of Hopf-tori in the $3$-sphere. J. Evol. Equ. {\bf 23}, Article no. 72 (2023), arXiv:2002.01006. 
	\bibitem{Ruben.MIWF.IV} Jakob, R.: Global existence and full convergence of the M\"obius-invariant Willmore flow 
	in the $3$-sphere. J. Geom. Anal. {\bf 34}, Article no. 24 (2024), arXiv:2101.00471.
	\bibitem{Jost.2006} Jost, J.: Compact Riemann surfaces. 
	An Introduction to Contemporary Mathematics. 3rd edition. Universitext. Springer publisher, Berlin-Heidelberg, 2006.    
	\bibitem{Jost.2017} Jost, J.: Riemannian Geometry and 
	Geometric Analysis. 7th edition. Universitext, 
	Springer Publisher. Berlin-Heidelberg, 2017.  
	\bibitem{Keller.Modino.Riviere} 
	Keller, L.G.A., Mondino, A., Rivi\`ere. T.: 
	Embedded Surfaces of Arbitrary Genus
	Minimizing the Willmore Energy Under
	Isoperimetric Constraint. 
	Arch. Rat. Mech. Anal. {\bf 212}, 645--682 (2012).
	\bibitem{Kuwert.Li.2012} Kuwert, E., Li, Y.: 
	$W^{2,2}$-conformal immersions of a closed
	Riemann surface into $\rel^n$. 
	Comm. Anal. Geom. {\bf 20}, 313--340 (2012).   
	\bibitem{Kuwert_Schaetzle_2001} Kuwert, E., Sch\"atzle, M.R.: The Willmore flow with small initial energy, J. Differential Geometry {\bf 57}, 409--441 (2001).
	\bibitem{Kuwert_Schaetzle_2002} Kuwert, E., Sch\"atzle, M.R.: Gradient flow for the Willmore functional. Comm. Anal. Geom. {\bf 10}, 307--339 (2002).
	\bibitem{Kuwert.Schaetzle.Annals} Kuwert, E., Sch\"atzle, M.R.: Removability of point singularities of Willmore surfaces. Ann. of Math. (2) {\bf 160}, 315--357 (2004).
	\bibitem{Kuwert.Schaetzle.2012} Kuwert, E., Sch\"atzle, 
	M.R.: Closed surfaces with bounds on their Willmore energy. Ann. Sc. Norm. Super. Pisa Cl. Sci. (5), {\bf 11}, 605--634 (2012). 
	\bibitem{Kuwert.Schaetzle.conf.class.2013}
	Kuwert, E., Sch\"atzle, M.R.: 
	Minimizers of the Willmore functional 
	under fixed conformal class. J. Diff. Geom.
	{\bf 93}, 471--530 (2013). 
	\bibitem{Kuwert.Scheuer.2020} 
	Kuwert, E., Scheuer, J.: Asymptotic Estimates for the Willmore Flow With Small Energy.
	International Mathematics Research Notices {\bf 2021}, 
	No. 18, 14252--14266 (2021).
	\bibitem{Langer.Singer.1984} Langer J., Singer, D.A.: The total squared curvature of closed curves. J. Differential Geometry {\bf 20}, 1--22 (1984).
	\bibitem{Langer.Singer.1987} Langer J., Singer, D.A.: Curve-straightening in Riemannian manifolds. Ann. Global Anal. Geom. {\bf 5}, No. 2, 133--150 (1987).
	\bibitem{Marques.Neves.2014} Marques, F.C., Neves, A.:
	Min-Max theory and the Willmore conjecture. Annals of Mathematics {\bf 179}, 683--782 (2014).
	\bibitem{Mayer} Mayer, M.: The inverse Willmore flow. Diploma thesis, Universit\"at T\"ubingen, (2009) arXiv:1508.07800.
	\bibitem{Mueller.Sverak.1995} M\"uller, S., Sver\'ak, V.: 
	On Surfaces of finite total curvature. 
	J. Diff. Geom. {\bf 42}, No.2, 229--258 (1995). 
	\bibitem{Ndiaye.Schaetzle.2014} Ndiaye, C.B., Sch\"atzle, R.M.: Explicit conformally constrained Willmore minimizers in arbitrary codimension, Calc. Var. {\bf 51}, 291--314 (2014).
	\bibitem{Palmurella.2022} Palmurella, F., Rivi\`ere. T.: 
	The parametric approach	to the Willmore flow. 
	Advances in Mathematics {\bf 400}, Paper No. 108257, 
	48 pages (2022). 
	\bibitem{Palmurella.2024} Palmurella, F., Rivi\`ere. T.: 
	The parametric Willmore flow. J. Reine Angew. Math. 
	{\bf 811}, 1--91, (2024).
	\bibitem{Pinkall} Pinkall, U.: Hopf-tori in $\sphere^3$. 
	Invent. math. {\bf 81}, 379--386 (1985).
    \bibitem{Polyak.1996} Polyak, M.: Invariants of 
    curves and fronts via Gauss diagrams. 
    Topology {\bf 37}, No. 5, 989--1009 (1998).
	\bibitem{Pruess.Simonett.2013} Pr\"uss, J., Simonett G.:
	On the manifold of closed hypersurfaces in $\mathbf{R}^n$. 
	Discrete Contin. Dyn. Syst. {\bf 33}, 5407--5428 (2013).
	\bibitem{Riviere.2008} Rivi\`ere, T.: 
	Analysis aspects of Willmore surfaces. 
	Invent. math. {\bf 174}, 1--45 (2008).
	\bibitem{Riviere.2011} Rivi\`ere, T.: 
	Conformally invariant $2$-dimensional variational problems. Cours joint de l'Institut Henri Poincar\'e, Paris XII Creteil, Novembre 2010.
	\bibitem{Riviere.Lip.conf.imm.2013} 
	Rivi\`ere, T.: Lipschitz conformal immersions from
	degenerating Riemann surfaces with $L^2$-bounded 
	second fundamental forms. Adv. Calc. Var. {\bf 6}, 
	1--31 (2013). 
	\bibitem{Riviere.Var.principle.2014}
	Rivi\`ere, T.: Variational principles for immersed surfaces with $L^2$-bounded second fundamental form. 
	J. reine angew. Math. {\bf 695}, 41--98 (2014).
	\bibitem{Riviere.Park.City.2013}
	Rivi\`ere, T.: Weak immersions of surfaces with
	$L^2$-bounded second fundamental form, Geometric
	Analysis, IAS/Park City Math. Ser., {\bf 22}, 
	Amer. Math. Soc., Providence, RI, 303--384 (2016), 
	MR 3524220. 
	\bibitem{Schaetzle.Lower.semicont.2009} 
	Sch\"atzle, M.R.: Lower semicontinuity of the 
	Willmore functional for currents.  
	J. Diff. Geom. {\bf 81}, 437--456 (2009). 
	\bibitem{Schaetzle.Conf.factor.2013}
	Sch\"atzle, M.R.: Estimation of the conformal factor under bounded Willmore energy. Math. Z., {\bf 274}, 1341--1383 (2013).
	\bibitem{Simon.1984} Simon, L.: Lectures on Geometric 
	Measure Theory. Proceedings of the Centre for Mathematical 
	Analysis, Australian National University, 1984. 
	\bibitem{Simon.1993} Simon, L.: Existence of surfaces minimizing the Willmore functional. CAG {\bf 1}, no. 2, 281--326 (1993).
	\bibitem{Skorzinski.2015} Skorzinski, F.: Local minimizers of the Willmore functional. Analysis {\bf 35}, 93--115 (2015).
	\bibitem{Tartar.1998} Tartar, L.: Imbedding theorems of Sobolev spaces into Lorentz spaces. Bollettino dell’Unione Matematica Italiana, Serie 8, {\bf 1-B}, 479--500 (1998).
	\bibitem{Teufel} Teufel, E.: The isoperimetric inequality and the total absolute curvature of closed curves in spheres, Manuscripta Math., {\bf 75} (1), 43--48 (1992).
	\bibitem{Tromba.Teichmueller.1992} Tromba, A.J.: 
	Teichm\"uller Theory in Riemannian Geometry. 
	2nd Edition, Springer Publisher, Basel AG, 1992.  
 	\bibitem{Weiner} Weiner, J. L.: On a problem of Chen, Willmore et al. Indiana University Mathematics Journal 
 	{\bf 27}, No.1, 19--35 (1978).
 	\bibitem{Wente.1969} Wente, H.: An Existence Theorem for 
 	Surfaces of Constant Mean Curvature. 
 	J. Math. Anal. Appl. {\bf 26}, 318--344 (1969). 
 \end{thebibliography}
\end{document}